\documentclass[11pt,reqno]{amsart}
\usepackage[T1]{fontenc}
\usepackage[utf8]{inputenc}
\usepackage{lmodern}
\usepackage[letterpaper,margin=1.05in]{geometry}
\usepackage{amsmath,amssymb,amsfonts,mathtools,mathrsfs}
\usepackage{microtype}
\usepackage{enumitem,booktabs,longtable,array}
\usepackage[dvipsnames]{xcolor}
\usepackage{hyperref}
\hypersetup{colorlinks=true,linkcolor=MidnightBlue,citecolor=MidnightBlue,urlcolor=MidnightBlue,
  pdftitle={Projective division and the homogeneous Arveson--Douglas conjecture},
  pdfsubject={Research manuscript on the homogeneous Arveson--Douglas conjecture}}
\numberwithin{equation}{section}
\newtheorem{theorem}{Theorem}[section]
\newtheorem{lemma}[theorem]{Lemma}
\newtheorem{proposition}[theorem]{Proposition}
\newtheorem{corollary}[theorem]{Corollary}
\newtheorem{thmA}{Theorem}
\newtheorem{thmB}{Theorem}
\newtheorem{thmC}{Theorem}
\newtheorem{thmD}{Theorem}
\newtheorem{thmE}{Theorem}
\newtheorem*{introB}{Theorem B}

\newtheorem*{introD}{Theorem D}
\newtheorem*{introE}{Theorem E}

\theoremstyle{definition}

\newtheorem{example}[theorem]{Example}
\theoremstyle{remark}
\newtheorem{remark}[theorem]{Remark}
\newcommand{\C}{\mathbb C}

\newcommand{\N}{\mathbb N}

\newcommand{\Pj}{\mathbb P}
\newcommand{\cO}{\mathcal O}
\newcommand{\cI}{\mathcal I}

\newcommand{\cS}{\mathcal S}
\newcommand{\FS}{\mathrm{FS}}
\newcommand{\DA}{\mathrm{DA}}
\newcommand{\dbar}{\bar\partial}
\newcommand{\eps}{\varepsilon}
\newcommand{\id}{\mathrm{id}}

\DeclareMathOperator{\im}{im}

\DeclareMathOperator{\supp}{supp}
\DeclareMathOperator{\Supp}{Supp}
\DeclareMathOperator{\Tr}{Tr}

\DeclareMathOperator{\Gr}{Gr}
\DeclareMathOperator{\Quot}{Quot}

\DeclareMathOperator{\dist}{dist}

\DeclareMathOperator{\Proj}{Proj}
\DeclareMathOperator{\Spec}{Spec}

\allowdisplaybreaks[2]
\setlist[enumerate]{label=\textup{(\roman*)},leftmargin=2.4em}
\setlist[itemize]{leftmargin=2em}
\title[The homogeneous Arveson--Douglas theorem]{The homogeneous Arveson--Douglas theorem}
\date{26 November 2026 (expanded exposition)}

\author{Lijia Ding}
\address{School of Mathematics and Statistics, Zhengzhou University, 450001, Zhengzhou, P.R.China}
\email{ljding@zzu.edu.cn}
\author{Xiang Fang}
\address{Department of Applied Mathematics, National Yang Ming Chiao Tung University, Hsinchu, Taiwan}
\email{xfang@nycu.edu.tw}
\author{Xiang Tang}
\address{Department of Mathematics, Washington University, St. Louis, MO, 63130, U.S.A.}
\email{xtang@math.wustl.edu}

\begin{document}
\begin{abstract}
We prove the homogeneous Arveson--Douglas conjecture for homogeneous
submodules of the finite-multiplicity Drury--Arveson shift with the
standard grading. For each fixed submodule, we establish a bound on the positive
defect square associated with quotient compression that decays
inversely with the degree. When the corresponding algebraic quotient
has positive support dimension, we deduce that its coordinate
commutators belong to the weak Schatten class whose exponent equals
that dimension.

We organize the proof in three steps. First, we realize homogeneous
polynomials as projective sections through an exact norm identity
and construct complete presentations of the corresponding ideal
sheaves that remain stable under normalized limits. We then combine
quantitative division, sampled Hermitian metrics, and Dolbeault
estimates to obtain holomorphic lifting estimates uniform in degree.
Finally, we transfer these estimates back to the original
Drury--Arveson norm, use an auxiliary-variable construction to pass
to vector coefficients, and count singular values in the original
quotient.
\end{abstract}
\maketitle
\tableofcontents

\section{Introduction}
\label{sec:introduction}

\subsection{The conjecture and the main theorem}
The homogeneous Arveson--Douglas conjecture concerns the effect of
polynomial relations on the multivariable shift. It predicts that the
coordinate operators on a homogeneous quotient of the Drury--Arveson
space are essentially normal, meaning that their commutators with one
another's adjoints are compact. Its quantitative form asks for Schatten
summability governed by the dimension of the algebraic support
\cite{Arveson2005,DouglasIndex}. In this paper, we prove the standard homogeneous
conjecture for finite coefficient multiplicity by establishing a
degreewise estimate for the part of each coordinate operator removed
by the quotient projection.

Let $m,r\ge1$, let $R=\C[z_1,\ldots,z_m]$, and let $H_m^2$ denote the
Drury--Arveson space on the unit ball $\mathbb B^m$, with monomial norms
\[
 \langle z^\alpha,z^\beta\rangle=0\quad(\alpha\ne\beta),
 \qquad
 \|z^\alpha\|_{\DA}^2=\frac{\alpha!}{|\alpha|!}.
\]
The coefficient space $\C^r$ has its standard Hermitian metric and is
placed in degree zero. For a homogeneous submodule
$M\subset R\otimes\C^r$, set
\[
 \mathcal M=\overline M\subset H_m^2\otimes\C^r,\qquad
 Q=\mathcal M^\perp,\qquad
 Q_n=(H_{m,n}^2\otimes\C^r)\ominus M_n.
\]
We write
\[
 Z_j=M_{z_j}\otimes I_r,\qquad
 S_j=P_QZ_j|_Q,\qquad V_j=Z_j|_{\mathcal M}.
\]
We use the algebraic support dimension of the finite graded module
$T=(R\otimes\C^r)/M$. More explicitly, let
$\operatorname{Ann}_R T=\{f\in R:fT=0\}$. The support is the closed
subset $V(\operatorname{Ann}_R T)$ of affine space (or its corresponding
closed subset of $\Spec R$), and its dimension is the Krull dimension
of $R/\operatorname{Ann}_R T$. Thus for scalar ideals it is the
complex dimension of the usual zero set, not the number of generators
or the dimension of the degree-zero space. We recall the algebraic
terminology and Hilbert-function theorem in Section~\ref{sec:operators}.

If $T=(R\otimes\C^r)/M$ is infinite dimensional and
$\delta=\dim\Supp T$, the dimension-sensitive assertion is
\[
 [S_i^*,S_j]\in\cS_p\qquad(p>\delta).
\]
Here $\cS_p$ is the Schatten class. We also use the weak class
$\cS_{p,\infty}$, consisting of compact operators whose singular values
satisfy $s_k=O(k^{-1/p})$.

Homogeneity separates this problem into an estimate on each finite
degree space and a count of its dimension. We define
\[
 B_{j,n}=P_{M_{n+1}}Z_j|_{Q_n}
\]
and the positive defect square
\begin{equation}\label{eq:intro-defect}
 A_n=\sum_{j=1}^m B_{j,n}^*B_{j,n}.
\end{equation}
We can now state our main result.

\begin{thmA}\label{thm:main}
For every homogeneous submodule $M\subset R\otimes\C^r$ there is a
constant $C_M$ such that
\begin{equation}\label{eq:main-degree}
 0\le A_n\le\frac{C_M}{n+1}I_{Q_n}\qquad(n\ge0).
\end{equation}
Consequently
\begin{equation}\label{eq:main-commutator}
 \|[S_i^*,S_j]|_{Q_n}\|\le\frac{C_M'}{n+1}
 \qquad(1\le i,j\le m).
\end{equation}
If $T=(R\otimes\C^r)/M$ is infinite dimensional and
$\delta=\dim\Supp T$, then
\[
 [S_i^*,S_j]\in\cS_{\delta,\infty}\subset\bigcap_{p>\delta}\cS_p.
\]
The column $B=(P_{\mathcal M}Z_j|_Q)_{j=1}^m$ belongs to
$\cS_{2\delta,\infty}$, as do the individual off-diagonal operators
$[P_{\mathcal M},Z_j]$ after extension by zero.  The restriction
commutators satisfy
\begin{equation}\label{eq:schatten=restriction}
 [V_i^*,V_j]\in\cS_{m,\infty}\subset\bigcap_{p>m}\cS_p.
\end{equation}
If $T$ is finite dimensional, the quotient commutators and off-diagonal
operators have finite rank; estimate (\ref{eq:schatten=restriction}) still holds.
\end{thmA}

The constant $C_M$ in Theorem \ref{thm:main} depends on the fixed module
and its original Hermitian geometry. In Lemma \ref{op:graded}, we show
that the polynomial formulation of Theorem \ref{thm:main} is equivalent
to its formulation for closed submodules invariant under
the standard circle action. Indeed, the homogeneous polynomial
part of every such submodule is finitely generated and dense.
For a scalar ideal $I$, the dimension parameter $\delta$
is the complex dimension of its zero set $V(I)$.
Throughout, we use the standard total-degree grading.

Our proof separates the analytic decay of the homogeneous commutator
blocks from the algebraic growth of their dimensions.
The estimate for $A_n$ is the main analytic input in
Theorem~\ref{thm:main}. Since $B_{j,n}^*B_{j,n}\le A_n$ for every $j$,
this estimate implies
$\|B_{i,n}^*B_{j,n}\|=O((n+1)^{-1})$.
The exact compression identity
\[
 [S_i^*,S_j]|_{Q_n}
 =P_{Q_n}[Z_i^*,Z_j]|_{Q_n}-B_{i,n}^*B_{j,n},
\]
together with the corresponding decay estimate for the ambient
commutator, therefore gives
\[
 \|[S_i^*,S_j]|_{Q_n}\|=O((n+1)^{-1}).
\]

To pass from this degreewise estimate to a singular-value bound,
we use the Hilbert polynomial theorem.
When $\delta>0$, the Hilbert function of
$T=(R\otimes\mathbb C^r)/M$ agrees, in all sufficiently large degrees,
with a polynomial of degree $\delta-1$.
Since $\dim Q_n=\dim_{\mathbb C}T_n$, it follows that
$\dim Q_n=O(n^{\delta-1})$, and hence
$\sum_{\ell=0}^{n}\dim Q_\ell=O(n^\delta)$.
Combining this dimension count with the degreewise norm estimate,
we obtain
$s_k([S_i^*,S_j])=O(k^{-1/\delta})$,
which gives $[S_i^*,S_j]\in\cS_{\delta,\infty}$.
Example \ref{ext:sharp} in Section~\ref{sec:extensions} shows that,
in general, this weak Schatten conclusion cannot be strengthened
to membership in $\cS_\delta$.

\subsection{Historical developments}
\label{sec:history}

We refer to the survey of Guo and Yi Wang~\cite{GuoYiWang}
for a comprehensive account of the background, applications,
and developments of the Arveson--Douglas conjecture up to 2020. Arveson  \cite{Arveson1998} established the multivariable shift's Schatten estimates and
Toeplitz extension, and \cite{Arveson2005}  proved essential normality
for homogeneous monomial submodules with finite vector multiplicity. In \cite{DouglasIndex}, Douglas connected the problem to $K$-homology and
proposed the dimension-sensitive summability threshold. Guo and K. Wang \cite{GuoWang} proved the homogeneous principal
case in arbitrary dimension and the scalar homogeneous conjecture
in two and three variables.

Shalit's stable division method \cite{ShalitDivision} and Eschmeier's
operator-inequality criterion \cite{EschmeierCriterion} identify
projected derivatives as the estimates that control the quotient
commutators. Y. Wang's local decomposition method extends this division
viewpoint in the Bergman setting \cite{YiWang}. Geometrically,
Engli\v{s} and Eschmeier \cite{EnglisEschmeier} proved the radical-ideal case for homogeneous
varieties smooth away from the origin, using generalized Toeplitz
operators. Their result includes the Drury--Arveson space and
explicitly gives the weak estimate $\cS_{\delta,\infty}$.

Douglas, Tang, and Yu~\cite{DouglasTangYu} developed Dolbeault
operator analysis and a bounded extension theorem,
extending Beatrous's work \cite{Beatrous}, for complete intersections
with isolated interior singularities and smooth
transverse intersections with the boundary of the ball. They obtained essential normality and identified the quotient Toeplitz
class with the boundary fundamental class in the ball Bergman
setting. Douglas, Jabbari, Tang, and Yu \cite{DouglasJabbariTangYu} used
finite resolutions by essentially normal Hilbert modules to compute
the $K$-homology classes of monomial quotients; their method also
applies to the Drury--Arveson space. Ding \cite{Ding26} developed related results for holomorphic Sobolev quotients on strongly pseudoconvex manifolds under geometric hypotheses, using embedded desingularization, and studied Schatten estimates for antisymmetrized products of compressed multiplication operators. 

\subsection{The proof in three steps}
In Sections~\ref{sec:projective}--\ref{sec:completion}, we prove Theorem \ref{thm:main} in three steps: 
\begin{enumerate}
\item transfer in Theorem \ref{thm:B} the scalar homogeneous problem to coherent ideal sheaves on projective space;
\item construct in Theorem \ref{thm:C}-\ref{geo:lifting} a bounded holomorphic lift by local division and $\dbar$ estimates on smooth presenting bundles; 
\item return to the original Drury-Arveson norm and pass to vector coefficients before counting singular values. 
\end{enumerate}
We use the following notation. Set $N=m-1$ and write
$\mathbb P^N$ for the space of complex lines in
$\mathbb C^{N+1}$, $\mathcal O$ for its structure sheaf
of holomorphic functions, and $\mathcal O(1)$ for the
dual of the tautological line bundle. A coherent ideal sheaf records local polynomial or
holomorphic ideal membership, not just the common zero set. Use
$\mathscr F(t)$ for $\mathscr F\otimes\cO(t)$, $H^0$ for global
sections, and $H^i$ for sheaf cohomology. The core of the argument is
an estimate on smooth bundles over the ambient projective space;
we do not need differential forms on a singular or nonreduced
subscheme itself. 

\subsubsection{Step 1: From the graded Hilbert module to projective sheaves}
We begin by describing the passage from homogeneous polynomials to
projective sections. Let $R=\mathbb C[z_1,\ldots,z_m]$, initially with
$m\ge2$, and let $R_n$ be the vector space of homogeneous polynomials
of total degree $n$. We write
$\mathbb P^{m-1}(\mathbb C)$ for the space of complex lines in
$\mathbb C^m$: a nonzero vector $z$ represents the point $[z]$,
and $[z]=[\lambda z]$ for every nonzero scalar $\lambda$. Let \(U_j\) be  the \(j\)-th standard coordinate chart of projective space:
\[
U_j=\bigl\{[z_1:\cdots:z_m]\in\mathbb P^{m-1}(\mathbb C):z_j\ne0\bigr\},
\qquad j=1,\ldots,m.
\]
It is easy to see that 
\[
\mathbb{P}^{m-1}=U_1\cup \cdots U_m
\]
Let $\mathcal O$ denote the sheaf of holomorphic functions on $\mathbb P^{m-1}(\mathbb C)$,
which assigns to each open set its ring of holomorphic functions.
We use $L=\mathcal O(1)$ for the hyperplane line bundle, whose
fiber over a line $\ell\subset\mathbb C^m$ is the dual vector space
$\ell^*$. We identify a holomorphic line bundle with its sheaf of
holomorphic sections and write
$L^n=L^{\otimes n}=\mathcal O(n)$ for its $n$-th tensor power.
Here $L^0=\mathcal O$.
The notation $H^0(\mathbb P^{m-1}(\mathbb C),L^n)$ means the vector space of global
holomorphic sections of $L^n$. For $n\ge0$, we have a natural identification
\[
 R_n \cong H^0(\mathbb P^{m-1}(\mathbb C),L^n).
\]
Concretely, on the coordinate chart
$U_j=\{[z]\in \mathbb P^{m-1}(\mathbb C):z_j\ne0\}$, the section corresponding to
$F_n\in R_n$ has local representative $F_n(z)/z_j^n$.
We use the same symbol $F_n$ for the polynomial and its
corresponding section.

We equip these spaces with two norms.
For
$F_n(z)=\sum_{|\alpha|=n}a_\alpha z^\alpha$,
the Drury--Arveson norm is
\[
 \|F_n\|_{\DA}^2
 =\sum_{|\alpha|=n}\frac{\alpha!}{n!}|a_\alpha|^2.
\]
Here $\alpha=(\alpha_1,\ldots,\alpha_m)$ is a tuple of
nonnegative integers,
$|\alpha|=\alpha_1+\cdots+\alpha_m$,
$\alpha!=\alpha_1!\cdots\alpha_m!$, and
$z^\alpha=z_1^{\alpha_1}\cdots z_m^{\alpha_m}$.
The projective Fubini--Study norm is
\[
 \|F_n\|_{\FS}^2
 =\int_{\mathbb{P}^N(\mathbb{C})}\frac{|F_n(z)|^2}{\|z\|^{2n}}\,
       dV_{\FS}([z]),
 \qquad
 \|z\|^2=\sum_{j=1}^m|z_j|^2.
\]
The quotient in the integrand is independent of the nonzero
representative $z$ of $[z]$. The measure $dV_{\FS}$ is the
Fubini--Study volume normalized to have total mass one;
equivalently, it is the image of normalized surface measure
on the unit sphere under $z\mapsto[z]$.
With these conventions, the polynomial--section correspondence
satisfies the exact norm identity
\begin{equation}\label{eq:intro-projective-norm}
 \|F_n\|_{\DA}^2=c_n\|F_n\|_{\FS}^2,\qquad
 c_n=\binom{n+m-1}{m-1},\qquad
 \frac{c_n}{c_{n+1}}=\frac{n+1}{n+m}.
\end{equation}
Thus, for each fixed $n$, the two inner products differ by a
single positive scalar. They consequently define the same
orthogonal complements and the same orthogonal projections
onto any fixed subspace of $R_n$.

We next describe how an ideal enters this correspondence.
Let $I\subset R$ be a homogeneous ideal, meaning an ideal
generated by homogeneous polynomials, and put $I_n=I\cap R_n$.
We write $\mathscr I=\widetilde I\subset\mathcal O$ for its
associated ideal sheaf, viewed holomorphically. 
On the $j$-th standard coordinate $U_j$, this sheaf is locally generated by the functions
$G(z)/z_j^{\deg G}$, where $G$ ranges over homogeneous
generators of $I$.
Its twist
\[
 \mathscr I(n):=\mathscr I\otimes_{\mathcal O}L^n
\]
is the subsheaf of sections of $L^n$ whose local coefficients
belong to $\mathscr I$.
Accordingly, the $0$-th sheaf cohomology $H^0(\mathbb P^{m-1}(\mathbb C),\mathscr I(n))$ is the space of
global holomorphic sections satisfying these local ideal
membership conditions.

Let $\mathfrak m=(z_1,\ldots,z_m)$ be the ideal of polynomials
vanishing at the origin. We define the saturation of $I$ by
\[
 I^{\mathrm{sat}}
 :=\bigl\{F\in R:
       \mathfrak m^kF\subset I
       \text{ for some integer }k\ge0\bigr\}.
\]
Here $\mathfrak m^kF\subset I$ means that $GF\in I$ for every
$G\in\mathfrak m^k$.
Under the polynomial--section identification, we have
\[
 H^0(\mathbb{P}^{m-1}(\mathbb{C}),\mathscr I(n))=(I^{\mathrm{sat}})_n
 \qquad(n\ge0),
\]
where $(I^{\mathrm{sat}})_n=I^{\mathrm{sat}}\cap R_n$.
Moreover, there is an integer $n_0$, depending on $I$, such that
\[
 (I^{\mathrm{sat}})_n=I_n
 \qquad(n\ge n_0).
\]
Thus $H^0(\mathbb{P}^{m-1}(\mathbb{C}),\mathscr I(n))=I_n$ in all sufficiently large
degrees, which is the meaning of the notation $n\gg1$.

The discrepancy $I^{\mathrm{sat}}/I$ is finite-dimensional
over $\mathbb C$ and is annihilated by a power of $\mathfrak m$.
In particular, it has finite length: it admits a finite
filtration by $R$-submodules whose nonzero successive quotients
are isomorphic to $R/\mathfrak m=\mathbb C$.
It therefore contributes in only finitely many homogeneous
degrees and is supported at the origin.
For a proper homogeneous ideal, the common zero set of its
polynomials is an affine cone, meaning that it is invariant
under scalar multiplication; the origin is its vertex.
Replacing $I$ by $I^{\mathrm{sat}}$ removes precisely the
submodule $I^{\mathrm{sat}}/I$ of $R/I$, which records this
finite-length information at the vertex.

At the projective level, however, the ideal sheaf is unchanged:
\[
 \widetilde{I^{\mathrm{sat}}}=\widetilde I=\mathscr I.
\]
Consequently, the quotient sheaf $\mathcal O/\mathscr I$,
which records local holomorphic functions modulo the defining
equations, is also unchanged.
This preserves all projective nonreduced and embedded structure.
Here nonreduced structure means the presence of a nonzero
local class $u$ such that $u^k=0$ for some positive integer $k$.
Embedded structure is recorded by nonminimal associated primes
of the local quotient rings: an associated prime of such a
ring $A$ is a prime ideal of the form
$\operatorname{Ann}_A(u)=\{a\in A:au=0\}$ for some nonzero
$u\in A$, and it is embedded if it strictly contains a
minimal prime of $A$.
A minimal prime is a prime ideal minimal under inclusion. We prove these facts in Section~\ref{sec:projective}.
There we also treat separately the zero ideal sheaf
$\mathscr I=0$, the unit ideal sheaf
$\mathscr I=\mathcal O$, and the one-variable case $m=1$,
for which $\mathbb P^0(\mathbb C)$ is a single point.

In constructing the projective presentation, we work first
in the category of coherent algebraic sheaves on
$\mathbb P^{m-1}$. We then pass to the corresponding
holomorphic sheaves by analytification. By the
algebraic-analytic comparison theorem for projective
varieties (GAGA), coherent algebraic sheaves and their
analytifications have the same cohomology groups.
Moreover, analytification preserves exact sequences.
Unless the distinction is necessary, we use the same
notation for an algebraic sheaf and its analytification.

With a projective presentation of $\mathscr I$, we study its
local matrices after centering at a point $p\in \mathbb{P}^{m-1}(\mathbb{C})$, rescaling
coordinates by $\varepsilon>0$, and choosing coefficient bases
adapted to the Gram form (and its local matrix expression).
We need the resulting normalized presentations, including
their coefficientwise limits, to retain the generation and
exactness properties required for division.
These algebraic properties allow us to obtain local lifting
estimates with constants independent of the center $p$ and
the scale $\varepsilon$.  We therefore use complete spaces of
sections in degrees chosen uniformly over proper parameter families.
For a nonzero proper $\mathscr I$, write
\[
 E_{-1}=\cO,\qquad D_{-1}=0,\qquad
 E_q=V_q\otimes L^{-D_q},
\]
and $d_0=\Phi$ as the composition of the evaluation map $
V_0\otimes L^{-D_0}\longrightarrow\mathscr I$
with the inclusion \(\mathscr I\hookrightarrow\mathcal O\). Next, set
\[
 \mathscr K_0=\mathscr I,\qquad \mathscr K_q=\ker d_{q-1}\quad(q\ge1),
\]
and $P_q(t)=\chi(\mathscr K_q(t))$, the holomorphic Euler characteristic of the sheaf $\mathscr{K}_q(t)$. 
A \emph{complete presentation} of $\mathscr{I}$ includes 
$V_q=H^0(\mathbb{P}^{m-1},\mathscr K_q(D_q))$ and the evaluation maps
\begin{equation}\label{eq:dq}
d_q:E_q\to E_{q-1}.
\end{equation}

We write
\[
\operatorname{Quot}(E_{q-1},P_{q-1})
\]
for the Grothendieck Quot scheme parameterizing coherent
quotients of the vector bundle $E_{q-1}$ on $\mathbb{P}^{m-1}$ with the
Hilbert polynomial being $P_{q-1}$. This projective parameter space carries a universal quotient
of the pullback of $E_{q-1}$. Its kernel is flat over the
parameter space. By uniform Castelnuovo--Mumford regularity,
we can choose an integer $D_q$ such that the twist by $D_q$
of every kernel in this family is generated by its global sections. In particular, the same $D_q$ works for
$\mathscr K_q=\ker d_{q-1}$ and for the kernels arising
from limiting quotient maps. We use this uniform choice
to construct the next stage of the presentation from the
complete space $H^0(\mathbb{P}^{m-1},\mathscr K_q(D_q))$.
We observe that the degree $D_0$ is regular for every ideal
sheaf with quotient Hilbert polynomial
$P_{\cO/\mathscr I}$; inductively, $D_q$ is regular for every
kernel in $\Quot(E_{q-1},P_{q-1})$.

To state the limiting property, choose a finite affine atlas
$V_\alpha\Subset U_\alpha$ whose smaller sets cover $\mathbb{P}^{m-1}(\mathbb{C})$. In the
standard polynomial frames let $d_{q,\alpha}$ denote the matrices
of the maps $d_q$ in (\ref{eq:dq}). For $p\in V_\alpha$ and a sufficiently small
$\eps>0$, define the sampled Gram forms recursively by
\begin{equation}
\label{eq:Gqalphaepsilon} G_{-1,\alpha,\eps}=1,\qquad
 G_{q,\alpha,\eps}(p)
 =\int_{B_1}d_{q,\alpha}(p+\eps z)^*
       G_{q-1,\alpha,\eps}(p)
       d_{q,\alpha}(p+\eps z)\,d\mu(z),
\end{equation}
where $\mu$ is probability volume on  the unit ball $\mathbb{B}^{N}\subset\C^N$. For every $q$, the forms $G_{q, \alpha, \epsilon}$ in (\ref{eq:Gqalphaepsilon}) are positive definite. With
$R_q=G_{q,\alpha,\eps}(p)^{1/2}$ and $R_{-1}=1$, set
\begin{equation}\label{eq:hatdqpep}
 \widehat d_{q,p,\eps}(\zeta)
   :=R_{q-1}d_{q,\alpha}(p+\eps\zeta)R_q^{-1}.
\end{equation}
Note that the matrices $R_q$ are constant as $\zeta$ varies.

\begin{introB}
For every nonzero proper coherent algebraic ideal sheaf
$\mathscr I\subset\cO_{\Pj^N}$, $N\ge1$, there are degrees
$0\le D_0\le\cdots\le D_N$, chosen with the uniform regularity
properties above, and a presentation
\[
 E_N\xrightarrow{d_N}\cdots\xrightarrow{d_1}E_0
       \xrightarrow{\Phi}\cO_{\mathbb{P}^N},
 \qquad \im d_q=\mathscr K_q.
\]
The normalized tuples
$(\widehat d_{0,p,\eps},\ldots,\widehat d_{N,p,\eps})$ have compact
closure in coefficient topology as $p$ ranges over the fixed 
atlas $\{V_\alpha\}$ and $0<\eps\le\eps_0$.
Homogenize each map to its prescribed degree
$e_q=D_q-D_{q-1}$. Every simultaneous limit, taken in a fixed chart,
is a complex exact on holomorphic stalks at $E_0,\ldots,E_{N-1}$,
and its first image is an
ideal sheaf of Hilbert polynomial $P_0$.
For every normalized or limiting map, the shifted homogeneous
polynomial image $J'_q$ satisfies
\[
 \dim(J'_q)_\ell=
 \begin{cases}
  0,&\ell<D_q,\\
  P_q(\ell),&\ell\ge D_q.
 \end{cases}
\]
The source and target shifts are $D_q$ and $D_{q-1}$, respectively.
\end{introB}

We remark that no injectivity at the last source
$E_N$ is required. The use of complete section spaces is decisive:
properness and uniform regularity make their Grassmannian images
closed, so limiting columns still span the complete sections of a
limiting sheaf. The resulting image Hilbert functions survive even
when the inverses $R_q^{-1}$ become unbounded.

The needed algebraic inputs are the projective correspondence,
Castelnuovo--Mumford regularity, properness of Hilbert and Quot spaces,
and cohomology and base change
\cite{SerreFAC,SerreGAGA,GrothendieckHilbert,StacksRegularity,
StacksQuot,ConradBaseChange}. We construct the presentation in the category of coherent
sheaves. Theorems \ref{thm:C}-\ref{geo:lifting} supply the analytic estimates needed
to apply this construction in the original Hilbert space.

\subsubsection{Step 2: Division, coefficient metrics, and the
\texorpdfstring{$\dbar$}{dbar} equation}
Having obtained compact families of polynomial matrices with fixed image Hilbert functions in Theorem \ref{thm:B}, we next establish uniform holomorphic division for these families. In its statement, we write
$S=\C[Z_0,\ldots,Z_N]$, $\mathbb{B}^N_s=\{z\in\C^N:|z|<s\}$, and
$A_{\mathrm{aff}}(z)=A(1,z)$.
Membership in the \emph{analytic image sheaf on a ball} means that
near every point of the ball the datum has a holomorphic
representation $F=A_{\mathrm{aff}}c$ for some vector valued holomorphic function $c$.

\begin{thmC}
    \label{thm:C}
Let $\mathcal A$ be a compact family, in coefficient topology, of
homogeneous matrices
\[
 A:S(-e)^a\longrightarrow S^b,\qquad e\ge0,
\]
with fixed sizes and degree. Suppose that
$\dim(\im A)_h$ is independent of $A\in\mathcal A$ for every
integer $h$. For every $s>0$, there are $R_s>s$ and $C_s<\infty$
such that every square-integrable holomorphic $F$ in the analytic
image sheaf of $A_{\mathrm{aff}}$ on $\mathbb{B}^N_{R_s}$ has a holomorphic
lift on $\mathbb{B}^N_s$ satisfying
\[
 A_{\mathrm{aff}}u=F|_{\mathbb{B}^N_s},\qquad
 \|u\|_{L^2(\mathbb{B}^N_s)}\le C_s\|F\|_{L^2(\mathbb{B}^N_{R_s})}.
\]
The constants are uniform in $A$ and $F$; for each fixed $A$, the
lift can be chosen complex-linearly in $F$.
The same affine conclusion holds under the weaker assumption that
the image modules have a common eventual Hilbert polynomial, without
a uniform degree at which their Hilbert functions stabilize.
\end{thmC}

We prove this theorem in Section~\ref{sec:division}.
The radius $R_s$ is a part of the conclusion, and image-sheaf membership
is required throughout that input ball. 
The all-degree hypothesis in Theorem \ref{thm:C} needed for each map in the presentation after
removing its common target shift is confirmed by Theorem \ref{thm:B}.

Galligo's parameter-dependent division and flatness results
\cite[Theorems~1.2.7--1.2.8]{Galligo} and Hauser--M\"uller's
quantitative complex-linear scissions
\cite[Theorems~5.1--5.2]{HauserMueller}
are direct antecedents of this division mechanism.
Starting with a standard basis at a reference matrix, we construct a bounded monomial reconstruction in weighted coefficient norms. We invert a small perturbation by a Neumann series and use constancy of image dimensions to remove the remainder. In the homogeneous construction here, the
homogenizing variable has weight $\rho_0=1$. Consequently increasing
the homogenization degree does not increase the coefficient norm or
the division constant. This supplies an estimate uniform both in
degree and on the compact families from Theorem \ref{thm:B}.

We now choose the norms on the coefficient bundles so that local division controls the derivative needed in Step 3. A fixed
smooth coefficient metric does not generally give a lift bounded
uniformly in degree, even for the row
$Z_0:\cO_{\Pj^1}(-1)\to\cO_{\Pj^1}$; the explicit calculation in
Section~\ref{ana:metrics} explains this obstruction. We retain the
original metric on the scalar datum and vary only the coefficient
metrics.

\begin{introD}
For the complete presentation in Theorem \ref{thm:B} and every sufficiently
small $\eps>0$, there are smooth positive Hermitian metrics
$H_{q,\eps}$ on the nonzero bundles $E_q$, with
$H_{-1,\eps}=1$, and a constant $C\ge1$ independent of $\eps$
and the point, such that on each smaller chart
$V_\alpha\Subset U_\alpha$,
\[
 C^{-1}G_{q,\alpha,\eps}(p)[v]
 \le H_{q,\eps}(p)[v]
 \le C G_{q,\alpha,\eps}(p)[v],
 \qquad p\in V_\alpha,\quad v\in(E_q)_p.
\]
The comparisons are on the actual fibers.
For $0\le q\le N$ the maps are contractive, and
\begin{align*}
 |d_qv|_{H_{q-1,\eps}}&\le |v|_{H_{q,\eps}},\\
 |(\partial_{z_j}\Phi)v|_{L^{-1}}
       &\le C\eps^{-1}|v|_{H_{0,\eps}},\\
 \|\Theta(H_{q,\eps})\|_{H_{q,\eps},\omega_{\FS}}
       &\le C\eps^{-2}.
\end{align*}
The derivatives are taken in the fixed homogeneous coordinates, so
$\partial_{z_j}\Phi:E_0\to L^{-1}$.
In particular the absolute Nakano curvature form is bounded by
$C\eps^{-2}$ times the squared norm of its arguments.
\end{introD}

In Theorem \ref{thm:D}, we assert uniform comparison on the smaller charts $V_\alpha$ and measure curvature in the same varying fiber metric that appears in the division and derivative estimates. We start from the sampled Gram forms, patch them on actual bundle fibers, and add graph terms to make the maps contractive. By differentiating polynomials in the sampling variables, we obtain relative first- and second-derivative bounds of orders $\eps^{-1}$ and $\eps^{-2}$, even when the smallest Gram
eigenvalues tend to zero.

Note that metrics and curvature bounds are central in Skoda's $L^2$ division and
bundle-lifting theorems \cite{SkodaDivision,SkodaBundles}, and in
the later work of Varolin and Ji \cite{Varolin,Ji}.
We retain the original scalar Fubini-Study norm on the datum and construct Theorem \ref{thm:D} to satisfy the three requirements of local division, differentiated-row control, and a curvature bound in one coefficient geometry.

We combine these ingredients to obtain the global lift.

\begin{introE}\label{thm:E}
Let $\mathscr I\subset\cO_{\mathbb{P}^N}$ be a nonzero proper coherent algebraic ideal
sheaf on $\Pj^N$, $N\ge1$, and choose
\[
 D=D_0,\qquad
 V=H^0(\Pj^N,\mathscr I\otimes L^D),\qquad
 \Phi:E_0=V\otimes L^{-D}\longrightarrow\cO
\]
as in Theorem \ref{thm:B}, with the metrics of Theorem \ref{thm:D}.
There are fixed constants $R_*>0$, $C_{\mathrm{lift}}$, and an
integer $n_0\ge1$ such that, for every integer $n\ge n_0$,
$\eps=R_*/\sqrt n$, and $h\in\{n,n+1\}$, there is a linear map
\[
 T_{n,h}:H^0(\Pj^N,\mathscr I\otimes L^h)
           \longrightarrow H^0(\Pj^N,E_0\otimes L^h)
\]
satisfying
\[
 \Phi T_{n,h}F=F,\qquad
 \|T_{n,h}F\|_{H_{0,\eps}\otimes L^h,L^2}
       \le C_{\mathrm{lift}}\|F\|_{L^h,L^2}.
\]
The right-hand norm is the original scalar Fubini--Study norm.
The presentation and constants are independent of $n,h,F$.
The same coefficient metric serves the two adjacent powers.
For the unit ideal sheaf, the row $1$ and the constant metric give
the same conclusion.
\end{introE}

In proving Theorem~\ref{geo:lifting}, we first apply local division to holomorphic overlap data and use exactness of the presentation to lift them successively through the syzygy maps. We then use a partition of unity to construct smooth forms of antiholomorphic degrees $0,\ldots,N$.
Each differentiated cutoff costs $O(\eps^{-1})$. For $\dbar$-closed
forms of positive antiholomorphic degree with values in the smooth
bundles $E_q\otimes L^h$, the positive-twist estimate has the form
\[
 \|\dbar_{\min}^{-1}\|
 \le \bigl(c_Lh-C_0-C_1\eps^{-2}\bigr)^{-1/2},
 \qquad c_Lh-C_0-C_1\eps^{-2}>0.
\]
We obtain this bound from the Bochner--Kodaira--Nakano estimate, with the usual
canonical-bundle adjustment for $(0,j)$-forms
\cite{Demailly,DemaillyL2}. We first fix the presentation and curvature constants, and then choose $R_*$ sufficiently large. At
$\eps=R_*/\sqrt n$ and $h=n$ or $n+1$, solving the $\dbar$-equation gains
$O(\eps)$. By a finite descent we cancel the powers lost to the cutoffs and obtain the bound in Theorem \ref{geo:lifting}. The last $(0,N)$-form is
automatically $\dbar$-closed, which explains why exactness at the
last source $E_N$ is unnecessary. We need no smallness assumption on the local division constants.

The descent belongs to the homological $\dbar$ tradition of
H\"ormander and Kelleher--Taylor
\cite{HormanderGenerators,KelleherTaylor}. Andersson's membership
theory \cite{AnderssonMembership} and, especially,
Andersson--Wulcan's projective resolution and Dolbeault construction
\cite[Section~6, Proposition~6.1]{AnderssonWulcan} are close
precedents. Theorem \ref{geo:lifting} adds the stated degree-uniform $L^2$ bound
using the coefficient metrics of Theorem \ref{thm:D}. We solve the equations on smooth bundles over the ambient projective space, retaining the singular and nonreduced data in their sheaf maps.

\subsubsection{Step 3: Back to the original Drury--Arveson space} 
For sufficiently large $n$, so that $n\ge n_0$ and
$I_n=(I^{\mathrm{sat}})_n$, $I_{n+1}=(I^{\mathrm{sat}})_{n+1}$,
take $F\in I_{n+1}$ and $a=T_{n,n+1}F$. The identity $F=\Phi a$ gives
\[
 \partial_{z_j}F=(\partial_{z_j}\Phi)a+\Phi\,\partial_{z_j}a.
\]
The second term belongs to $I_n$, so by orthogonal projection we obtain
\[
 P_{Q_n}\partial_{z_j}F
   =P_{Q_n}\bigl((\partial_{z_j}\Phi)a\bigr).
\]
Theorems \ref{thm:D} and \ref{geo:lifting}, applied in the same coefficient metric, imply
\[
 \|P_{Q_n}\partial_{z_j}F\|_{\FS,n}
 \le C\eps^{-1}\|F\|_{\FS,n+1}.
\]
Since $Z_j^*=(n+1)^{-1}\partial_{z_j}$ on degree $n+1$, the exact
norm ratio in \eqref{eq:intro-projective-norm} gives
\[
 \|B_{j,n}^*\|
 \le\frac{C}{\eps\sqrt{(n+1)(n+m)}}=O(n^{-1/2}).
\]
Summing the squared bounds, we obtain the scalar estimate
$A_n\le C_I(n+1)^{-1}I$; we absorb the remaining finitely many degrees into $C_I$. This is the division-and-projection
mechanism of Shalit and Eschmeier, with the required derivative
estimate furnished by the projective lifting construction.

For the vector module $M\subset R\otimes\C^r$, we introduce
$y_1,\ldots,y_r$ and the homogeneous ideal
\[
 \lambda(f)=\sum_{b=1}^r f_b(z)y_b,\qquad
 J_M=(\lambda(M))+(y_1,\ldots,y_r)^2
       \subset R[y_1,\ldots,y_r].
\]
The auxiliary-variable realization has a direct antecedent in
\cite[Section~5, equation~(5.1)]{ShalitDivision}. With
\[
 U_nq=\sqrt{n+1}\sum_{b=1}^r q_b(z)y_b,
\]
the linear auxiliary-degree sector is unitarily identified with
$Q_n$, and the square-zero relations give the exact identity
\[
 U_n^*A_{n+1}^{J_M}U_n
   =\frac{n+1}{n+2}A_n^M+\frac{r+1}{n+2}I_{Q_n}.
\]
Thus we transfer the scalar estimate to the original vector quotient before counting singular values. For infinite-dimensional
$T=(R\otimes\C^r)/M$, Hilbert--Serre gives
$\dim Q_n=O(n^{\delta-1})$ and hence the weak Schatten exponents
in Theorem \ref{thm:main}. Counting in the auxiliary scalar quotient would use
its support dimension $m$ and could lose the sharper exponent
$\delta$.

We link these estimates through a specific norm requirement. With the complete-section construction we preserve the algebraic image spaces at normalized limits; with Theorem \ref{thm:C} we obtain division uniform in degree; with Theorem \ref{thm:D} we control both differentiation and curvature; and with Theorem \ref{geo:lifting} we produce a lift measured in that same coefficient metric with the original scalar norm on the datum. We finally transfer the estimates back to the given shift by means of the projected derivative.

We are currently working to extend the methods developed here to the Arveson–Douglas problem for general, nonhomogeneous polynomial ideals. Quantitative division, compatible Hermitian metrics, and \(\bar\partial\) estimates in this paper provide a natural starting point for this investigation. The main difficulty is the absence of the homogeneous decomposition, which in this paper allows us to combine degreewise analytic estimates with algebraic dimension growth. For nonhomogeneous ideals, the geometry near the boundary of the unit ball becomes particularly relevant, including possible singularities, nonreduced structures, and the interaction of different components. Our ongoing work focuses on adapting the tools developed in the proof of Theorem \ref{thm:main} to neighborhoods of the boundary of the unit ball, with the aim of obtaining uniform local estimates that replace the degreewise bounds used in the homogeneous setting. The progress will be reported in the near future. 

The paper is organized as follows. 
In Section~\ref{sec:operators} we prove the compression and counting lemmas, and in Section~\ref{sec:projective} we establish the projective norm and saturation identifications. We prove Theorems \ref{thm:B}-\ref{geo:lifting} in Sections~\ref{sec:algebra}--\ref{ana:global-section} and complete the proof of Theorem \ref{thm:main} in Section~\ref{sec:completion}. In Section~\ref{sec:extensions} we treat other standard ball kernels, analytic ideals generated by homogeneous functions, and exact examples.\\

\noindent{\bf Notations}: We use $\N=\{0,1,2,\ldots\}$. All support dimensions are complex
dimensions. We also view algebraic sheaf maps on projective space as holomorphic maps after analytification. Constants may depend on
the fixed ideal, presentation, charts, and radii, but not on the
large degree or the datum. We omit zero coefficient stages and all subsequent stages, and treat the elementary one-variable scalar case separately.

\medskip

\noindent\textbf{Note added on November 30, 2026.}
We began circulating an initial version of this paper among a number
of experts in the field on November 20, 2026. While continuing to
refine the exposition, we made a version of the manuscript, together
with a research announcement, publicly available on
\href{https://github.com/Xiang-Tang-Math/Arveson-Douglas-Conjecture}{GitHub}
on November 25, 2026.
On November 30, 2026, some of us received an email announcement
of another proof of the homogeneous Arveson--Douglas conjecture
by Chris Bruce and Adam Dor-On, presented in their preprint
\emph{A proof of the Arveson--Douglas conjecture}
(\href{https://arxiv.org/abs/2611.37464v1}{arXiv:2611.37464}),
submitted to arXiv on November 27, 2026.
We welcome this development and, in light of it, have decided
to make our manuscript available on arXiv as well.


\section{Graded quotients and coordinate defects}
\label{sec:operators}

We first separate the operator estimate from the algebraic count.  This
also fixes the original Hilbert-space norms in which every later
comparison will be made.

\subsection{Algebraic terminology and the dimension parameter}
\label{op:algebra-background}
We recall the elementary algebra used here. A graded module is a
module $T=\bigoplus_{n\in\mathbb Z}T_n$ over the graded polynomial
ring $R=\bigoplus_{n\geq0}R_n$, with $R_aT_b\subset T_{a+b}$.
Our polynomial quotients are generated in degree zero. The shift
$T(-d)$ is the same module with $T(-d)_n=T_{n-d}$. A ring is
Noetherian if its ideals satisfy the ascending chain condition, or,
equivalently, if every ideal is finitely generated. Hilbert's basis
theorem makes $R$ Noetherian; submodules of a finite free module
are then finitely generated. We use this fact, not any assumption of
finite generation of a closed Hilbert module by a bounded division
operator.

For an ideal $J\subset R$, its radical is
$\sqrt J=\{f\in R:f^k\in J\text{ for some }k\geq1\}$. The affine
scheme defined by $J$ is $\Spec(R/J)$, the space of prime ideals of
$R/J$ equipped with the structure sheaf of local rings. Its closed
complex points are the common zeros of $J$. Its reduced scheme is
$\Spec(R/\sqrt J)$. An element of $R/J$ whose positive power is zero
is called nilpotent, and the scheme is reduced precisely when it
has no nonzero nilpotent local sections. For example, $(x)$ and
$(x^2)$ have the same zero set, but the class of $x$ survives in the
second quotient as a nilpotent element. The two ideals give different
Hilbert modules. We retain this distinction throughout.

For a finite $R$-module $T$, its annihilator is an ideal and its
algebraic support consists of the primes $\mathfrak p$ for which
$T_{\mathfrak p}\ne0$. Here localization at $\mathfrak p$ allows
all elements of $R\setminus\mathfrak p$ as denominators. To see the
support description, if an annihilator element lies outside
$\mathfrak p$, it becomes invertible and kills $T_{\mathfrak p}$.
Conversely, if $T_{\mathfrak p}=0$, each of finitely many generators
is killed by an element outside $\mathfrak p$; their product kills
$T$ and is still outside $\mathfrak p$. Hence the support is exactly
$V(\operatorname{Ann}_R T)$. Its Krull dimension is the supremum of
lengths of strict chains of prime ideals in the quotient ring.
Over $\C$, for these finitely generated rings this agrees with the
complex algebraic dimension of the closed-point variety. We use no
dimension convention for the zero module: its quotient Hilbert space
is simply zero.

An associated prime of $T$ is a prime equal to the annihilator of
some nonzero element of $T$. It is embedded if it is not a minimal
prime of the support. Such primes can occur even when the underlying
zero set has just one irreducible component. We give a concrete defect
calculation for an embedded prime in Section~\ref{sec:extensions}.

We call $h_T(n)=\dim_\C T_n$ the Hilbert function. The Hilbert--Serre
theorem for finite graded modules over a standard polynomial ring
asserts that $h_T(n)$ is eventually a polynomial, of degree
$\dim\Supp T-1$ and positive leading coefficient if the support
has positive dimension. For a nonzero zero-dimensional module the
function is eventually zero. This is the form of the classical theorem
used below; the eventual polynomial is called the Hilbert polynomial.
The Hilbert series $\sum_n h_T(n)t^n$ is rational with denominator a
power of $1-t$ and, in positive support dimension, its pole at $1$
has order $\dim\Supp T$. In particular, a bound on the Hilbert function
is a count of the ranks of the blocks and is not an estimate on their
operator norms. We refer to \cite{SerreFAC,StacksHilbert,StacksNoetherian}
for these algebraic facts.

\subsection{The grading and its polynomial part}
For a multi-index $\alpha$, we put $|\alpha|=\sum_j\alpha_j$ and
$\alpha!=\prod_j\alpha_j!$.  The vectors
\[
 e_\alpha=\sqrt{\frac{|\alpha|!}{\alpha!}}z^\alpha
\]
form an orthonormal basis of $H_m^2$.  If $|\alpha|=n$, then
\begin{equation}\label{eq:op-shift}
 Z_je_\alpha=\sqrt{\frac{\alpha_j+1}{n+1}}e_{\alpha+e_j},
 \qquad
 Z_j^*z^\alpha=\frac{\alpha_j}{n}z^{\alpha-e_j}\quad(n>0),
\end{equation}
with the second expression interpreted as zero if $\alpha_j=0$.
In particular $\|Z_j\|\le1$.  Tensoring with $\C^r$ preserves these
formulas.  Let $H_n=H_{m,n}^2\otimes\C^r$ and let $P_n$ denote its
orthogonal projection.

\begin{lemma}\label{op:graded}
A closed subspace $\mathcal M\subset H_m^2\otimes\C^r$ invariant under
the coordinates and under $U_t f(z)=f(e^{it}z)$ is the closure of a
finitely generated homogeneous submodule $M\subset R\otimes\C^r$.
Its polynomial part is exactly $M$.  Conversely, the closure of a
homogeneous polynomial submodule has these properties.
\end{lemma}
\begin{proof}
The circle action is unitary and strongly continuous, as follows first
on polynomials and then by density.  We express the Fourier projection as the norm-convergent Bochner integral
\[
 P_nf=\frac1{2\pi}\int_0^{2\pi}e^{-int}U_tf\,dt.
\]
For $f\in\mathcal M$, its Riemann sums belong to $\mathcal M$, so
$P_nf\in\mathcal M\cap H_n$.  The partial sums $\sum_{n=0}^NP_nf$
converge to $f$.  Thus $M=\mathcal M\cap(R\otimes\C^r)$ is homogeneous
and dense in $\mathcal M$.  Coordinate invariance makes it an
$R$-submodule.  Since $R$ is Noetherian and $R\otimes\C^r$ is finite
free, $M$ is finitely generated.  We obtain homogeneous generators by replacing each generator by its finitely many homogeneous parts; all those parts
belong to $M$ and their sums recover the original generators.

Conversely, boundedness of the coordinates and unitarity of $U_t$ preserve
both invariances upon taking closure.  If a polynomial $f$ belongs to
$\overline M$, then $P_nf$ lies in the closure of $M_n$ for each $n$.
The latter space is finite dimensional and closed.  Hence every
homogeneous part of $f$ belongs to $M$, and their finite sum does too.
\end{proof}

We write $Q=\mathcal M^\perp=\widehat\bigoplus_{n\ge0}Q_n$.
Invariance of $\mathcal M$ implies $Z_j^*Q\subset Q$.  Indeed,
$\langle Z_j^*q,u\rangle=\langle q,Z_ju\rangle=0$ for
$q\in Q$ and $u\in\mathcal M$.  Thus the quotient tuple
$S_j=P_QZ_j|_Q$ has adjoint $Z_j^*|_Q$.

\subsection{Compression identities}
The degree-$n$ defect column is
\[
 B_n:Q_n\longrightarrow M_{n+1}^{\oplus m},\qquad
 B_nq=(B_{1,n}q,\ldots,B_{m,n}q),
\]
so $A_n=B_n^*B_n$.  Its direct sum is a bounded operator because
each coordinate is bounded.

\begin{lemma}\label{op:commutators}
For $n\ge0$,
\begin{align}
 [Z_i^*,Z_j]|_{H_n}
   &=\frac{\delta_{ij}I-Z_jZ_i^*}{n+1}\bigg|_{H_n},
       \label{eq:op-ambient}\\
 [S_i^*,S_j]|_{Q_n}
   &=P_{Q_n}[Z_i^*,Z_j]|_{Q_n}-B_{i,n}^*B_{j,n}.
       \label{eq:op-quotient}
\end{align}
For $n\ge1$, the restriction tuple satisfies
\begin{equation}\label{eq:op-restriction}
 [V_i^*,V_j]|_{M_n}
 =P_{M_n}[Z_i^*,Z_j]|_{M_n}+B_{j,n-1}B_{i,n-1}^*.
\end{equation}
The last term is zero in degree zero.  In particular, a bound
$A_n\le C(n+1)^{-1}I$ gives $O((n+1)^{-1})$ bounds for all
quotient and restriction commutator blocks.
\end{lemma}
\begin{proof}
For $n>0$, we apply \eqref{eq:op-shift} to a homogeneous polynomial:
\[
 Z_i^*Z_jf=\frac{\delta_{ij}f+z_j\partial_if}{n+1},\qquad
 Z_jZ_i^*f=\frac{z_j\partial_if}{n}.
\]
Taking the difference, we obtain \eqref{eq:op-ambient}; constants give the same
identity in degree zero.  Relative to $\mathcal M\oplus Q$, we write
\[
 Z_j=\begin{pmatrix}V_j&B_j\\0&S_j\end{pmatrix}.
\]
Multiplying the two matrices gives
\[
 Z_i^* Z_j=\begin{pmatrix}
 V_i^*V_j&V_i^*B_j\\ B_i^*V_j&B_i^*B_j+S_i^*S_j
 \end{pmatrix},\qquad
 Z_jZ_i^*=\begin{pmatrix}
 V_jV_i^*+B_jB_i^*&B_jS_i^*\\ S_jB_i^*&S_jS_i^*
 \end{pmatrix}.
\]
We subtract these matrices and solve for the diagonal commutators.
The term $B_i^*B_j$ acts in $Q_n$ through $M_{n+1}$, whereas
$B_jB_i^*$ acts in $M_n$ through $Q_{n-1}$; this gives exactly
\eqref{eq:op-quotient} and \eqref{eq:op-restriction}. In degree zero
there is no incoming block.

Since $B_{j,n}^*B_{j,n}\le A_n$, each
$\|B_{i,n}^*B_{j,n}\|$ is at most $\|A_n\|$.  The ambient block
has norm at most $2/(n+1)$.  The same argument in incoming degree
$n-1$ gives the restriction bound, after increasing a constant to include
degree zero.
\end{proof}

The defect is an operator-norm question.  Its trace alone does not give
the needed bound in higher dimension.  In fact an exact identity is
available without any division estimate.

\begin{proposition}\label{op:trace}
If $h_n=\dim Q_n$, then
\begin{equation}\label{eq:op-trace}
 \Tr A_n=\frac{n+m}{n+1}h_n-h_{n+1}.
\end{equation}
\end{proposition}
\begin{proof}
The forward weights in \eqref{eq:op-shift} give
$\sum_jZ_j^*Z_j|_{H_n}=(n+m)(n+1)^{-1}I$.
Splitting each forward product orthogonally into its ideal and quotient parts, we obtain
\[
 A_n=\frac{n+m}{n+1}I_{Q_n}-\sum_j S_{j,n}^*S_{j,n}.
\]
On positive degrees the backward weights give
$\sum_jZ_jZ_j^*=I$.  Since $Q$ is coinvariant,
$\sum_jS_{j,n}S_{j,n}^*=I_{Q_{n+1}}$.  Taking finite-dimensional traces, we obtain the formula.
\end{proof}

Suppose that the quotient is infinite dimensional, so that its support
dimension satisfies $1\leq\delta\leq m$, and write $a>0$ for the leading
coefficient of its eventual Hilbert polynomial (of degree $\delta-1$
by Hilbert--Serre; see Lemma~\ref{op:hilbert}). Rewriting \eqref{eq:op-trace}, we obtain
\begin{equation}\label{eq:op-trace-asymptotic}
 \Tr A_n=\frac{m-1}{n+1}h_n-(h_{n+1}-h_n)
       =(m-\delta)a n^{\delta-2}+O(n^{\delta-3}).
\end{equation}
Indeed, if $\delta\geq2$, the polynomial difference has leading term
$(\delta-1)a n^{\delta-2}$ and remainder $O(n^{\delta-3})$; if
$\delta=1$, the difference is eventually zero and the same expansion
holds. Thus $\Tr A_n\sim(m-\delta)a n^{\delta-2}$ when $\delta<m$.
When $\delta=m$ the leading coefficient vanishes; the trace is then
$O(n^{m-3})$ and may be identically zero, as for $M=0$. In all cases the
average eigenvalue, $\Tr A_n/h_n$, is $O(n^{-1})$. For $\delta=1$, $h_n=a$
eventually, so positivity gives the full estimate
$\|A_n\|\leq (m-1)a/(n+1)$ for large $n$. This is the dimension-one
trace argument of Guo and Wang for finite coefficient multiplicity
\cite[Proposition~4.1, p.~920]{GuoWang}. In higher support dimension the
trace estimate generally permits concentration in a small subspace.
Below we obtain the required uniform operator-norm bound by the lifting argument, without a restriction on the support dimension.

\subsection{Counting singular values in the original grading}
For a compact operator $T$, let $s_k(T)$ be its singular values in
nonincreasing order, starting at $k=1$.  For $p>0$, write
\[
 T\in\cS_{p,\infty}\quad\Longleftrightarrow\quad
 \sup_{k\ge1} k^{1/p}s_k(T)<\infty.
\]

\begin{lemma}\label{op:counting}
Suppose $T=\bigoplus_{n\ge0}T_n$ has mutually orthogonal domain and
range blocks, with
\[
 \dim\operatorname{dom}T_n\le C_0(n+1)^{\delta-1},\qquad
 \|T_n\|\le C_1(n+1)^{-a},
\]
where $a,\delta>0$.  Then $T$ is compact,
$T\in\cS_{\delta/a,\infty}$, and $T\in\cS_p$ for $p>\delta/a$.
\end{lemma}
\begin{proof}
Finite block truncations are finite rank and converge in norm because
the block norms tend to zero.  If a singular value of $T_n$ exceeds
$t>0$, then $n+1<(C_1/t)^{1/a}$.  We therefore have the following estimate for the counting function:
\[
 \#\{k:s_k(T)>t\}
 \le\sum_{n+1<(C_1/t)^{1/a}}C_0(n+1)^{\delta-1}
 \le C_2t^{-\delta/a}
\]
for small $t$, which is the weak Schatten assertion.  Alternatively,
for $p>\delta/a$,
\[
 \sum_n\|T_n\|_{\cS_p}^p
 \le C_0C_1^p\sum_n(n+1)^{\delta-1-ap}<\infty.
\]
\end{proof}

\begin{lemma}\label{op:hilbert}
Let $T=(R\otimes\C^r)/M$.  If $T$ is infinite dimensional and
$\delta=\dim\Supp T$, then
\[
 \dim Q_n=\dim T_n=O((n+1)^{\delta-1}).
\]
If $\delta=0$, then $T$ is finite dimensional.
\end{lemma}
\begin{proof}
In degree $n$, orthogonal projection gives an isomorphism
$T_n\simeq Q_n$ of finite-dimensional vector spaces.  The
Hilbert--Serre theorem states that the Hilbert function of a nonzero
finite graded module of positive support dimension is eventually a
polynomial of degree $\delta-1$, with positive leading coefficient.
Equivalently, its Hilbert series has a pole of order $\delta$ at $1$.
This applies to $T$ because it is a quotient of a finite free standard
graded module.  We absorb the finitely many earlier degrees into a constant.  In dimension zero its Hilbert function is eventually zero,
so its total vector-space dimension is finite; see
\cite{SerreFAC,StacksHilbert}.
\end{proof}

Applying Lemma~\ref{op:counting} with $a=1$ to
\eqref{eq:main-commutator} proves the quotient assertions of
Theorem~\ref{thm:main}, once \eqref{eq:main-degree} is established.
For the column $B$, we use $\|B_n\|^2=\|A_n\|$ and $a=1/2$.
The commutator $[P_{\mathcal M},Z_j]$ is the block $B_j$ extended by
zero, so it has the same nonzero singular values.  For the restriction tuple we use \eqref{eq:op-restriction} and
\[
 \dim M_n\le r\binom{n+m-1}{m-1}=O((n+1)^{m-1}).
\]
We obtain its weak exponent $m$.  If $T$ is finite dimensional then
$Q$ is finite dimensional, and the assertions involving $Q$ have finite
rank.  Notice that the restriction exponent and the quotient exponent
come from different dimension counts.


\section{Step 1: Projective norms and the lifting problem}
\label{sec:projective}

We now treat a scalar homogeneous ideal, making explicit the relation between the holomorphic Hilbert space and the geometric $L^2$ spaces. In particular, we later solve $\bar\partial$ equations for bundle-valued forms on projective space.

\subsection{Projective coordinates, sheaves, and line bundles}
\label{geo:background}
We recall the projective construction explicitly. For $N=m-1$, we
write homogeneous coordinates as $[Z_0:\cdots:Z_N]$, a nonzero vector
modulo common nonzero scalar multiplication. On the chart $Z_k\ne0$
we use coordinates $w_j=Z_j/Z_k$ for $j\ne k$. This is a copy of
$\C^N$ with polynomial coordinate ring in the algebraic category.
The structure sheaf $\cO_{\mathbb{P}^N(\mathbb{C} )}$ assigns regular functions to algebraic
open sets (and holomorphic functions to analytic open sets after
analytification). A sheaf assigns sections to open sets, with
restriction maps and the property that compatible local sections
glue uniquely. Its stalk $\mathscr F_p$ consists of germs of
sections near $p$; two representatives define the same germ if they
agree on a smaller neighborhood. A sheaf of $\cO_{\mathbb{P}^N}$-modules is
coherent when locally it admits a finite presentation by free
$\cO_{\mathbb{P}^N}$-modules. We use the standard coherence theorems, which imply
that kernels, images, and cokernels of maps between coherent sheaves
are coherent. Exactness of a sequence of sheaves means exactness of
its sequence of stalk modules, not exactness after taking global
sections or evaluating in point fibers.

A holomorphic vector bundle is locally a product with a complex
vector space, with holomorphic invertible changes of frame. Its
sheaf of holomorphic sections is locally free: locally it is
isomorphic to $\cO_{\mathbb{P}^N}^a$ for a fixed rank $a$. Conversely such a
sheaf determines a holomorphic vector bundle. A line bundle has
rank one. The tautological bundle $\cO_{\mathbb{P}^N}(-1)$ has fiber $\C Z$ at
$[Z]$, and $L=\cO_{\mathbb{P}^N}(1)$ is its dual. For an integer $t$ we write
$L^t$ for the tensor power if $t\geq0$ and the corresponding power
of the dual for $t<0$; $L^0=\cO$. The coordinate functional $Z_k$
is a nowhere-zero local frame of $L$ on $Z_k\ne0$. A homogeneous
polynomial $f$ of degree $t\geq0$ gives a section of $L^t$ whose
coefficient in this frame is $f/Z_k^t$. Conversely, in the algebraic category a section over $Z_k\ne0$
is represented by an element of $R[Z_k^{-1}]$ homogeneous of
degree $t$ (we multiply its regular frame coefficient by $Z_k^t$).
Compatibility on overlaps means these elements are the same rational
function in the fraction field of $R$. If $N\ge1$, the intersection
of $R[Z_k^{-1}]$ over all $k$ is $R$: in lowest terms the denominator
would have to divide a power of every $Z_k$, and two distinct
coordinate variables are relatively prime. The common function
therefore belongs to $R_t$. This also proves vanishing for negative
$t$. When $N=0$ and $t\ge0$ the assertion reduces to the unique
monomial $Z_0^t$ spanning the sections on a point.
GAGA, recalled below, gives the same descriptions for holomorphic
sections. We write $\mathscr F(t)=\mathscr F\otimes L^t$ for the
corresponding twist of a sheaf. Tensoring by a line bundle is exact,
since in each trivialization it is tensoring by a free rank-one
module.

For a homogeneous ideal $I$, its sheafification $\widetilde I$ on
$\mathbb{P}^N(\mathbb{C})=\Proj R$ is defined on $Z_k\ne0$ by the degree-zero part of the
localized ideal $I[Z_k^{-1}]$; localization allows powers of $Z_k$
in denominators. If $f_1,\ldots,f_a$ are homogeneous generators of
degrees $d_1,\ldots,d_a$, its local generators are the polynomials
$f_j/Z_k^{d_j}$ in the $w$ coordinates. The localized definitions
agree on overlaps, so $\widetilde I$ is a coherent ideal sheaf. Its
quotient $\cO/\widetilde I$ is the structure sheaf of the projective
subscheme determined by $I$. In particular the sheaf stores
relations and nilpotent structure, rather than only a subset of
projective space. It need not be the sheaf of all functions
vanishing on its set of zeros. For a local holomorphic section $F$
to belong to an ideal or image sheaf means that its germ has a
holomorphic representation at every point. For instance, $z$ and
$z^2$ vanish at the same point, but $z$ does not belong to the
analytic ideal $(z^2)$ at zero. This local distinction is essential
to our division theorem.

We use Serre's GAGA comparison \cite{SerreGAGA} in two forms.
Analytification of coherent algebraic sheaves on the complex
projective space is exact, and their algebraic cohomology groups
identify with the analytic cohomology groups of the analytified
sheaves. In particular the polynomially defined sheaf maps have
the same exactness on holomorphic stalks, and their holomorphic
global sections on projective space have the stated algebraic
description. These facts allow us to construct our complexes
algebraically and solve differential equations on their analytic
bundles.

\subsection{Saturation changes only finitely many blocks}
Let $I\subset R=\C[z_1,\ldots,z_m]$ be homogeneous and let
$\mathfrak m=(z_1,\ldots,z_m)$.  We put
\[
 I^{\mathrm{sat}}=\bigcup_{k\ge0}(I:\mathfrak m^k),\qquad
 (I:\mathfrak m^k)=\{f\in R:\mathfrak m^k f\subset I\}.
\]
The homogeneous maximal ideal $\mathfrak m$ is called the irrelevant
ideal for the projective construction. Saturation removes exactly the
classes of $R/I$ killed by some power of $\mathfrak m$; equivalently,
$I^{\mathrm{sat}}/I$ is the submodule of such classes. Since some
coordinate is invertible on each standard chart, $I$ and its
saturation give identical localized ideals there.

\begin{lemma}\label{geo:saturation}
The graded module $I^{\mathrm{sat}}/I$ is finite dimensional.
Consequently $I_n=I_n^{\mathrm{sat}}$ for all sufficiently large $n$,
and a bound $A_n\le C(n+1)^{-1}I$ for the saturation gives the same
type of bound for $I$ after increasing $C$.
\end{lemma}
\begin{proof}
The ascending chain of colon ideals stabilizes, since $R$ is Noetherian.
Thus $\mathfrak m^k I^{\mathrm{sat}}\subset I$ for one $k$.
The quotient $I^{\mathrm{sat}}/I$ is finitely generated and is annihilated
by $\mathfrak m^k$, so it is a finite module over the finite-dimensional
algebra $R/\mathfrak m^k$.  It has only finitely many nonzero degree
spaces.  In high degree the orthogonal projections onto $I_n$ and onto
its saturation are literally equal in the original inner product.
The defect in degree $n$ uses only degrees $n$ and $n+1$, so those
defects coincide eventually.  We absorb the norms of the finitely many remaining finite-dimensional blocks into $C$.
\end{proof}

For example, in $\C[x,y]$ the ideal $I=(x^2,xy)$ has saturation $(x)$.
Indeed $\mathfrak m x\subset I$, and $(x)$ is saturated since
$\C[x,y]/(x)\simeq\C[y]$ has no $\mathfrak m$-torsion. The quotient
$(x)/I$ is the one-dimensional space spanned by the class of $x$.
In this example the projective scheme is the reduced point $[0:1]$;
the finite-length discrepancy occurs only at the cone vertex. In general,
saturation preserves the entire projective ideal sheaf, including any
nonreduced structure supported on projective space.

For $m=1$, a scalar homogeneous ideal is $0$ or $(z_1^a)$ for some
$a\ge0$.  The first has zero defect; the others have finite-dimensional
quotient.  Likewise, $I^{\mathrm{sat}}=R$ gives a finite quotient in
any dimension.  We may therefore assume $m\ge2$ and that the saturated
ideal sheaf is nonzero and proper.  The zero ideal has zero defect.

We write $N=m-1$ and $L=\cO_{\mathbb{P}^N}(1)$.  The saturation determines
an algebraic coherent ideal sheaf $\cI\subset\cO_{\mathbb{P}^N}$, whose
analytification is also denoted $\cI$.  In nonnegative degree,
\begin{equation}\label{eq:geo-sections}
 I^{\mathrm{sat}}_n=H^0(\mathbb{P}^N,\cI\otimes L^n)\subset H^0(\mathbb{P}^N,L^n)=R_n.
\end{equation}
For completeness, an algebraic section represented by $p\in R_n$
belongs to $\cI(n)$ precisely when each standard localization places
it in the localized ideal.  Thus for each $i$ some $z_i^{k_i}p$ lies
in $I$, with $k_i$ chosen positive.  Every monomial of degree greater than
$\sum_i(k_i-1)$ is divisible by one of these powers, so
$\mathfrak m^kp\subset I$ for a sufficiently large $k$.
This is saturation.  The converse is immediate by localization.
We use GAGA to identify these algebraic sections with holomorphic global sections.
Analytification also preserves exact sequences of coherent sheaves;
these are the exactness and global-section comparisons used below
\cite{SerreGAGA}.

\subsection{The original norm as a projective integral}
We give $L$ its Fubini--Study metric and normalize the Fubini--Study volume
$dV$ on $\mathbb{P}^N(\mathbb{C})$ to have total mass one.  A degree-$n$ homogeneous polynomial
$p$ defines a section of $L^n$ with pointwise norm
\[
 |p([z])|_{L^n}^2=\frac{|p(z)|^2}{|z|^{2n}}.
\]
The quotient is unchanged under $z\mapsto\lambda z$.

We recall the metric and volume conventions that we will also use for
differential forms. The standard Euclidean metric on $\C^{N+1}$
induces the metric of the tautological line and its dual. In the
frame $Z_k$ the squared norm is $(1+|w|^2)^{-1}$. Its potential is
$\varphi(w)=\log(1+|w|^2)$. The Fubini--Study K\"ahler form can be
normalized as $\omega_{\FS}=i\partial\bar\partial\varphi$; any fixed
positive normalization only changes fixed constants below. Here
$\partial$ and $\bar\partial$ are the $(1,0)$ and $(0,1)$ parts of
the exterior derivative. The Hermitian matrix of this form is
\[
 g_{a\bar b}(w)=
 \frac{(1+|w|^2)\delta_{ab}-\overline w_a w_b}{(1+|w|^2)^2}.
\]
By Cauchy--Schwarz it is positive definite. Its associated Riemannian
metric supplies geodesic distance and the norms of covectors and
forms; $dV$ is the volume proportional to $\omega_{\FS}^N/N!$
with total mass one. In the chart, on each relatively compact
subset, this metric and its volume density are bounded above and
below by positive constants relative to Euclidean ones. We keep
this base metric fixed, and later vary only the metrics of the
coefficient bundles. This convention prevents a change of the
original norm on the datum.

\begin{lemma}\label{geo:norm}
For $p\in R_n$,
\begin{equation}\label{eq:geo-norm}
 \|p\|_{\DA,n}^2=c_n\|p\|_{\FS,n}^2,
 \qquad
 c_n=\frac{(n+m-1)!}{(m-1)!\,n!}.
\end{equation}
Moreover,
\begin{equation}\label{eq:geo-ratio}
 \frac{c_n}{c_{n+1}}=\frac{n+1}{n+m}.
\end{equation}
Orthogonal projections within a fixed degree are the same for the two
inner products.
\end{lemma}
\begin{proof}
The Hopf projection pushes normalized sphere measure $d\sigma$ on
$S^{2m-1}$ to $dV$.  Hence the projective squared norm of $p$ is
$\int_{S^{2m-1}}|p|^2d\sigma$.  For $|\alpha|=n$, the Gaussian
integral factors as
\[
 \int_{\C^m}|z^\alpha|^2e^{-|z|^2}\,d\lambda(z)=\pi^m\alpha!.
\]
In polar coordinates its radial factor, including the area of the unit
sphere, is $\pi^m(n+m-1)!/(m-1)!$.  Therefore
\[
 \|z^\alpha\|_{\FS,n}^2
 =\frac{(m-1)!\,\alpha!}{(n+m-1)!}.
\]
Coordinatewise rotations make distinct monomials orthogonal for this
norm.  By comparison with $\alpha!/n!$ we obtain \eqref{eq:geo-norm}, and taking the ratio gives \eqref{eq:geo-ratio}.  A scalar change of inner product leaves
orthogonal projections unchanged.
\end{proof}

\begin{remark}\label{geo:dbar-space}
The Drury--Arveson norm is not an ordinary positive-measure norm on the
closed ball when $m\ge2$.  Such a measure would have mass
$\|1\|^2=1$, but would also satisfy
$\int|z|^2d\mu=\sum_j\|z_j\|^2=m>1$.
Lemma~\ref{geo:norm} instead identifies its homogeneous degree spaces
with $H^0(\mathbb{P}^N,L^n)$; for $f=\sum_{n\geq0}f_n$ its full squared norm is
\[
 \|f\|_{\DA}^2=\sum_{n\geq0}c_n\|f_n\|_{\FS,n}^2.
\]
The geometric ambient spaces used in this proof are
$L^2_{0,j}(\mathbb{P}^N,E_q\otimes L^n)$ for finite-rank coefficient bundles
$E_q$.  Intermediate solutions need not be holomorphic.  After the
correction, a holomorphic section of a nonnegative power of $L$ is an
actual homogeneous polynomial, to which \eqref{eq:geo-norm} applies.
We do not use a $\bar\partial$ operator on a space consisting only of holomorphic functions.
\end{remark}

We record one further consequence of the norm correspondence.
For any $n$ the finite spaces $R_n$ with their two norms are
identified by multiplication by $c_n^{1/2}$ when mapping from
Drury--Arveson to Fubini--Study space. Thus if $D:R_{n+1}\to R_n$
is a linear map, its operator norm in the original norm is exactly
$\sqrt{c_n/c_{n+1}}$ times its Fubini--Study operator norm. The same
identity holds for restrictions to subspaces. This factor matters
when we estimate a projected derivative; comparison only within
one degree would not suffice.

\subsection{The geometric task}
We recover the scalar estimate in Section~\ref{sec:completion} from the compatible metrics of Theorem~\ref{thm:D} and the lifting bound of Theorem~\ref{geo:lifting}. In particular, if $\Phi a=F$ is a
holomorphic lift of $F\in I_{n+1}$ through a fixed complete row of ideal
generators, then
\[
 P_{Q_n}\partial_{z_j}F=P_{Q_n}((\partial_{z_j}\Phi)a)
\]
for $n$ beyond the saturation discrepancy. This follows from the ordinary
product rule. To control the right side, we need only construct a
lift and a differentiated-row estimate in one compatible coefficient
metric. The rest of Step 1 constructs the algebraic presentation on which
those estimates will be based.


\section{Complete presentations and their normalized limits}
\label{sec:algebra}

Throughout this section, \(L=\mathcal {O}_{\mathbb{P}^N}(1)\), and
\(R=\mathbb C[Z_0,\ldots,Z_N]\). We assume \(N\geq1\) and fix a nonzero
proper coherent algebraic ideal sheaf \(\mathscr I\subset\mathcal O_{\mathbb{P}^N}\). We carry out all constructions involving Hilbert and Quot
schemes in the algebraic category. After constructing
the presentation and establishing its limiting properties,
we pass to holomorphic sheaves for the analytic arguments. The zero and unit ideals, and the one-variable operator problem, require no
presentation of the kind constructed below.

The local estimates will involve changes of center, changes of scale, and
possibly ill-conditioned constant changes of basis. An arbitrary
finite set of generators is poorly suited to taking these limits. We
instead use complete spaces of sections in degrees chosen uniformly over
proper parameter spaces. This choice controls the dimensions of the
polynomial images in every degree, including at a limiting matrix.

\paragraph{Resolution and the choice of degrees.}
The existence of a partial locally free resolution by sufficiently
positive evaluation spaces follows from the classical generation
and regularity theorems. In the context of essential normality,
\cite{DouglasJabbariTangYu} gives a resolution by analytic Hilbert
modules for monomial ideals; in polynomial membership problems,
\cite{AnderssonWulcan} uses graded free resolutions and projective
Dolbeault complexes. Our choice here serves a different quantitative
purpose. A presentation must still describe the correct image sheaves
when affine dilations and coefficient normalizations degenerate. The
entire proper quotient family, not just the initial fiber, determines
the generating degree. The uniformity proved in this section is the
input to the fixed-scale division theorem.

\subsection{Cohomology, families, and complete section spaces}
\label{alg:background}
We first recall the sheaf-cohomological and parameter-space
terminology. For a sheaf $\mathscr F$, $H^0(\mathbb{P}^N,\mathscr F)$ is its
space of global sections. The higher groups $H^i(\mathbb{P}^N,\mathscr F)$
are the right-derived functors of global sections: one takes a
resolution by injective sheaves, applies global sections, and forms
the cohomology of the resulting complex. In particular a short
exact sequence $0\to\mathscr K\to E\to\mathscr F\to0$ gives a
long exact sequence beginning
\[
 0\to H^0(\mathscr K)\to H^0(E)\to H^0(\mathscr F)
 \to H^1(\mathscr K)\to H^1(E)\to\cdots.
\]
We omit $\mathbb{P}^N(\mathbb{C})$ from this formula. The first obstruction to lifting a
global section of $\mathscr F$ is the indicated $H^1$ class. A
vanishing theorem for that class gives existence of a lift, not a
bound on its norm. The groups of a coherent sheaf on projective
space are finite dimensional and vanish in degrees greater than $N$.
The Euler characteristic of a twist is
$\chi(\mathscr F(t))=\sum_{i=0}^N(-1)^i\dim H^i(\mathbb{P}^N,\mathscr F(t))$.
The long exact sequence gives additivity of $\chi$ in short exact
sequences. Serre's projective generation and vanishing theorems
assert that for sufficiently large $t$ the sheaf $\mathscr F(t)$
is generated by global sections and its higher cohomology vanishes.
Generation means the evaluation morphism
\[
 H^0(\mathbb{P}^N,\mathscr F(t))\otimes\cO\longrightarrow\mathscr F(t)
\]
is surjective on stalks. The resulting Hilbert polynomial is the
polynomial $P_{\mathscr F}(t)=\chi(\mathscr F(t))$; for large $t$
it equals $h^0(\mathbb{P}^N,\mathscr F(t))$. In particular $P_{\cO}(t)$ is
the polynomial $\binom{t+N}{N}$, by counting homogeneous monomials.
We use these standard statements from \cite{SerreFAC,Hartshorne}.

It is necessary to consider entire families of quotient sheaves, not
only a single evaluation row. A family over a parameter scheme
$T$ is described on $\mathbb{P}^N(\mathbb{C})\times T$. We say that a sheaf there is flat
over $T$ if each stalk is a flat module over the local ring at its
parameter. A module is flat when tensor product with it preserves
short exact sequences. Thus flatness of a quotient in a family is
what guarantees that its defining kernel sequence stays exact
on passage to each fiber; it is not a statement about smoothness
or constant pointwise matrix rank. The fiber at $t$ is obtained
by restriction to $\mathbb{P}^N(\mathbb{C})\times\{t\}$, or algebraically by tensoring
with the residue field of $t$ in the parameter direction.

For a fixed bundle $E$ on $\mathbb{P}^N(\mathbb{C})$ and a fixed polynomial $P$, the
Quot functor assigns to a parameter scheme the flat coherent
quotients of the pullback of $E$ whose fibers have Hilbert polynomial
$P$. Two quotients are identified if there is an isomorphism of
quotient sheaves intertwining the two surjections from $E$. Hence
for a single parameter point the datum is equivalently the kernel
as a subsheaf of $E$. Grothendieck's Quot theorem represents this
functor by a parameter space $\Quot_{\mathbb{P}^N}(E,P)$ with a universal quotient
on its product with $\mathbb{P}^N$. For the fixed polynomial and projective
space used here this parameter space is proper over $\C$.
The Hilbert scheme is the special case $E=\cO$: a coherent quotient
of the structure sheaf defines a closed subscheme with the indicated
Hilbert polynomial. A proper map is a separated, finite-type,
universally closed map in algebraic geometry. In particular its image
is closed under any morphism to a separated algebraic variety, and
such a closed subset is also closed in the ordinary complex topology.
We only use the closed-image consequence. For our split source
$E=\cO_{\mathbb{P}^N}(-D)^r$, properness follows from the case $\cO_{\mathbb{P}^N}^r$ by twisting
all quotients by $D$, which replaces $P(t)$ by $P(t+D)$.
We cite \cite{GrothendieckHilbert,StacksQuot} for this theorem.

A Grassmannian $\Gr(a,B)$, for a finite-dimensional complex vector
space $B$, is the space of $a$-dimensional linear subspaces of $B$.
It is compact and projective; in the usual topology a sequence of
subspaces converges if their orthogonal projections (for any fixed
inner product on $B$) converge in norm. A subbundle of the trivial
bundle $B\otimes\cO_T$ of constant rank $a$ defines a morphism from
$T$ to this Grassmannian. In a chart where one $a\times a$ minor
of a basis is invertible, this morphism records the subspace as a
graph of a linear map. The definition does not depend on the basis.
This is precisely how we will retain column spans even when the
changes of coefficient basis do not converge to invertible matrices.

We recall why regularity is needed in addition to properness. For
a family kernel $\mathscr K_t\subset E$, we want to use the full
space $H^0(\mathbb{P}^N(\mathbb{C}),\mathscr K_t(D))$ to define its Grassmannian point.
Neither the dimension of this space nor the surjectivity of its
evaluation map follows from flatness at a small, arbitrary twist.
A common regularity degree, defined next, provides both properties.
Cohomology and base change then guarantees that these spaces form
one subbundle over the parameter space. We use the entire proper
family when choosing $D$; choosing a degree only for the original
kernel would not justify passage to a limit.

\subsection{Regularity and complete evaluation maps}

For a coherent sheaf \(\mathscr F\) on \(\mathbb{P}^N\), we write
\(\mathscr F(t)=\mathscr F\otimes L^t\), and let
\[
 P_{\mathscr F}(t)=\chi(\mathbb{P}^N,\mathscr F(t))
\]
be its Hilbert polynomial. We use the convention
\(R(-D)_\ell=R_{\ell-D}\), with \(R_j=0\) for \(j<0\).

We recall that \(\mathscr F\) is \(a\)-regular if
\[
 H^i(\mathbb{P}^N,\mathscr F(a-i))=0,\qquad 1\leq i\leq N.
\]
We will use the following precise consequences of
Castelnuovo--Mumford regularity.

\begin{lemma}[Complete sections in a regular degree]
\label{alg:regular-sections}
If \(\mathscr F\) is \(a\)-regular, then for every \(t\geq a\):
\begin{enumerate}
\item \(\mathscr F\) is \(t\)-regular, and \(\mathscr F(t)\) is generated
by its global sections;
\item \(H^i(\mathbb{P}^N,\mathscr F(t))=0\) for \(i>0\), so that
\(h^0(\mathbb{P}^N,\mathscr F(t))=P_{\mathscr F}(t)\);
\item multiplication induces a surjection
\[
 H^0(\mathbb{P}^N,L)\otimes H^0(\mathbb{P}^N,\mathscr F(t))
 \longrightarrow H^0(\mathbb{P}^N,\mathscr F(t+1)).
\]
\end{enumerate}
Consequently, for \(b\geq a\) and \(\ell\geq b\), the map
\[
 R_{\ell-b}\otimes H^0(\mathbb{P}^N,\mathscr F(b))
 \longrightarrow H^0(\mathbb{P}^N,\mathscr F(\ell))
\]
is surjective.
\end{lemma}

\begin{proof}
Persistence of regularity, generation, and the multiplication assertion
are the Castelnuovo--Mumford regularity lemmas
\cite[Tags 08A6--08A8]{StacksRegularity}. Persistence up to \(t+i\),
followed by the defining vanishing in cohomological degree \(i\), gives
\(H^i(\mathbb{P}^N,\mathscr F(t))=0\). The equality with the Hilbert polynomial is
then the Euler characteristic formula. Iterating multiplication gives a
surjection from the corresponding tensor power of \(H^0(\mathbb{P}^N,L)\).
Multiplication of linear forms is symmetric, so this map factors through
\(R_{\ell-b}\), proving the last assertion.
\end{proof}

\begin{lemma}[One degree for a whole quotient family]
\label{alg:uniform-regularity}
Fix \(N,r,D\) with \(r>0\), and a polynomial \(P\). There is an integer
\(a\), depending only on these data, such that the kernel of every
quotient
\[
 \mathcal O_{\mathbb{P}^N}(-D)^r\twoheadrightarrow\mathscr F,
 \qquad P_{\mathscr F}=P,
\]
is \(a\)-regular.
\end{lemma}

\begin{proof}
Twisting by \(D\) gives a quotient of \(\mathcal O_{\mathbb{P}^N}^r\) whose Hilbert
polynomial is \(P(t+D)\). The uniform regularity theorem for kernels of
such quotients gives a regularity bound depending only on \(N,r\), and
\(P(t+D)\) \cite[Tag 08AG]{StacksRegularity}. Untwisting adds \(D\) to
that bound.
\end{proof}

We choose the degree in this lemma for the entire quotient family,
including fibers with additional torsion or embedded structure, rather
than just for the sheaf with which we begin. The uniform theorem
\cite[Tag 08AG]{StacksRegularity} does not require the quotients to be
locally free, reduced, or supported on a fixed subvariety. Its
hypotheses are precisely a fixed projective dimension, a fixed
trivial source rank after twisting, and a fixed quotient Hilbert
polynomial. By persistence of regularity we may enlarge the degree
arbitrarily. This permits us to require $D_q\ge D_{q-1}$ at each
stage, so all matrix entries have nonnegative prescribed degree.
The family for the next stage is then determined by this fixed source
bundle and the Hilbert polynomial of its quotient, not by a particular
matrix of a later normalization.

\begin{proposition}[Complete presentation through dimension \(N\)]
\label{alg:complete}
There are integers
\[
 0\leq D_0\leq D_1\leq\cdots\leq D_N
\]
and finite-dimensional spaces \(V_q\) giving a complex
\begin{equation}
 E_N\xrightarrow{d_N}\cdots\xrightarrow{d_2}E_1
 \xrightarrow{d_1}E_0\xrightarrow{d_0=\Phi}\mathcal O_{\mathbb{P}^N},
 \qquad E_q=V_q\otimes L^{-D_q},
 \label{alg:complex}
\end{equation}
with the following properties. Put \(E_{-1}=\mathcal O_{\mathbb{P}^N}\),
\(D_{-1}=0\), \(r_q=\dim V_q\), \(r_{-1}=1\), and
\[
 \mathscr K_0=\mathscr I,\qquad
 \mathscr K_q=\ker d_{q-1}\quad(q\geq1).
\]
Then
\[
 \operatorname{im}d_q=\mathscr K_q,
 \qquad
 V_q=H^0(\mathbb{P}^N,\mathscr K_q(D_q)).
 \tag{\(\ast\)}
\]
Moreover, \(D_0\) is a regular degree for every ideal sheaf with quotient
Hilbert polynomial \(P_{\mathcal O/\mathscr I}\); for \(q\geq1\),
\(D_q\) is a regular degree for every kernel in
\[
 \operatorname{Quot}_{\mathbb{P}^N}(E_{q-1},P_{\mathscr K_{q-1}}).
 \tag{\(\ast\ast\)}
\]
Every nonzero matrix \(d_q\) has columns linearly independent over
constant scalars. Its entries are homogeneous polynomials of degree
\(e_q=D_q-D_{q-1}\).
\end{proposition}

\begin{proof}
We set \(P_0=P_{\mathscr I}\) and
\(P_X=P_{\mathcal O/\mathscr I}\). We apply
Lemma~\ref{alg:uniform-regularity} to the quotients of
\(\mathcal O_{\mathbb{P}^N}\) with polynomial \(P_X\), and choose a common regular
degree \(D_0\geq0\). The complete evaluation map
\[
 H^0(\mathbb{P}^N,\mathscr I(D_0))\otimes L^{-D_0}
 \longrightarrow \mathscr I
\]
is surjective by Lemma~\ref{alg:regular-sections}. Its composition with
\(\mathscr I\hookrightarrow\mathcal O_{\mathbb{P}^N}\) is \(d_0=\Phi\).

Suppose \(E_{q-1}\) and \(d_{q-1}\) have been constructed. Their kernel
fits into an exact sequence
\begin{equation}
 0\longrightarrow\mathscr K_q\longrightarrow E_{q-1}
 \longrightarrow\mathscr K_{q-1}\longrightarrow0.
 \label{alg:kernel-sequence}
\end{equation}
We write \(P_q=P_{\mathscr K_q}\). Additivity gives
\begin{equation}
 P_q(t)=r_{q-1}\binom{t-D_{q-1}+N}{N}-P_{q-1}(t).
 \label{alg:polynomial-recursion}
\end{equation}
The binomial is interpreted as a polynomial in \(t\). The source
\(E_{q-1}\) and the quotient polynomial \(P_{q-1}\) are now fixed.
We choose \(D_q\geq D_{q-1}\) regular for every kernel in
\((\ast\ast)\), and put
\[
 V_q=H^0(\mathbb{P}^N,\mathscr K_q(D_q)),\qquad
 E_q=V_q\otimes L^{-D_q}.
\]
The complete evaluation onto \(\mathscr K_q\), followed by its inclusion
in \(E_{q-1}\), defines \(d_q\). Its image is precisely the preceding
kernel. If a kernel is zero, all subsequent bundles may be taken to be
zero.

We choose a basis of each nonzero \(V_q\). The columns of \(d_q\) are exactly
these basis sections inside
\(H^0(\mathbb{P}^N,E_{q-1}(D_q))\), so a constant linear relation between the
columns is zero only when its coefficients vanish. Since all summands
of the source and target have common twists, the entries have degree
\(D_q-D_{q-1}\). This completes the induction.
\end{proof}

We do not require the presentation to terminate with an injective map \(d_N\). In the analytic correction we use forms only through antiholomorphic degree
\(N\), so exactness at \(E_0,\ldots,E_{N-1}\) is sufficient. Nor do we
assert that the associated graded polynomial complex is exact in small
degrees.

\paragraph{The last map need not be injective.}
For example, on $\Pj^1$ take $\mathscr I=(Z_0)\simeq\cO(-1)$,
$D_0=2$ and $D_1=4$. These degrees meet the uniform regularity
requirements: the ideals of length-one subschemes are isomorphic to
$\cO(-1)$, and the kernels in the next Quot family are torsion-free
rank-one sheaves of degree $-3$, hence isomorphic to $\cO(-3)$.
The complete evaluations are
\[
 d_0=(Z_0^2, Z_0Z_1),\qquad
 d_1=\begin{pmatrix}-Z_1\\ Z_0\end{pmatrix}(Z_0,Z_1).
\]
Here $E_0=\cO(-2)^2$ and $E_1=\cO(-4)^2$. The row $d_0$ has
kernel $\cO(-3)$: the map $(-Z_1,Z_0)^{\mathsf T}$ identifies
that line bundle with the kernel, because the two coordinates
have no common projective zero. The row $(Z_0,Z_1)$ is a surjection
$\cO(-4)^2\to\cO(-3)$ obtained from the complete section space
of $\cO(1)$. Its kernel is $\cO(-5)$ by the same calculation.
Thus $\im d_1=\ker d_0$ but $\ker d_1\simeq\cO(-5)\ne0$.
This is a partial presentation of the permitted kind, not a complex
asserted exact at its left endpoint.

\begin{lemma}[The actual graded image]
\label{alg:evaluation-image}
Let \(J_q\) be the image of the homogeneous polynomial matrix
\[
 R(-D_q)^{r_q}\longrightarrow R(-D_{q-1})^{r_{q-1}}
\]
representing \(d_q\). Then
\begin{equation}
 (J_q)_\ell=
 \begin{cases}
  0,&\ell<D_q,\\
  H^0(\mathbb{P}^N,\mathscr K_q(\ell)),&\ell\geq D_q,
 \end{cases}
 \qquad
 \dim(J_q)_\ell=
 \begin{cases}
  0,&\ell<D_q,\\
  P_q(\ell),&\ell\geq D_q.
 \end{cases}
 \label{alg:image-function}
\end{equation}
The same formulas hold for the complete evaluation of every kernel in
the corresponding parameter family.
\end{lemma}

\begin{proof}
Below \(D_q\), the source homogeneous component vanishes. At degree
\(D_q\), its image is the entire section space by construction. In every
higher degree, Lemma~\ref{alg:regular-sections} identifies the polynomial
multiples of these columns with all sections of
\(\mathscr K_q(\ell)\). Including those sections in
\(H^0(\mathbb{P}^N,E_{q-1}(\ell))\) gives the stated equality. The dimension is its
Hilbert polynomial because all higher cohomology vanishes. The same
regularity bound was imposed on every member of the family.
\end{proof}

\subsection{Sampled Gram forms and compact normalization}

We use standard projective affine charts and their standard line frames.
For a chart \(\alpha\), let \(d_{q,\alpha}(w)\) be the polynomial
coefficient matrix of \(d_q\). Centers \(p\) range in fixed relatively
compact subsets of these charts. We choose \(\varepsilon_0>0\) so that all
balls used in the following definition remain in the corresponding
chart whenever \(0<\varepsilon\leq\varepsilon_0\). Let \(\mu\) be
probability Lebesgue measure on the unit ball \(B_1\subset\mathbb C^N\).

For \(v\in V_q\), we define
\begin{equation}
\begin{split}
 |v|_{G_{q,\alpha,\varepsilon}(p)}^2
 =\int_{(B_1)^{q+1}}
 \bigl|d_{0,\alpha}(p+\varepsilon z_0)\cdots
 d_{q,\alpha}(p+\varepsilon z_q)v\bigr|^2
 \prod_{j=0}^{q}d\mu(z_j).
\end{split}
\label{alg:gram-product}
\end{equation}
The sample points in the different factors are independent. Sampling
all factors at the same point would give zero whenever two consecutive
maps occur, because their composition is zero.

\begin{lemma}[Positivity and recursion]
\label{alg:gram}
The forms \(G_{q,\alpha,\varepsilon}(p)\) are smooth and positive
definite at every nonzero stage. With \(G_{-1,\alpha,\varepsilon}=1\),
they satisfy the frozen-center recursion
\begin{equation}
 G_{q,\alpha,\varepsilon}(p)=
 \int_{B_1}d_{q,\alpha}(p+\varepsilon z)^*
 G_{q-1,\alpha,\varepsilon}(p)
 d_{q,\alpha}(p+\varepsilon z)\,d\mu(z).
 \label{alg:gram-recursion}
\end{equation}
\end{lemma}

\begin{proof}
We integrate first over \(z_0,\ldots,z_{q-1}\) in
\eqref{alg:gram-product}; the result is exactly
\eqref{alg:gram-recursion}. Inductively, if its value on \(v\) is zero,
the nonnegative continuous integrand vanishes on \(B_1\). Positivity of
the preceding form implies
\[
 d_{q,\alpha}(p+\varepsilon z)v=0\qquad(z\in B_1).
\]
Since \(\varepsilon>0\), this is a polynomial identity on a nonempty
open set. It is a constant linear relation between the columns of
\(d_q\), which are independent by Proposition~\ref{alg:complete}.
Therefore \(v=0\). Smoothness follows directly from integration of a
polynomial in \(p,\overline p\), and \(\varepsilon\) over a fixed compact
domain.
\end{proof}

For a fixed center and scale, we abbreviate
\[
 R_q=G_{q,\alpha,\varepsilon}(p)^{1/2},
 \qquad R_{-1}=1,
\]
and define
\begin{equation}
 A_q(\zeta)
 =
 R_{q-1}\,d_{q,\alpha}(p+\varepsilon\zeta)\,R_q^{-1}.
 \label{alg:normalized}
\end{equation}
The \(R_q\) in this formula are constant matrices: the center is frozen
while \(\zeta\) varies.

\begin{thmB}[Complete presentations stable under normalized limits]\label{thm:B}
The complete presentation of Proposition~\ref{alg:complete}, equipped
with the Gram forms \eqref{alg:gram-product}, has the following properties.
The tuples \eqref{alg:normalized}, as the center ranges over the relatively
compact sets of a finite buffered affine atlas and $0<\varepsilon\le
\varepsilon_0$, have compact coefficient closure. Homogenize each map to
its prescribed degree $e_q=D_q-D_{q-1}$. In the coordinates of the corresponding chart (or after the fixed
coordinate permutation identifying that chart with $\{Z_0\neq0\}$), each
simultaneous limit is a complex, exact at $E_0,\ldots,E_{N-1}$, with
image of its first row an ideal sheaf of Hilbert polynomial $P_0$.
Since the atlas is finite, a convergent sequence may be analyzed on a
subsequence in a single chart. For every normalized or limiting
map, its homogeneous polynomial image $J'_q$ satisfies
\[
 \dim (J'_q)_\ell=\begin{cases}0,&\ell<D_q,\\
 P_q(\ell),&\ell\ge D_q.\end{cases}
\]
Here $P_q$ is the polynomial of $\mathscr K_q$ and the degrees in the
image are measured with the prescribed source and target shifts.
\end{thmB}
We establish coefficient compactness, retention of complete sections, exactness, and the all-degree image formula separately.

\begin{proposition}[Compactness of normalized tuples]
\label{alg:compact}
The tuples \((A_0,\ldots,A_N)\) have compact closure in the
finite-dimensional spaces of their polynomial coefficients. At every
nonzero stage,
\begin{equation}
 \int_{B_1}A_q(\zeta)^*A_q(\zeta)\,d\mu(\zeta)=I_{r_q}.
 \label{alg:orthonormal-columns}
\end{equation}
Every simultaneous limiting tuple has linearly independent columns at
each such stage and satisfies
\[
 A_{q-1}^{\infty}(\zeta)A_q^{\infty}(\zeta)=0.
\]
\end{proposition}

\begin{proof}
Substitution into the integral in
\eqref{alg:orthonormal-columns} gives
\[
 R_q^{-1}
 \left(\int_{B_1}d_{q,\alpha}(p+\varepsilon\zeta)^*
 G_{q-1,\alpha,\varepsilon}(p)
 d_{q,\alpha}(p+\varepsilon\zeta)\,d\mu(\zeta)\right)
 R_q^{-1}
 =I.
\]
The columns thus form an orthonormal frame in the fixed
finite-dimensional space of \(\mathbb C^{r_{q-1}}\)-valued polynomials
of degree at most \(e_q\), with its probability-ball inner product.
The space of such frames is compact. We take its product over the finitely
many stages. Orthonormality survives coefficient limits, so no column
dependence appears. The intermediate factors \(R_{q-1}^{-1}R_{q-1}\)
cancel in a product of consecutive normalized matrices. The original
composition is zero, and the resulting polynomial identity survives
every simultaneous limit.
\end{proof}

Compactness and independence of the columns alone would not control
their polynomial multiples. For example, the homogeneous row
\[
 (X^2,\ XY+tY^2)
\]
has two independent columns even at \(t=0\), but its degree-three image
has dimension four for \(t\neq0\) and dimension three at zero. We use complete presentations to rule out this loss of image dimension.

\begin{example}[Why regularity must be uniform over the whole family]
\label{alg:degeneration-test}
On \(\mathbb P^2\), with coordinates \([Z:X:Y]\), consider
\[
 I_s=(X^2,XY,XZ-sY^2),\qquad s\ne0.
\]
On \(Z=1\) its quotient is
\(\mathbb C[x,y]/(x-sy^2,y^3)\), and it has no points on \(Z=0\).
Thus it defines a length-three scheme. Its quadratic generators are
its complete quadratic section space. Indeed, restriction of the six
quadratic monomials to the local quotient has image spanned by
\(1,y,y^2\), so the kernel has dimension three.
Nevertheless, the column-independent limit of this quadratic row is
\((X^2,XY,XZ)\), which generates the line ideal sheaf \((X)\).
Its cubic image has dimension six, whereas the preceding cubic images
have dimension seven.

The omitted direction is explicit:
\[
 Y(XZ-sY^2)-Z(XY)=-sY^3.
\]
Adding \(Y^3\) to the family of generators does not change its fibers
for $s\ne0$. The resulting quotient sheaf is supported in the chart
$Z\ne0$, since $X^2$ and $Y^3$ vanish on its support. Its restriction
there corresponds to
$\mathbb C[s,x,y]/(x-sy^2,y^3)\simeq\mathbb C[s,y]/(y^3)$,
which is free over $\mathbb C[s]$. It is therefore flat over the
parameter line, and its ideal kernel is flat as well, by the exact
sequence with the flat ambient structure sheaf. At $s=0$ the ideal
sheaf is $(X,Y^3)$. Its complete degree-two section space generates
only $(X)$, so the ideal sheaf is not generated in degree two and
cannot be two-regular. A common regularity degree for the \emph{ideal
sheaves} of all length-three subschemes must therefore be at least
three. The quotient sheaves themselves have no higher cohomology in
any twist, since their support is zero dimensional. At degree three the complete section
space of the limiting ideal sheaf has dimension seven and retains the \(Y^3\) direction at the
limit. This example arises from the local dilation
\((x,y)\mapsto(sx,sy)\) applied to \(I_1\); it tests precisely the
kind of rescaling used in \eqref{alg:normalized}. Completeness for
one fiber alone is insufficient. The degree choice in
Proposition~\ref{alg:complete} is what also controls its limits.
\end{example}

\subsection{Complete section spaces form a closed family}

We work first in the chart \(Z_0\neq0\). The affine change of variables
\(\zeta\mapsto p+\varepsilon\zeta\) is induced by the projective
automorphism
\begin{equation}
 g_{p,\varepsilon}[Z_0:Z']
 =
 [Z_0:pZ_0+\varepsilon Z'].
 \label{alg:projective-change}
\end{equation}
Its matrix is invertible for \(\varepsilon>0\). Homogenizing an entry
of \(d_{q,\alpha}(p+\varepsilon\zeta)\) to the prescribed degree
\(e_q\) gives its pullback by \(g_{p,\varepsilon}\), even if its
observed affine degree is smaller. Coordinate permutations give the
same description in the other standard charts.

Pullback identifies \(g_{p,\varepsilon}^*L\) with \(L\) through this
homogeneous substitution. The target matrix \(R_{q-1}\) acts by a
global automorphism of \(E_{q-1}\): all its summands have the same
line twist. The source factor \(R_q^{-1}\) changes the basis without
changing its span. Thus the homogenized columns of \(A_q\) span
\begin{equation}
 H^0(\mathbb{P}^N,\mathscr K_q'(D_q))
 \subset H^0(\mathbb{P}^N,E_{q-1}(D_q)),
 \qquad
 \mathscr K_q'=R_{q-1}g_{p,\varepsilon}^*\mathscr K_q.
 \label{alg:transformed-sections}
\end{equation}
For \(q\geq1\), the quotient \(E_{q-1}/\mathscr K_q'\) has polynomial
\(P_{q-1}\). At \(q=0\), the quotient polynomial is \(P_X\).

\begin{proposition}[Complete sections at every normalized limit]
\label{alg:complete-limit}
Let \((A_0^\infty,\ldots,A_N^\infty)\) be a simultaneous coefficient
limit, and homogenize \(A_q^\infty\) to degree \(e_q\). Its columns
span the complete space
\[
 H^0(\mathbb{P}^N,\mathscr K_q'(D_q))
 \subset H^0(\mathbb{P}^N,E_{q-1}(D_q))
\]
for a kernel \(\mathscr K_q'\) in the parameter family used to choose
\(D_q\). Their evaluation, after untwisting, generates
\(\mathscr K_q'\).
\end{proposition}

\begin{proof}
For \(q=0\), let \(T\) be the Hilbert scheme of quotients of
\(\mathcal O_{\mathbb{P}^N}\) with polynomial \(P_X\); for \(q\geq1\), let
\[
 T=\operatorname{Quot}_{\mathbb{P}^N}(E_{q-1},P_{q-1}).
\]
These are proper parameter spaces
\cite{GrothendieckHilbert,StacksQuot}. We put \(E=E_{q-1}\), and write
their universal quotient sequence on \(\mathbb{P}^N\times T\) as
\begin{equation}
 0\longrightarrow\mathscr U
 \longrightarrow\operatorname{pr}_{\mathbb{P}^N}^*E
 \longrightarrow\mathscr F
 \longrightarrow0.
 \label{alg:universal-kernel}
\end{equation}
The middle term and $\mathscr F$ are flat over $T$: the first is
locally free on a product flat over $T$, and the second is flat by
the defining property of the universal quotient. At a local ring
$A$ of the parameter, apply $-\otimes_A M$ to the sequence, for an
arbitrary $A$-module $M$. Flatness of $\mathscr F$ gives vanishing
of its first Tor group, so the resulting sequence starts with
$0\to\mathscr U\otimes_A M$. Here $\operatorname{Tor}_i^A$ denotes
the derived functors of tensor product; all its positive degrees
vanish when one argument is flat. The long exact Tor sequence in
degree one gives $\operatorname{Tor}_1^A(\mathscr U,M)=0$, since
$\operatorname{Tor}_2^A(\mathscr F,M)=0$ and the middle term is flat.
This proves flatness of $\mathscr U$. Taking $M$ to be the residue
field of the parameter proves exactness on its fiber and identifies
$\mathscr U_t$ with the actual quotient kernel. In particular the
regularity bound applies to this fiber, even at a nonreduced
parameter point.

Let \(\pi:\mathbb{P}^N(\mathbb{C})\times T\to T\). The choice of \(D_q\) gives
\[
 H^i(\mathbb{P}^N,\mathscr U_t(D_q))=0\quad(i>0),\qquad
 h^0(\mathbb{P}^N,\mathscr U_t(D_q))=r_q
\]
for every fiber. We use the following form of cohomology and base change: for a proper
morphism over a locally Noetherian base and a coherent sheaf flat
over the base, vanishing of the fiber cohomology in positive degrees
ensures that the direct image in degree zero is locally free and its
formation commutes with taking fibers. The degree-zero direct image
has sections over a parameter open set given by the sections of the
sheaf on its inverse image. In our situation the morphism is the
projective projection $\pi$, the sheaf $\mathscr U(D_q)$ is coherent
and flat, and the required vanishings hold for every fiber by the
common regularity degree. The cohomology-and-base-change theorem
therefore makes
\[
 \mathscr V=\pi_*(\mathscr U(D_q))
\]
a locally free sheaf of rank \(r_q\), with fiber exactly
\(H^0(\mathbb{P}^N,\mathscr U_t(D_q))\). This is valid over a nonreduced base:
the vanishing of fiber \(H^1\), flatness, and properness are the needed
hypotheses \cite[Theorem~1.1 and Corollary~1.2]{ConradBaseChange};
see also \cite[III, Theorem~12.11]{Hartshorne}.

The inclusion in \eqref{alg:universal-kernel} gives a map
\[
 \mathscr V\longrightarrow
 H^0(\mathbb{P}^N,E(D_q))\otimes\mathcal O_T
\]
whose fiber maps are injective. It is a subbundle injection. Indeed,
at any point choose a nonzero \(r_q\) by \(r_q\) minor of the fiber
matrix. Its determinant is invertible on a neighborhood, and projecting
onto those coordinates and inverting the minor gives a local left
inverse. This argument also applies to nonreduced bases.

The subbundle defines a morphism
\[
 T\longrightarrow
 \operatorname{Gr}\bigl(r_q,H^0(\mathbb{P}^N,E(D_q))\bigr).
\]
Its image is closed by properness. It is therefore closed in the usual
complex topology as well. Equation \eqref{alg:transformed-sections}
places every positive-radius column span in this image.
Proposition~\ref{alg:compact} says that the limiting columns remain
independent. More explicitly, in any fixed inner product on the
finite-dimensional homogeneous section space let $C_k$ be the matrix
with these columns for the $k$th normalized tuple. The orthogonal
projection onto their span is
$C_k(C_k^*C_k)^{-1}C_k^*$. The columns converge to an independent
matrix $C_\infty$; its Gram matrix is invertible. Therefore these
projections converge to the projection onto the limiting span, which
proves convergence in the Grassmannian. It is consequently the complete section space of
one of the family kernels. That kernel is \(D_q\)-regular, so these
sections generate it. The case \(r_q=0\) is the corresponding zero
subbundle and needs no Grassmannian coordinate calculation.
\end{proof}

The Gram recursion \eqref{alg:gram-recursion} and boundedness of the
original polynomial matrices on the buffered sampling regions give
uniform upper bounds on $G_{q-1}$, hence on $R_{q-1}=G_{q-1}^{1/2}$.
However, $R_{q-1}$ can approach a singular matrix, and its inverse
need not be uniformly bounded. We do not require an invertible limit of the changes of basis: each positive-radius normalization defines
a point of the same proper parameter space, and its column span is
retained by the closed Grassmannian image.

\begin{theorem}[Exactness of every limiting complex]
\label{alg:exact-limit}
The fixed-degree homogenizations of every simultaneous normalized limit
satisfy
\[
 \operatorname{im}A_0^\infty=\mathscr I^\infty
\]
for an ideal sheaf with polynomial \(P_0\) and quotient polynomial
\(P_X\), and
\[
 \operatorname{im}A_q^\infty=\ker A_{q-1}^\infty
 \qquad(1\leq q\leq N).
\]
\end{theorem}

\begin{proof}
We obtain the first assertion from Proposition~\ref{alg:complete-limit}: the columns are the complete
regular-degree sections of an ideal in the Hilbert family and generate
that ideal. Suppose inductively that the actual image of
\(A_{q-1}^\infty\) has Hilbert polynomial \(P_{q-1}\).
The exact sequence onto that image gives
\[
 P_{\ker A_{q-1}^\infty}(t)
 =
 r_{q-1}\binom{t-D_{q-1}+N}{N}-P_{q-1}(t)
 =
 P_q(t).
 \]
Independently, Proposition~\ref{alg:complete-limit} identifies
\(\operatorname{im}A_q^\infty\) with a kernel from the corresponding
Quot family, also of polynomial \(P_q\).

Homogenizing the composition identity in Proposition~\ref{alg:compact} to the sum of the two prescribed degrees, we obtain an inclusion
\[
 \operatorname{im}A_q^\infty\subseteq\ker A_{q-1}^\infty.
\]
Their coherent quotient has zero Hilbert polynomial. Such a quotient
must vanish. To see this directly, a sufficiently positive twist of
any nonzero coherent sheaf is generated by global sections and has
zero higher cohomology. A nonzero stalk then forces a nonzero global
section, whereas zero Hilbert polynomial would force its space of
sections to have dimension zero. This is a contradiction.
Thus the inclusion is equality, completing the induction.
\end{proof}

Analytification preserves these exact sequences of coherent algebraic
sheaves \cite[Proposition~10]{SerreGAGA}. Accordingly, for
$1\leq q\leq N$, a holomorphic section of $E_{q-1}$ in the kernel
of $d_{q-1}$ (or of the corresponding limiting map) has a holomorphic
preimage under $d_q$ (respectively the limiting $q$th map) locally at
every point of each affine chart. The first row generates its
corresponding analytic ideal sheaf. These assertions supply membership
in the relevant analytic image sheaves;
the estimates for such preimages are established in the next section.

\begin{theorem}[The Hilbert functions do not change at normalized limits]
\label{alg:hilbert}
For every normalized matrix, including every matrix in a simultaneous
coefficient limit, let \(J_q'\) denote its homogeneous polynomial image
from \(R(-D_q)^{r_q}\) to \(R(-D_{q-1})^{r_{q-1}}\). Then, for every
integer \(\ell\),
\begin{equation}
 \dim(J_q')_\ell=
 \begin{cases}
  0,&\ell<D_q,\\
  P_q(\ell),&\ell\geq D_q.
 \end{cases}
 \label{alg:limit-hilbert}
\end{equation}
\end{theorem}

\begin{proof}
Proposition~\ref{alg:complete-limit} gives the complete degree-\(D_q\)
evaluation of a uniformly regular family kernel at each limit.
Lemma~\ref{alg:evaluation-image} applies to that evaluation and
gives \eqref{alg:limit-hilbert}. The same lemma applies to every
positive-radius member.
\end{proof}

\begin{proof}[Proof of Theorem~\ref{thm:B}]
Proposition~\ref{alg:complete} constructs the presentation with the
required common generating degrees. Lemma~\ref{alg:gram} proves that
its Gram forms are positive. Proposition~\ref{alg:compact} proves
compactness and the composition identities at a simultaneous limit.
Proposition~\ref{alg:complete-limit} identifies the column span of each
limiting map with the complete section space of a kernel in its proper
parameter family. Theorem~\ref{alg:exact-limit} proves exactness, and
Theorem~\ref{alg:hilbert} gives the displayed image formula, completing
the proof of all the assertions.
\end{proof}

Removing the common target shift \(D_{q-1}\) presents the matrix as
\[
 R(-e_q)^{r_q}\longrightarrow R^{r_{q-1}}.
\]
Its degree-\(h\) image is the degree-\(h+D_{q-1}\) image in
\eqref{alg:limit-hilbert}. Thus its image dimension is constant for
every \(h\). Taking the coefficient closure of the tuples in each of
finitely many charts, projecting onto stage \(q\), and taking their
finite union gives a compact matrix family with exactly this property.

\begin{remark}
Keeping the scheme structure is essential in this compactness
argument. For example, the sections \(X(X+\varepsilon Y)\) of
\(\mathcal O_{\mathbb P^1}(2)\), after harmless scalar normalization,
can tend to \(X^2\). The limiting doubled point has the same quotient
Hilbert polynomial as the two preceding points. Replacing the limit
by its reduction would change that polynomial and destroy the
dimension argument.
\end{remark}

\section{Step 2: Uniform holomorphic division}
\label{sec:division}

We have converted rescaled algebraic data into a compact family with a fixed Hilbert function. We now show that these two properties imply a uniform division estimate between two balls. We use a finite Gr\"obner basis only at one reference matrix and then perform a bounded perturbation on a coefficient space containing every degree at once. The leading-term rescaling is related
to the stable division argument of Shalit
\cite[Proposition~3.2, Lemmas~3.4--3.5, and Theorem~3.6]{ShalitDivision}.
Here the parameter uniformity and the treatment of the homogenizing
variable must both be retained.

In this section \(R=\mathbb C[Z_0,\ldots,Z_N]\). We write \(B_s=\{z\in\mathbb C^N:|z|<s\}\), with Lebesgue measure.
For a homogeneous matrix \(A\), we put \(A_{\mathrm{aff}}(z)=A(1,z)\).
A holomorphic vector \(F\) belongs to its analytic image sheaf on
\(B_R\) if, near every point of that ball, it has a holomorphic
representation \(F=A_{\mathrm{aff}}c\). This is stronger than
membership at the center, or membership in the image of the evaluated
matrix at each point.

\begin{thmC}[Uniform holomorphic division between two balls]
\label{div:uniform}
Let \(\mathcal A\) be a compact family, in coefficient topology, of
homogeneous matrices
\[
 A:R(-e)^a\longrightarrow R^b,
 \qquad e\geq0,
\]
with fixed sizes and degree. Suppose that for every integer \(h\),
\begin{equation}
 \dim(\operatorname{im}A)_h
 \quad\text{is independent of }A\in\mathcal A.
 \label{div:fixed-function}
\end{equation}
For every \(s>0\), there are \(R>s\) and \(C<\infty\) such that, for
all \(A\in\mathcal A\), every square-integrable holomorphic \(F\) in
the analytic image sheaf on \(B_R\) has a holomorphic solution
\[
 A_{\mathrm{aff}}u=F|_{B_s},
 \qquad
 \|u\|_{L^2(B_s)}\leq C\|F\|_{L^2(B_R)}.
 \tag{\(\dagger\)}
\]
For each fixed $A$, the solution can be chosen linearly in $F$.
In fact the same affine conclusion holds if \eqref{div:fixed-function}
is weakened to constancy of the eventual Hilbert polynomial of
$\operatorname{im}A$; in this weaker formulation no uniform threshold
for equality of the Hilbert functions is assumed.
\end{thmC}

A nonzero homogeneous submodule of a finite free polynomial module
contains a shifted free rank-one submodule, and therefore has nonzero
Hilbert polynomial. Thus if the common polynomial is zero, every
matrix is zero and the theorem is immediate. We prove it otherwise,
for $N\ge1$. In dimension zero the homogeneous image dimensions
in large degrees equal the rank of a constant matrix after setting
$Z_0=1$. Compactness and constant rank bound its positive singular
values away from zero, which gives the assertion directly.

\paragraph{Relation to parameter-dependent division and scissions.}
We adapt to homogeneous matrices the classical bounded splitting construction of \cite{Galligo,HauserMueller}. In particular,
we may view the perturbation of a reference splitting $\sigma_0$ as
inverting
\[
 U_A=A\sigma_0+(1-A_0\sigma_0)=1+(A-A_0)\sigma_0.
\]
On each degree block $U_A$ carries the reference image injectively into
the new image. Equal dimensions then remove the remainder. We give the
monomial-reconstruction version in detail to verify both the relevant
norms and preservation of homogeneous degree. The selected inverse
need only be complex-linear, not $\mathcal O$-linear. The classical
small-polydisc formulation is compatible with the present large affine
radii: if $r_j$ are homogeneous coefficient radii, write
$\rho_j=r_j/r_0$. For degree-$h$ data and degree-$(h-e)$ coefficients,
\[
 \|F_h\|_{r,\eta}=r_0^h\|F_h\|_{\rho,\eta},\qquad
 \|u_{h-e}\|_r=r_0^{h-e}\|u_{h-e}\|_\rho.
\]
A splitting estimate at fixed $r$ therefore changes only by a factor
$r_0^e$, independent of $h$. In the construction below we choose
$\rho_0=1$ and the affine radii directly. In contrast to
Beatrous' restriction theorem \cite{Beatrous} and its extension in
\cite{DouglasTangYu}, the datum here is already an ambient holomorphic
section in the image sheaf. The loss of domain between $B_R$ and $B_s$
is part of the estimate.

\paragraph{Algebraic and analytic membership.}
We distinguish the polynomial image $J=A R^a\subset R^b$ (with the
specified source shift) from the analytic image sheaf of the affine
matrix. The first requires a polynomial preimage; the second requires
holomorphic preimages on neighborhoods of all points of the domain.
A pointwise preimage of the evaluated matrix is insufficient. For a
homogeneous matrix, $(\im A)_h$ is the image of the finite-dimensional
map $R_{h-e}^a\to R_h^b$. We call $A$ a degree-$e$ matrix when all
its entries are homogeneous polynomials of degree $e$, with zero
entries allowed. Coefficient topology means ordinary Euclidean
convergence of the finite list of coefficients of these entries,
including zero coefficients. In particular multiplication by a fixed
polynomial and homogenization to a fixed degree are continuous maps
in this topology.

We recall the algebra used in polynomial division. A module monomial
is $Z^\alpha e_\nu$, with $e_\nu$ the standard basis vector in the
target. One monomial divides another if their component indices are
equal and every exponent of the first is at most the corresponding
exponent of the second. A monomial submodule is generated by a set
of such monomials and consists of their polynomial multiples. A
monomial order is a total well-order preserved by multiplication
by scalar monomials. For a nonzero polynomial vector, its leading
monomial is the largest monomial with nonzero coefficient (the
coefficient itself is excluded); monic means the leading coefficient
is one. The initial module of $J$ is generated by the leading
monomials of all nonzero elements of $J$. A finite set of elements
of $J$ whose leaders generate this module is a Gr\"obner basis.
The Noetherian property gives finitely many leaders and their
representatives; for a graded $J$, the homogeneous part containing
each leader still has that leader, so we may take the representatives
homogeneous. We need a basis only at a single reference matrix and
will not assume that this Gr\"obner basis persists as the matrix
varies.

\subsection{A reference initial module and coefficient weights}

We fix \(A_0\in\mathcal A\). We order module monomials \(Z^\alpha e_\nu\)
by lexicographic comparison in the variable order
\[
 Z_1,\ldots,Z_N,Z_0,
\]
using a fixed component order to break ties. We choose homogeneous
vectors
\[
 H_i^0=A_0V_i,\qquad 1\leq i\leq k,
\]
with monic leading module monomials \(m_i\) that generate the initial
module of \(\operatorname{im}A_0\). Such a finite collection exists
because monomial submodules of a finite free polynomial module are
finitely generated. The source vector \(V_i\) can be chosen
homogeneous of degree \(\deg m_i-e\), by taking that component of
any polynomial representation of \(H_i^0\).

At each degree \(h\), the module monomials not in this initial module
form a basis of the quotient by \((\operatorname{im}A_0)_h\).
Indeed, in degree $h$ we process terms in descending monomial order.
If the leading term is divisible by $m_i$, we subtract the unique
scalar monomial multiple of $H_i^0$ with that same leading term;
all the new terms are strictly smaller by compatibility of the
order with multiplication. If it is not divisible, we move it to
the remainder. Each operation strictly decreases the largest term
still to be processed. The degree-$h$ block has finitely many
monomials, so the procedure terminates. It expresses the original
polynomial as a combination of the $H_i^0$ plus a combination of
the complementary monomials. Such a nonzero combination cannot lie in the image,
because its leader would then belong to the initial module.
Consequently the number of degree-\(h\) monomials in the initial
module equals \(\dim(\operatorname{im}A_0)_h\).

For positive numbers \(\rho_j\) and \(\eta_\nu\), we define weighted
coefficient norms by
\[
 \|f\|_{\rho,\eta}
 =
 \sum_{\alpha,\nu}|f_{\alpha,\nu}|\rho^\alpha\eta_\nu,
 \qquad
 \rho^\alpha=\prod_{j=0}^{N}\rho_j^{\alpha_j}.
 \]
Scalar series have the norm \(\|\cdot\|_\rho\). Source-vector norms
are sums of the scalar norms of their components.

\begin{lemma}[Dominance of the reference leaders]
\label{div:weights}
For the prescribed \(s>0\), there are weights with
\[
 \rho_0=1,\qquad \rho_j>s\quad(1\leq j\leq N),
 \qquad \eta_\nu>0,
\]
and a coefficient neighborhood \(U\) of \(A_0\), such that
\[
 \|H_i^0-m_i\|_{\rho,\eta}\leq\tfrac14\|m_i\|_{\rho,\eta},
 \qquad
 \|AV_i-H_i^0\|_{\rho,\eta}\leq\tfrac14\|m_i\|_{\rho,\eta}
 \tag{\(\ddagger\)}
\]
for every \(i\) and every \(A\in U\).
\end{lemma}

\begin{proof}
We assign integer scores \(c_\nu\) increasing with component priority.
We set \(w_0=0\), and choose positive integers
\(w_N,w_{N-1},\ldots,w_1\) in that order. Only the finitely many
comparisons between a leader of an \(H_i^0\) and a term in its tail
need to be considered. For one whose first different affine exponent
occurs at \(Z_j\), say leader \((\alpha,\nu)\) and tail
\((\beta,\mu)\), require
\[
 w_j>
 \left|\sum_{\ell=j+1}^{N}w_\ell(\alpha_\ell-\beta_\ell)
       +c_\nu-c_\mu\right|.
\]
This is possible at each step. The leader has a larger exponent of
\(Z_j\), so its weighted exponent is strictly larger. If all affine
exponents agree, homogeneity forces equality of the \(Z_0\)-exponents,
and the component scores decide the comparison.

We put \(\rho_j=T^{w_j}\) and \(\eta_\nu=T^{c_\nu}\).
Each tail-to-leader monomial weight ratio is at most \(T^{-1}\).
For sufficiently large \(T\), the sum of its finitely many
coefficient-weighted ratios is at most \(1/4\) for each \(H_i^0\),
and every affine radius exceeds \(s\). The first inequality follows.
With these weights and the \(V_i\) fixed, continuity of
\(A\mapsto AV_i\) in a finite-dimensional coefficient space supplies
one neighborhood \(U\) for all the second inequalities.
\end{proof}

The leading monomial of \(AV_i\) is allowed to change in this
neighborhood. We use the old monomial partition and control its reconstruction map by a norm perturbation.

\subsection{One bounded inverse for all homogeneous degrees}

Let \(\mathcal W_{\rho,\eta}\) be the Banach space of target formal
series with finite weighted coefficient norm. For a scalar series
\(v\) and target series \(f\), absolute convolution gives
\[
 \|vf\|_{\rho,\eta}\leq\|v\|_\rho\|f\|_{\rho,\eta}.
 \tag{\(\star\)}
\]
Multiplication by a single module monomial gives equality with its
monomial weight.

We partition all target module monomials as follows. A monomial lies
in \(\Delta_i\) if \(m_i\) is the first leader in the chosen list
that divides it; it lies in \(\Delta\) if none does. Divisibility
requires the same target component. Let \(\mathcal Q_i\) consist of
scalar series \(v_i\) for which \(v_im_i\) is supported in
\(\Delta_i\), and let \(\mathcal R_\Delta\) consist of target series
supported in \(\Delta\). These are closed coefficient subspaces.
We give
\[
 \mathcal D=
 \left(\bigoplus_{i=1}^k\mathcal Q_i\right)\oplus\mathcal R_\Delta
\]
the norm
\[
 \|(v,r)\|_{\mathcal D}
 =
 \sum_i\|m_i\|_{\rho,\eta}\|v_i\|_\rho+\|r\|_{\rho,\eta}.
\]
The map
\[
 T_{\mathrm{mon}}(v,r)=\sum_i v_im_i+r
\]
is an isometric bijection onto \(\mathcal W_{\rho,\eta}\), since
its summands occupy disjoint monomial sets.

\begin{lemma}[A uniform, degree-preserving inverse]
\label{div:neumann}
For \(A\in U\), the map
\[
 T_A(v,r)=\sum_i v_iAV_i+r
\]
is invertible and satisfies
\[
 \|T_A^{-1}\|_{\mathcal W_{\rho,\eta}\to\mathcal D}\leq2.
 \]
For homogeneous input of degree \(h\), its inverse has \(r\)
homogeneous of degree \(h\), and \(v_i\) homogeneous of degree
\(h-\deg m_i\).
\end{lemma}

\begin{proof}
Lemma~\ref{div:weights} gives
\(\|AV_i-m_i\|_{\rho,\eta}\leq\frac12\|m_i\|_{\rho,\eta}\).
By \((\star)\),
\[
 \|(T_A-T_{\mathrm{mon}})(v,r)\|_{\rho,\eta}
 \leq\tfrac12\|(v,r)\|_{\mathcal D}.
\]
Thus
\[
 K_A=(T_A-T_{\mathrm{mon}})T_{\mathrm{mon}}^{-1}
 \quad\text{satisfies}\quad \|K_A\|\leq\tfrac12,
\]
and
\[
 T_A^{-1}
 =
 T_{\mathrm{mon}}^{-1}\sum_{\ell=0}^{\infty}(-K_A)^\ell.
\]
The series converges in operator norm and has norm at most two.

Every \(AV_i\) is homogeneous of degree \(\deg m_i\).
Consequently \(T_A\) intertwines the projection onto total degree
\(h\) in the target with the projection onto degrees
\(h-\deg m_i\) and \(h\) in its domain. Its inverse intertwines the
same projections. Those homogeneous blocks are finite dimensional,
so the asserted components are polynomials.
\end{proof}

\begin{lemma}[Vanishing of the remainder on the image]
\label{div:polynomial}
For \(A\in U\cap\mathcal A\) and any degree \(h\) in which
\(\dim(\operatorname{im}A)_h=\dim(\operatorname{im}A_0)_h\), every
homogeneous \(F\in(\operatorname{im}A)_h\) admits a homogeneous source vector
\(u\) of degree \(h-e\) satisfying
\begin{equation}
 Au=F,\qquad
 \|u\|_\rho\leq C_0\|F\|_{\rho,\eta},
 \qquad
 C_0=2\max_i\frac{\|V_i\|_\rho}{\|m_i\|_{\rho,\eta}}.
 \label{div:polynomial-bound}
\end{equation}
The same \(C_0\) works in every degree satisfying this equality; under
\eqref{div:fixed-function} it works in all degrees.
\end{lemma}

\begin{proof}
Let \(\mathcal Q_h\) be the degree-\(h\) part of
\(\bigoplus_i\mathcal Q_i\), where \(v_i\) is assigned shifted
degree \(\deg v_i+\deg m_i\).
Its dimension is the number of degree-\(h\) monomials in the
reference initial module, namely
\(\dim(\operatorname{im}A_0)_h\).
The map \(T_A\) sends \((\mathcal Q_h,0)\) injectively into
\((\operatorname{im}A)_h\), since each \(AV_i\) is a column
combination of \(A\). The stated dimension equality gives
\[
 \dim\mathcal Q_h
 =
 \dim(\operatorname{im}A_0)_h
 =
 \dim(\operatorname{im}A)_h.
\]
The image is therefore the whole target image. The unique inverse
of \(F\) under \(T_A\) has the form \((v,0)\).

We set \(u=\sum_i v_iV_i\). Its summands all have degree \(h-e\), and
\[
 \|u\|_\rho
 \leq\sum_i\|v_i\|_\rho\|V_i\|_\rho
 \leq
 \left(\max_i\frac{\|V_i\|_\rho}{\|m_i\|_{\rho,\eta}}\right)
 \|(v,0)\|_{\mathcal D}.
\]
Lemma~\ref{div:neumann} proves the asserted bound.
\end{proof}

The all-degree hypothesis gives this estimate in every degree. It is
not equivalent to constancy of the eventual Hilbert polynomial. For
example, $A_t=(X^2,tXY,Y^2)$ on $\mathbb P^1$, $|t|\le1$, has the
same image Hilbert polynomial for all $t$: in degrees $h\ge3$ the
image is the entire homogeneous space. But in degree two the image
has dimension three for $t\ne0$ and two for $t=0$. Representing $XY$
in degree two requires coefficient $1/t$ on the middle entry. No
uniform all-degree estimate in the fixed coefficient norm is possible.
Over $\C\{t\}$ the quotient by $(X^2,tXY,Y^2)$ has the nonzero
$t$-torsion class of $XY$ and is not flat. This explains why one
cannot simply invoke Galligo's analytic flatness criterion from
constancy of an eventual projective Hilbert polynomial. The affine
conclusion of C nevertheless holds, as the next argument shows.

\subsection{Dehomogenization without a degree-dependent loss}

We write \(\rho'=(\rho_1,\ldots,\rho_N)\). For an affine polynomial vector \(f(z)=\sum_\beta f_\beta z^\beta\) and \(h\geq\deg f\), we define its degree-\(h\) homogenization by
\[
 f^{[h]}(Z_0,Z')
 =
 \sum_\beta f_\beta Z_0^{h-|\beta|}(Z')^\beta.
\]
Because \(\rho_0=1\),
\begin{equation}
 \|f^{[h]}\|_{\rho,\eta}=\|f\|_{\rho',\eta}.
 \label{div:homogenization-norm}
\end{equation}
There is exactly one \(Z_0\)-power for each affine monomial in a
fixed homogeneous degree, so dehomogenization causes no merging of
distinct coefficients.

\begin{lemma}[Affine polynomial division]
\label{div:affine}
For \(A\in U\cap\mathcal A\), every affine polynomial vector
\(F\in A_{\mathrm{aff}}\mathbb C[z]^a\) has a polynomial preimage
\(u\) satisfying
\[
 A_{\mathrm{aff}}u=F,\qquad
 \|u\|_{\rho'}\leq C_0\|F\|_{\rho',\eta}.
 \]
\end{lemma}

\begin{proof}
We choose a polynomial representation \(F=A_{\mathrm{aff}}c\) and
an integer \(h\geq e+\deg c\). If only the eventual Hilbert
polynomial is fixed, increase \(h\) so that the image dimensions of
\(A\) and \(A_0\) agree in degree \(h\); such a degree exists for this
particular \(A\). We need no common threshold. Both sides of
\[
 F^{[h]}=A\,c^{[h-e]}
\]
are homogeneous of degree \(h\) and dehomogenize to the same
polynomial, so the equality holds. We apply
Lemma~\ref{div:polynomial} and dehomogenize its preimage. Equation
\eqref{div:homogenization-norm}, also with unweighted source
components, gives the bound. The choice of \(h\) introduces no
factor into the estimate.
\end{proof}

The representation degree can be much larger than the degree of
the datum. For example, the row \((z,z+1)\) maps \((-1,1)^T\) to
the constant \(1\), while adding any multiple of
\((z+1,-z)^T\) produces representations of arbitrarily high degree.
The unit weight on the homogenizing variable is what makes the
estimate insensitive to this choice.

\subsection{Holomorphic data and the two-ball estimate}

For a positive vector of radii \(\tau\), we put
\[
 P_\tau=\{z\in\mathbb C^N:|z_j|<\tau_j\text{ for every }j\}.
\]
For the reference neighborhood \(U\), we set
\[
 R_U=2|\rho'|_2+1.
\]
The domains satisfy
\[
 B_s\subset P_{\rho'}\Subset P_{2\rho'}\Subset B_{R_U}.
\]

We recall that a Stein manifold is a complex manifold with enough
global holomorphic functions to separate points and give local
coordinates, and which is holomorphically convex. Holomorphic
convexity means that the hull of a compact set under all global
holomorphic functions is compact. Euclidean balls are Stein: the
coordinate functions separate points and give coordinates, and the
hull of a compact subset of a ball of radius $R$ lies in the smaller
closed ball of radius $r<R$ containing that compact set, by testing
linear functions of norm one. We use Oka coherence (kernels and
images of matrices of holomorphic functions are coherent) and
Cartan's theorem B (positive-degree cohomology of a coherent
analytic sheaf on a Stein manifold vanishes), as in
\cite{GrauertRemmert}. Together with the long exact sequence, these
supply an unestimated holomorphic preimage. The quantitative work
has already been done in the polynomial coefficient estimate;
we do not claim that Cartan's theorem itself supplies a uniform
norm bound.

\begin{lemma}[Holomorphic division in the coefficient norm]
\label{div:holomorphic}
For \(A\in U\cap\mathcal A\), every holomorphic section \(F\) of
the analytic image sheaf on \(B_{R_U}\) has a holomorphic preimage
\(u\) on \(P_{\rho'}\) with
\[
 A_{\mathrm{aff}}u=F,\qquad
 \|u\|_{\rho'}\leq C_0\|F\|_{\rho',\eta}.
 \]
Here the norms are sums of the Taylor coefficients at the origin.
\end{lemma}

\begin{proof}
On \(B_{R_U}\), let \(\mathscr J_A\) be the analytic image and
\(\mathscr Z_A\) the kernel of the holomorphic matrix. Oka
coherence makes \(\mathscr Z_A\) coherent, and
\[
 0\longrightarrow\mathscr Z_A\longrightarrow\mathcal O^a
 \xrightarrow{A_{\mathrm{aff}}}\mathscr J_A
 \longrightarrow0
\]
is exact. The ball is Stein, so Cartan's theorem B gives
\(H^1(B_{R_U},\mathscr Z_A)=0\). We therefore obtain a surjective global-section map: some holomorphic \(c\) on \(B_{R_U}\) satisfies
\(A_{\mathrm{aff}}c=F\). These are the standard coherence and
Stein vanishing theorems; see \cite{GrauertRemmert}.

Let \(c^{(k)}\) be the Taylor polynomial of \(c\) through total
degree \(k\). The closure of \(P_{2\rho'}\) is compact in the
input ball. Cauchy's estimate bounds each component coefficient
by a constant times \((2\rho')^{-\beta}\). Since
\[
 \sum_{\beta\in\mathbb N^N}2^{-|\beta|}=2^N,
\]
the sequence \(c^{(k)}\) converges to \(c\) in
\(\|\cdot\|_{\rho'}\). Multiplication by the fixed polynomial
matrix is continuous in the coefficient norms; hence
\[
 F^{(k)}=A_{\mathrm{aff}}c^{(k)}
 \longrightarrow F
 \quad\text{in }\|\cdot\|_{\rho',\eta}.
\]
These approximants lie in the polynomial image.

We apply Lemma~\ref{div:affine} separately to each \(F^{(k)}\),
obtaining polynomials \(u^{(k)}\) with
\[
 A_{\mathrm{aff}}u^{(k)}=F^{(k)},\qquad
 \|u^{(k)}\|_{\rho'}\leq C_0\|F^{(k)}\|_{\rho',\eta}.
\]
Every coefficient is bounded. By a diagonal argument we pass to a subsequence converging coefficientwise to a series \(u\), and Fatou's lemma for the
weighted sums gives
\[
 \|u\|_{\rho'}\leq C_0\|F\|_{\rho',\eta}.
\]
For \(0<\theta<1\), the tail of total degree greater than \(d\)
on \(P_{\theta\rho'}\) is bounded by
\[
 \theta^{d+1}\sup_k\|u^{(k)}\|_{\rho'}.
\]
The same holds for the limit. Finite coefficient convergence
therefore gives locally uniform convergence on \(P_{\rho'}\).
Its limit is holomorphic and satisfies the required equation.

We have used no quantitative bound on the initial Cartan-B lift \(c\). We use it only to obtain polynomial image approximants, and then construct bounded lifts afresh from their data.
\end{proof}

Approximating the datum itself by Taylor polynomials would not be
valid here. For example, if
\[
 A_{\mathrm{aff}}(z)=z-1,\qquad F(z)=(z-1)e^z,
\]
then the total Taylor polynomial \(T_kF\) has value \(-1/k!\) at
\(z=1\), and so is not in the polynomial image. The approximants
\((z-1)T_k(e^z)\) used by the proof do belong to that image.

\begin{lemma}[Conversion to ball norms]
\label{div:ball}
For every \(A\in U\cap\mathcal A\), the conclusion of
Theorem~\ref{div:uniform} holds with input radius \(R_U\) and
constant
\[
 C_U=
 C_0\,2^N\left(\sum_{\nu=1}^b\eta_\nu\right)
 \left(\frac{\operatorname{Vol}(B_s)}
 {\operatorname{Vol}(B_{1/2})}\right)^{1/2}.
\]
\end{lemma}

\begin{proof}
Cauchy's estimate on \(P_{2\rho'}\) gives
\[
 \|F\|_{\rho',\eta}
 \leq
 2^N\left(\sum_{\nu=1}^b\eta_\nu\right)
 \sup_{P_{2\rho'}}|F|.
\]
For every point in the closure of this polydisc, its ball of
radius \(1/2\) is contained in \(B_{R_U}\). The submean inequality
for the subharmonic function \(|F|^2\) bounds the supremum by
\[
 \operatorname{Vol}(B_{1/2})^{-1/2}\|F\|_{L^2(B_{R_U})}.
\]
On \(B_s\subset P_{\rho'}\), the Euclidean norm of the preimage
from Lemma~\ref{div:holomorphic} is at most its source
coefficient norm. Hence
\[
 \|u\|_{L^2(B_s)}
 \leq\operatorname{Vol}(B_s)^{1/2}\|u\|_{\rho'}.
\]
Combining the three inequalities, we obtain the claim.
\end{proof}

\begin{proof}[Completion of the proof of Theorem~\ref{div:uniform}]
The neighborhoods \(U\) obtained from all reference matrices
cover the compact family \(\mathcal A\). We choose a finite
subcover \(U_1,\ldots,U_M\), and take
\[
 R=\max_j R_{U_j},\qquad C=\max_j C_{U_j}.
\]
For \(A\in U_j\), restrict an admissible datum from \(B_R\) to
\(B_{R_{U_j}}\). Analytic image membership is preserved, and
its norm can only decrease. Lemma~\ref{div:ball} supplies the
solution on the same prescribed output ball \(B_s\).

For a fixed \(A\), multiplication by \(A_{\mathrm{aff}}\) is a
bounded map between the appropriate vector-valued Bergman
spaces on \(B_s\), since its entries are polynomials. Its
kernel is closed. The nonempty fiber over \(F|_{B_s}\) is an
affine translate of this kernel; it has a unique element
perpendicular to the kernel. This element minimizes the norm,
so it obeys the same estimate. Sums and scalar multiples of
these perpendicular solutions solve the corresponding equations
and remain perpendicular to the fixed kernel. Uniqueness proves
linearity on the space of admissible data on \(B_R\).
We neither assert nor need closedness of the range.
\end{proof}

\begin{remark}[What the theorem requires]
The theorem supplies a sufficiently large input radius for a
prescribed output radius. It does not assert the estimate for
every prescribed pair of radii. It also requires analytic image
membership throughout the input ball. For instance, \(1\) is
divisible by \(z-t\) near the origin when \(t\neq0\), but has no
holomorphic quotient on a ball containing \(t\). Finally,
pointwise membership is insufficient: \(z\) belongs to the
image of the scalar matrix \(z^2\) at each point, including
zero, but has no holomorphic quotient near zero.

The compactness hypothesis alone, even with independent limiting columns,
does not ensure uniform division. We consider the homogeneous rows
$(X^2,XY+tY^2)$, $|t|\leq1$, on $\Pj^1$. On $Y=1$ the row is
$(z^2,z+t)$. For $t\ne0$, the constant $1$ has a polynomial
representation, for instance
$1=t^{-2}(z^2-(z-t)(z+t))$. If $0<|t|<s/2$, every holomorphic
solution of $z^2u+(z+t)v=1$ on $B_s$ satisfies $u(-t)=t^{-2}$.
The submean inequality on the disk of radius $s/2$ centered at $-t$
gives $\|u\|_{L^2(B_s)}\geq\sqrt\pi(s/2)|t|^{-2}$ for ordinary
area measure. No uniform two-ball constant exists as $t\to0$.
Indeed the projective image is the unit ideal sheaf for $t\ne0$,
but at $t=0$ it is the ideal sheaf of a point. Its Hilbert polynomial
changes, precisely the kind of degeneration excluded by
Theorem~\ref{thm:B} for the normalized families used here.
\end{remark}

For the normalized presentations, Theorem~\ref{alg:hilbert}
verifies the fixed-Hilbert-function hypothesis, and
Proposition~\ref{alg:compact} supplies compactness. Each original
or limiting presentation is exact on holomorphic stalks at
$E_0,\ldots,E_{N-1}$, and its first row generates its corresponding
ideal sheaf. No injectivity of the last map is asserted. Consequently
both ideal data and the holomorphic cycles in these stages admit the
uniform local division just proved. Proposition~\ref{ana:microscopic-division}
uses exactly these stalks; the final correction form has degree $(0,N)$
and is automatically $\bar\partial$-closed. We next transport this estimate back from the frozen Euclidean coefficient norms to metrics on the actual bundles.


\section{Metrics adapted to the presentation}
\label{ana:metrics}

We now turn the algebraic compactness of the preceding section into
estimates in the original scalar line norm.  Throughout this section
$N\geq1$ and $L=\cO_{\mathbb{P}^N}(1)$ has its Fubini--Study metric.
We fix the complete presentation of Proposition~\ref{alg:complete}, through
stage $N$.  Thus $E_{-1}=\cO_{\mathbb{P}^N}$, $E_q=V_q\otimes L^{-D_q}$,
$d_0=\Phi$, and $d_q:E_q\to E_{q-1}$.  Thus $\im d_0=\cI$ and
$\im d_q=\ker d_{q-1}$ for $1\leq q\leq N$ on holomorphic stalks;
there is no exactness requirement at the last source $E_N$.
The bundles, their bases, their degrees, and the maps are fixed before any
analytic parameter is chosen.  Zero coefficient spaces are allowed, and
all forms and maps at a zero stage are understood to be zero.

We choose a finite collection of buffered affine charts
$V_\alpha\Subset U_\alpha$, with the $V_\alpha$ covering $\mathbb{P}^N$.
We use standard projective line frames on these charts, so the matrices
$d_{q,\alpha}$ are polynomials of fixed degrees.  Every compact set on
which a chart is used has a fixed positive buffer inside that chart.
After decreasing $\varepsilon_0>0$, all sampling balls of radius
$0<\varepsilon\leq\varepsilon_0$ centered in these sets remain in their
charts.  Further fixed multiples of this radius will be accommodated
by decreasing $\varepsilon_0$ when they are introduced.

We recall the local Gram matrices of Lemma~\ref{alg:gram}:
\begin{equation}\label{ana:gram-recall}
 G_{-1,\alpha,\varepsilon}=1,\qquad
 G_{q,\alpha,\varepsilon}(p)
 =\int_{B_1}d_{q,\alpha}(p+\varepsilon z)^*
 G_{q-1,\alpha,\varepsilon}(p)
 d_{q,\alpha}(p+\varepsilon z)\,d\mu(z),
\end{equation}
where $\mu$ is probability volume on the unit ball in $\C^N$.
Equivalently, if
\[
\mathcal P_{q,\alpha,\varepsilon}(p;z_0,\ldots,z_q)
:=d_{0,\alpha}(p+\varepsilon z_0)
\cdots d_{q,\alpha}(p+\varepsilon z_q),
\]
then
\begin{equation}\label{ana:product-gram}
 G_{q,\alpha,\varepsilon}(p)[v]
 =
 \int_{(B_1)^{q+1}}
 \left|
 \mathcal P_{q,\alpha,\varepsilon}
 (p;z_0,\ldots,z_q)v
 \right|^2\,d\mu^{q+1}.
\end{equation}
The sample points are independent. For $q\geq1$, evaluating every map
at the same point would make the products vanish by the complex identities.
The constant independence of the columns of each complete evaluation
instead makes \eqref{ana:gram-recall} positive definite at every
positive radius, as proved in Lemma~\ref{alg:gram}.

\subsection{Hermitian metrics and curvature: conventions}
\label{ana:metric-background}
A smooth Hermitian metric on a complex vector bundle $E$ is a
smoothly varying positive definite Hermitian inner product on its
fibers. In a local holomorphic frame, if a section has coefficient
column $v$, its squared norm is $v^*H(p)v$ for a positive Hermitian
matrix $H(p)$. Under a change of frame $e'=eC$ the coefficient column
becomes $C^{-1}v$ and the metric matrix becomes $C^*HC$, so the
length of the actual fiber vector does not change. We use $H[v]$
for the squared norm and $|v|_H$ for its square root. The operator
norm of a bundle map with matrix $D$ between fibers with metrics
$H$ and $K$ is the norm of $K^{1/2}D H^{-1/2}$ in the Euclidean
inner product. In particular a map is contractive if $D^*KD\le H$
in the order of Hermitian forms. Tensoring both source and target
by the same metrized line preserves this inequality.

We recall the Chern connection, the unique connection on a
holomorphic Hermitian bundle whose $(0,1)$ part is its holomorphic
$\bar\partial$ operator and which preserves the metric. In a
holomorphic frame, writing the connection as $d+\Gamma$, these
properties give $\Gamma^{0,1}=0$ and $\partial H=H\Gamma^{1,0}$.
Hence $\Gamma^{1,0}=H^{-1}\partial H$. The curvature (the square of
the connection) has type $(1,1)$ and in this frame equals
\[
 \Theta(H)=\bar\partial(H^{-1}\partial H).
\]
Indeed the $(2,0)$ part cancels by differentiating $H^{-1}H=I$.
For a line metric $H=e^{-\varphi}$ the formula gives
$\Theta=\partial\bar\partial\varphi$. We use the convention that
$i\Theta$ is positive for the Fubini--Study line $L$, whose potential
is given in Section~\ref{geo:background} and the discussion following
Lemma~\ref{geo:saturation}. Curvature of a tensor product of bundles
is the sum of the curvatures acting on the respective factors; in
particular twisting by $L^h$ adds $h\Theta(L)$ times the identity.
These are the connection and curvature conventions of our estimates
(see \cite{Demailly}).

At a point choose an orthonormal complex tangent frame for the fixed
K\"ahler metric. Denote the curvature coefficients in this frame by
$\Theta_{a\bar b}$, so that $i\Theta(H)=i\sum_{a,b}
\Theta_{a\bar b}\,dz_a\wedge d\bar z_b$ at the point. They satisfy
$\Theta_{a\bar b}^{*H}=\Theta_{b\bar a}$, where $*H$ denotes the
adjoint in the fiber metric. The Nakano quadratic form of the curvature is
\[
  \sum_{a,b=1}^N \langle \Theta_{a\bar b}\xi_a,\xi_b\rangle_H,
 \qquad (\xi_1,\ldots,\xi_N)\in E_p^N.
\]
Nakano positivity means positivity of this quadratic form on all
such tuples. The estimate we need is an absolute upper bound, not
positivity, for the coefficient-bundle curvature. If every
coefficient has operator norm at most $B\varepsilon^{-2}$, the
absolute value of the displayed form is bounded by
$NB\varepsilon^{-2}\sum_a|\xi_a|_H^2$ by Cauchy--Schwarz. We will
then make the tensor product curvature positive by a large power
of $L$. It is essential that all curvature norms here be computed
in the varying fiber metric, rather than in an unrelated fixed
Euclidean metric.

\begin{thmD}[Compatible metrics with controlled curvature]\label{thm:D}
For the presentation in Theorem~\ref{thm:B} there are smooth positive
Hermitian metrics $H_{q,\varepsilon}$ on the nonzero bundles $E_q$, with
$H_{-1,\varepsilon}=1$, for all sufficiently small $\varepsilon>0$.
On each smaller buffered chart $V_\alpha\Subset U_\alpha$,
they are uniformly comparable on the actual fibers to the sampled
Gram metrics $G_{q,\alpha,\varepsilon}$: a common constant $C\geq1$
satisfies
\[
 C^{-1}G_{q,\alpha,\varepsilon}(p)[v]
 \leq H_{q,\varepsilon}(p)[v]
 \leq C G_{q,\alpha,\varepsilon}(p)[v]
 \qquad(p\in V_\alpha,\ v\in(E_q)_p).
\]
We claim no such uniform comparison on the full affine charts.
Each map $d_q$ is contractive from its source metric to its target metric.
With curvature measured in the indicated fiber and Fubini--Study tangent metrics,
\begin{align}
 |\Phi v|&\le |v|_{H_{0,\varepsilon}},\label{eq:geo-contract}\\
 |(\partial_{z_j}\Phi)v|_{L^{-1}}
 &\le C\varepsilon^{-1}|v|_{H_{0,\varepsilon}},\label{eq:geo-differentiate}\\
 \|\Theta(H_{q,\varepsilon})\|_{H_{q,\varepsilon},\omega_{\FS}}
 &\le C\varepsilon^{-2}.\label{eq:geo-curvature}
\end{align}
The constants depend on the fixed presentation and geometry, not on the
scale or the point. In particular the absolute value of each Nakano form
is bounded by a constant times $\varepsilon^{-2}$ times the squared
fiber norm of its arguments.
\end{thmD}

\paragraph{Why the metric varies.}
We take the presentation $Z_0:\cO_{\mathbb P^1}(-1)\to\cO$ of the point
$Z_0=0$. For $F_n=Z_0Z_1^n$ the holomorphic lift is unique and equals
$a_n=Z_1^n$. In the standard Fubini--Study metrics,
\[
 \|a_n\|^2=\frac1{n+1},\qquad
 \|F_n\|^2=\frac1{(n+1)(n+2)}.
\]
Consequently no fixed smooth positive metric on the coefficient line can
make the division estimate uniform: all such metrics are comparable to
the standard one on compact projective space. But the metric
\[
 |v|_{H_\varepsilon}^2=|Z_0v|^2+\varepsilon^2|v|_{\rm std}^2
\]
gives, after twisting by $L^{n+1}$,
\[
 \|a_n\|_{H_\varepsilon\otimes L^{n+1}}^2
 = (1+\varepsilon^2(n+2))\|F_n\|^2.
\]
The scale $\varepsilon\asymp n^{-1/2}$ repairs the division estimate.
Such a change alone is not enough: it must also control derivatives of
the map out of the coefficient bundle and curvature of the new metric.
For this reason we use the sampling construction below.

\paragraph{Comparison with singularly weighted division.}
Skoda's $L^2$ theorem \cite[Theorem~1]{SkodaDivision}, for a row $g$ of
holomorphic functions on a pseudoconvex domain in $\C^N$, uses on its
right side the weight $|g|^{-2(\alpha q+1)}e^{-\psi}$, where
$q=\min(N,r-1)$, $r$ is the number of generators, $\alpha>1$, and $\psi$
is plurisubharmonic, meaning that its restriction to each complex
line is subharmonic (for a smooth function its complex Hessian is
nonnegative). A pseudoconvex domain admits a plurisubharmonic
exhaustion; the balls used in this proof have the strictly
plurisubharmonic exhaustion $-\log(R^2-|z|^2)$. Ordinary ideal membership need not imply finiteness
of this integral. For example, $g=(z_1^2,z_2^2)$, $F=z_1^2$ and $\psi=0$
in $\C^2$ give, on a set of angular directions of positive measure, the
radial integral $\int_0^1 t^{3-4\alpha}\,dt=\infty$. The desired
right side in Theorem E is instead the original scalar norm of $F$.
This does not rule out other approaches using refined $L^2$ division
\cite{SkodaBundles,Varolin,Ji}; it identifies the specific norm problem
that the present metrics address.

\paragraph{Sampling as a finite-jet metric.}
On one affine chart, we expand the polynomial row about $p$ and use the
orthogonality of monomials on the probability unit ball. The exact identity
is
\begin{equation}\label{eq:sampled-jets}
 G_{0,\varepsilon}(p)[v]=
 \sum_{|\alpha|\le D_0}
 \frac{N!\,\varepsilon^{2|\alpha|}}{\alpha!(N+|\alpha|)!}
 |(\partial^\alpha\Phi)(p)v|^2.
\end{equation}
Here $\Phi$ denotes the dehomogenized row in that chart. Indeed the
Taylor coefficient of $z^\alpha$ is $\varepsilon^{|\alpha|}
(\partial^\alpha\Phi)(p)v/\alpha!$, and the unit-ball moment is
$N!\alpha!/(N+|\alpha|)!$. Thus the metric sees all the finite jets at
the sampling scale, even where the value of the row vanishes. The
orders zero and one immediately control evaluation and first derivative
evaluation. At higher stages one must use independent sample points,
since same-point sampling of consecutive differentials vanishes.
For instance the two maps $(x,y)$ and $(-y,x)^T$ on $\C^2$ give
$G_{0,\varepsilon}(0)=\varepsilon^2 I_2/3$ and
$G_{1,\varepsilon}(0)=2\varepsilon^4/9$, by the second moments
of independent points in the unit ball. These small eigenvalues are
not an obstruction: the derivative estimates will be measured relative
to the metrics themselves.

\subsection{Polynomial estimates and relative derivatives}

\begin{lemma}\label{ana:polynomial-estimates}
Fix a bound for the degree and two positive radii.  The $L^2$ norms on
the corresponding balls are equivalent on the space of polynomials of
that degree.  Evaluation at a fixed point in a fixed bounded set,
differentiation as a map between polynomial spaces, and derivative
evaluation are bounded in these norms.  The constants remain valid
for polynomials with values in any finite-dimensional Hilbert space,
independently of its inner product.  The same assertions hold on
products of balls with fixed degrees in each group of variables.
\end{lemma}

\begin{proof}
Each ball integral is a norm on one finite-dimensional scalar polynomial
space: vanishing of the integral forces the polynomial to vanish on an
open set and hence identically.  Norm equivalence and boundedness of
linear maps between finite-dimensional spaces give the scalar claims.
For a vector-valued polynomial we choose an orthonormal basis of the given
target Hilbert space, apply the scalar estimate to every component,
and sum the squared inequalities.  This procedure introduces no
constant depending on the target metric or its smallest eigenvalue.
We apply the same argument successively in the different variable groups
to obtain the assertions on products.
\end{proof}

\begin{lemma}\label{ana:relative-gram}
In a fixed chart we put $G(p)=G_{q,\alpha,\varepsilon}(p)$ and
$|v|_p^2=v^*G(p)v$.  For constant coefficient vectors $v,w$ and affine
coordinate indices $a,b$,
\begin{align}
 |v^*(\partial_aG)w|&\leq C\varepsilon^{-1}|v|_p|w|_p,
 \label{ana:first-relative}\\
 |v^*(\partial_{\bar b}\partial_aG)w|
 &\leq C\varepsilon^{-2}|v|_p|w|_p.
 \label{ana:second-relative}
\end{align}
If two centers in a buffered chart are joined by a segment of length
at most $A\varepsilon$, their Gram norms, expressed in that frame,
are comparable with constants depending only on $A$ and the fixed
presentation.  Furthermore,
\begin{align}
 |d_q(p)v|_{G_{q-1}(p)}&\leq C|v|_{G_q(p)},
 \label{ana:local-map}\\
 |(\partial_a d_q)(p)v|_{G_{q-1}(p)}
 &\leq C\varepsilon^{-1}|v|_{G_q(p)}.
 \label{ana:local-map-derivative}
\end{align}
At $q=0$ the frozen target norm is the usual absolute value in the
local frame of the trivial line bundle.
\end{lemma}

\begin{proof}
Fix the chart $\alpha$ and the scale $\varepsilon$.
We write $\mathcal P_q$ for
$\mathcal P_{q,\alpha,\varepsilon}$ and define the linear map
$F(p)$ by
\[
 F(p)v=\mathcal P_q(p;\cdot)v,
\]
where the target of $F(p)$ is the fixed finite-dimensional
Hilbert space of product polynomials equipped with its
probability $L^2$-norm. By (\ref{ana:product-gram}),  $G(p)=F(p)^*F(p)$ and $F$ depends holomorphically
on $p$.  Each appearance of $p_a$ in the $j$th factor is through
$p_a+\varepsilon z_{j,a}$. We differentiate the matrix product with
$v$ fixed. By the product rule the derivative in $p_a$ is the sum
of the derivatives of its individual factors, leaving all others
unchanged and in their original order. Differentiation in the $j$th
sample variable differentiates only that factor and multiplies its
center derivative by $\varepsilon$. Consequently
\[
 \partial_{p_a}(F(p)v)
 =\varepsilon^{-1}\sum_{j=0}^q
 \partial_{z_{j,a}}(F(p)v).
\]
The degree bounds and the number of groups are fixed, so
Lemma~\ref{ana:polynomial-estimates} gives
\[
 \|(\partial_aF(p))v\|\leq C\varepsilon^{-1}\|F(p)v\|.
\]
Holomorphic dependence gives the exact matrix identities
\[
 \partial_aG=F^*\partial_aF,
 \qquad
 \partial_{\bar b}\partial_aG=(\partial_bF)^*\partial_aF.
\]
For instance the first bilinear form equals
$\langle F(p)v,\partial_a F(p)w\rangle$ in the fixed polynomial
Hilbert space, and its absolute value is at most
$C\varepsilon^{-1}|v|_p|w|_p$. The mixed derivative is estimated by
the product of the two first-derivative bounds. This proves
\eqref{ana:first-relative} and \eqref{ana:second-relative}. Taking
$v=G^{-1/2}x$ and $w=G^{-1/2}y$ for Euclidean unit vectors shows
that these are bounds on $G^{-1/2}G_aG^{-1/2}$ and
$G^{-1/2}G_{a\bar b}G^{-1/2}$. No unweighted inverse bound is
needed.

For a nonzero constant vector $v$, the first estimate and its conjugate
bound the absolute derivative of $\log(v^*G(\gamma(t))v)$ along a
segment by $C|\gamma'(t)|/\varepsilon$.  Integration gives factors
$e^{-CA}$ and $e^{CA}$ for the squared norms.  The buffers ensure
that the segments under consideration stay in the region where the
Gram matrices are defined.

Finally, we give the target of the polynomial
$z\mapsto d_q(p+\varepsilon z)v$ the frozen norm $G_{q-1}(p)$.
Its squared $L^2(B_1)$ norm is $G_q(p)[v]$ by
\eqref{ana:gram-recall}.  Evaluation at zero proves
\eqref{ana:local-map}.  Derivative evaluation, followed by division
by $\varepsilon$, proves \eqref{ana:local-map-derivative}.
The vector-valued assertion of Lemma~\ref{ana:polynomial-estimates}
is what makes these estimates uniform even when a Gram eigenvalue
tends to zero with $\varepsilon$.
\end{proof}

\subsection{Comparison on actual bundle fibers}

We must compare the local quadratic forms on a vector in $E_q$, not
on two unrelated coordinate columns.  We set $D_{-1}=0$.  On the overlap
of two line frames write
\[
 e_\beta(x)=t(x)e_\alpha(x),
\]
where $t$ is holomorphic and nowhere zero.  If $v\in(E_q)_p$ has
column $v_\alpha$ in the first frame, it has column
$v_\beta=t(p)^{D_q}v_\alpha$ in the second.  The matrix transition is
\begin{equation}\label{ana:map-transition}
 d_{j,\beta}(x)=t(x)^{D_{j-1}-D_j}d_{j,\alpha}(x),
\end{equation}
with the corresponding change of affine coordinates understood.

\begin{proposition}\label{ana:fiber-comparison}
After a common decrease of $\varepsilon_0$, the Gram forms defined
in any two charts satisfy, for $p\in V_\alpha\cap V_\beta$,
\[
 C^{-1}G_{q,\alpha,\varepsilon}(p)[v]
 \leq G_{q,\beta,\varepsilon}(p)[v]
 \leq C G_{q,\alpha,\varepsilon}(p)[v]
\]
for every actual fiber vector $v\in(E_q)_p$.  The constant is
uniform in $p$, $v$, and $0<\varepsilon\leq\varepsilon_0$.
\end{proposition}

\begin{proof}
The closures of the smaller charts are compact inside their larger
charts.  There are finitely many chart pairs.  We may therefore require
that both systems of sample points at a point in a smaller overlap
remain in a common buffered overlap.  The transitions, their inverses,
and the coordinate Jacobians are uniformly bounded there.

At independent sample points $x_0,\ldots,x_q$, formula
\eqref{ana:map-transition} and the transformation of the input vector
show that the two product evaluations differ by the scalar
\begin{equation}\label{ana:independent-transition}
 T(p;x_0,\ldots,x_q)
 =t(p)^{D_q}\prod_{j=0}^q t(x_j)^{D_{j-1}-D_j}.
\end{equation}
All factors are evaluated at their indicated points; no telescoping
at one point is being asserted.  Both $|T|$ and $|T|^{-1}$ have
uniform bounds on these overlaps.

We change each sampling variable in the $\beta$ integral to $\alpha$
coordinates.  Its coordinate ball of radius $\varepsilon$ contains
a ball of radius $c\varepsilon$ about the same center and is contained
in a ball of radius $C\varepsilon$, with uniform Jacobian bounds.
The probability-volume normalizations contribute $\varepsilon^{-2N}$
in each variable.  After bounding \eqref{ana:independent-transition}
and the densities, we have upper and lower bounds by integrals of the
same $\alpha$-coordinate product polynomial over products of those
two balls.  Rescale every variable by $\varepsilon$ about the center.
Lemma~\ref{ana:polynomial-estimates} compares the resulting integrals
on the fixed products of balls with the unit-product integral in
\eqref{ana:product-gram}.  This proves both inequalities.
\end{proof}

\paragraph{Why the buffers are needed.}
We consider again $Z_0:\cO_{\Pj^1}(-1)\to\cO$, now on the two full
standard affine charts $U_0=\{Z_0\ne0\}$ and $U_1=\{Z_1\ne0\}$.
For $p$ in their overlap we write $w=Z_0/Z_1$ and take the actual
fiber vector $v=Z_0^{-1}(p)=w^{-1}Z_1^{-1}(p)$. The matrix of
$Z_0$ is $1$ in the $U_0$ frame and $w$ in the $U_1$ frame.
With probability area $\mu$ on the unit disk (the frame labels
$0,1$ below denote the two charts), one has
\[
 \begin{split}
 G_{0,0,\varepsilon}(p)[v]&=1,\\
 G_{0,1,\varepsilon}(p)[v]
 &=|w|^{-2}\int_{|\zeta|<1}|w+\varepsilon\zeta|^2\,d\mu(\zeta)
 =1+\frac{\varepsilon^2}{2|w|^2}.
 \end{split}
\]
For fixed $\varepsilon>0$ the ratio is unbounded as $w\to0$.
Thus uniform comparability to a common global metric on both full
charts would be false. Restriction to buffered charts is essential.

We choose a fixed smooth partition of unity $\{\rho_\alpha\}$,
with $\supp\rho_\alpha\subset V_\alpha$, and define metrics on the
actual bundles by
\begin{equation}\label{ana:patched-metric}
 \widehat G_{q,\varepsilon}
 =\sum_\alpha\rho_\alpha G_{q,\alpha,\varepsilon},
 \qquad \widehat G_{-1,\varepsilon}=1.
\end{equation}
These are smooth positive metrics.  Proposition~\ref{ana:fiber-comparison}
makes the Gram metric of every chart $V_\alpha$ uniformly comparable
to the patched metric on $V_\alpha$: every nonzero summand at $p$
is from a chart $V_\beta$ containing $p$. The same pairwise comparison
applies on $V_\alpha\cap\supp\rho_\beta$, including a point where
$\rho_\beta$ itself is zero, since $\supp\rho_\beta\subset V_\beta$.

\begin{proposition}\label{ana:patched-derivatives}
In any of the fixed smaller affine charts, put
$K=\widehat G_{q,\varepsilon}$.  Then
\begin{align*}
 |v^*(\partial_aK)w|&\leq C\varepsilon^{-1}|v|_K|w|_K,\\
 |v^*(\partial_{\bar b}\partial_aK)w|
 &\leq C\varepsilon^{-2}|v|_K|w|_K.
\end{align*}
With $K_- =\widehat G_{q-1,\varepsilon}$ one also has
\[
 |d_qv|_{K_-}\leq C|v|_K,
 \qquad
 |(\partial_a d_q)v|_{K_-}
 \leq C\varepsilon^{-1}|v|_K.
\]
\end{proposition}

\begin{proof}
In an $\alpha$ frame, a second local metric has matrix
\[
 B_\beta(z)=|t(z)|^{2D_q}
 G_{q,\beta,\varepsilon}(y(z)),
\]
where $y$ is the holomorphic coordinate transition.  Its first and
mixed second derivatives are combinations of the corresponding
derivatives of $G_{q,\beta,\varepsilon}$ and derivatives of the
fixed coordinate and scalar frame factors.  All the latter factors
are bounded on the buffered overlap.  Holomorphicity of $y$ preserves
the mixed derivative type.  Lemma~\ref{ana:relative-gram} therefore
gives the asserted orders, measured initially in $B_\beta$ and then
in $K$ by Proposition~\ref{ana:fiber-comparison}.

The mixed derivative of a weighted summand is explicitly
\[
 (\partial_{\bar b}\partial_a\rho_\beta)B_\beta
 +(\partial_a\rho_\beta)\partial_{\bar b}B_\beta
 +(\partial_{\bar b}\rho_\beta)\partial_aB_\beta
 +\rho_\beta\partial_{\bar b}\partial_aB_\beta.
\]
The derivatives of this fixed partition are bounded.  The support
comparison applies even to the first term when $\rho_\beta$ vanishes.
Since $\varepsilon\leq1$, terms of lower order are bounded by the
same $C\varepsilon^{-2}$.  Summing the finitely many terms proves
the metric estimates.  The map estimates follow from
\eqref{ana:local-map}--\eqref{ana:local-map-derivative} and the same
comparison in source and target.  A differentiation of a line-frame
factor adds a bounded multiple of $d_q$, already covered by the first
map estimate.
\end{proof}

\subsection{Graph metrics and curvature}

We define $H_{-1,\varepsilon}=1$ and, successively on actual fibers,
\begin{equation}\label{ana:graph-definition}
 |v|_{H_{q,\varepsilon}}^2
 =|v|_{\widehat G_{q,\varepsilon}}^2
  +|d_qv|_{H_{q-1,\varepsilon}}^2.
\end{equation}
The identity $d_{q-1}d_q=0$ reduces this recursive expression to
\begin{equation}\label{ana:graph-matrix}
 H_{q,\varepsilon}
 =\widehat G_{q,\varepsilon}
   +d_q^*\widehat G_{q-1,\varepsilon}d_q.
\end{equation}
For $q=0$ the last metric is the scalar metric $1$.

\begin{proposition}\label{ana:graph-curvature}
The metrics $H_{q,\varepsilon}$ are smooth and positive, and
\[
 \widehat G_{q,\varepsilon}\leq H_{q,\varepsilon}
 \leq C\widehat G_{q,\varepsilon}.
\]
Every $d_q:(E_q,H_{q,\varepsilon})\to
(E_{q-1},H_{q-1,\varepsilon})$ is contractive.  The Chern curvature
coefficients, measured in the fiber metric $H_{q,\varepsilon}$ and
unitary Fubini--Study tangent coordinates, have operator norm at most
$B\varepsilon^{-2}$.  In particular their Nakano quadratic form has
absolute value at most
\[
 NB\varepsilon^{-2}\sum_{a=1}^N|\xi_a|_{H_{q,\varepsilon}}^2
\]
on a tuple of fiber vectors $(\xi_1,\ldots,\xi_N)$.
\end{proposition}

\begin{proof}
We obtain the comparison from \eqref{ana:graph-matrix} and Proposition~\ref{ana:patched-derivatives}, and contractivity directly from \eqref{ana:graph-definition}.  It persists after a
common line twist and after tensoring with the same exterior-form
metric.

For the curvature estimate we write, in a fixed holomorphic frame,
$H=J+D^*KD$, with $D=d_q$,
$J=\widehat G_{q,\varepsilon}$, and
$K=\widehat G_{q-1,\varepsilon}$.  We use $K=1$ at stage zero.
If subscripts denote chart derivatives, holomorphicity of $D$ gives
\begin{align*}
 (D^*KD)_a&=D^*K_aD+D^*KD_a,\\
 (D^*KD)_{a\bar b}
 &=D_b^*K_aD+D^*K_{a\bar b}D
   +D_b^*KD_a+D^*K_{\bar b}D_a.
\end{align*}
We have $|Dv|_K\leq|v|_H$ and $|v|_J\leq|v|_H$.
Moreover $|D_av|_K\leq C\varepsilon^{-1}|v|_H$.
Every first-order term above is consequently bounded, when tested
against $v,w$, by $C\varepsilon^{-1}|v|_H|w|_H$; every mixed term
is bounded by $C\varepsilon^{-2}|v|_H|w|_H$.
The same holds for the derivatives of $J$.  Equivalently,
\begin{equation}\label{ana:graph-relative}
 \|H^{-1/2}H_aH^{-1/2}\|\leq C\varepsilon^{-1},
 \qquad
 \|H^{-1/2}H_{a\bar b}H^{-1/2}\|
 \leq C\varepsilon^{-2}.
\end{equation}

The connection matrix is $H^{-1}\partial H$.  In the order
$d\bar z_b\wedge dz_a$, its curvature coefficient is
$\partial_{\bar b}(H^{-1}H_a)$.  Conjugating by $H^{1/2}$ to compute the operator norm in $H$, we obtain
\[
 H^{-1/2}H_{a\bar b}H^{-1/2}
 -(H^{-1/2}H_{\bar b}H^{-1/2})
  (H^{-1/2}H_aH^{-1/2}).
\]
Both terms have norm $O(\varepsilon^{-2})$ by
\eqref{ana:graph-relative}.  Passing to unitary tangent coordinates
costs only a fixed atlas constant.  Taking a maximum over the finitely
many stages gives $B$.  If $R_{a\bar b}$ are the curvature coefficients in such coordinates, then
\[
 \left|\sum_{a,b}\langle R_{a\bar b}\xi_a,\xi_b\rangle_H\right|
 \leq B\varepsilon^{-2}\left(\sum_a|\xi_a|_H\right)^2
 \leq NB\varepsilon^{-2}\sum_a|\xi_a|_H^2.
\]
We use no lower bound on a metric eigenvalue in a fixed Euclidean frame in this calculation.
\end{proof}

\begin{proposition}\label{ana:homogeneous-derivative}
For each original homogeneous coordinate $Z_j$, the differentiated
row defines a bundle map
\[
 \partial_{Z_j}\Phi:E_0\longrightarrow L^{-1},
\]
and
\begin{equation}\label{ana:row-bound}
 |(\partial_{Z_j}\Phi)v|_{L^{-1}}
 \leq C\varepsilon^{-1}|v|_{H_{0,\varepsilon}}.
\end{equation}
The target has its original Fubini--Study line metric.
\end{proposition}

\begin{proof}
We put $D=D_0$.  If $D=0$ the derivative vanishes.  Otherwise, in the
chart $Z_k\neq0$ write $w_\ell=Z_\ell/Z_k$ and
$\Phi(Z)=Z_k^D\Phi_k(w)$.  The row $\partial_{Z_j}\Phi$ is
$Z_k^{D-1}$ times
\[
 P_{j,k}(w)=
 \begin{cases}
  \partial_{w_j}\Phi_k(w),&j\neq k,\\
  D\Phi_k(w)-\displaystyle\sum_{\ell\neq k}
       w_\ell\partial_{w_\ell}\Phi_k(w),&j=k.
 \end{cases}
\]
The coordinates are bounded on the smaller charts.
Proposition~\ref{ana:patched-derivatives}, together with
$\widehat G_0\leq H_0$ and $\varepsilon\leq1$, bounds each of these
rows by $C\varepsilon^{-1}$ in the indicated coefficient metric.
Its homogeneous degree is $D-1$, so its target as a map out of
$V_0\otimes L^{-D}$ is exactly $L^{-1}$.
The squared norm of the local target frame is $|Z_k|_L^{-2}$,
bounded on the chosen smaller chart.  This proves
\eqref{ana:row-bound}.  Equivalently, after tensoring by $L^{n+1}$,
the ratio of output and input line weights is this same fixed factor;
it is independent of $n$.
\end{proof}

\begin{proof}[Proof of Theorem~\ref{thm:D}]
The relative Gram estimates are proved in Lemma~\ref{ana:relative-gram}.
Propositions~\ref{ana:fiber-comparison} and~\ref{ana:patched-derivatives}
show that the chart metrics patch with the same bounds on actual fibers
on each smaller buffered chart $V_\alpha$.
The graph adjustment in Proposition~\ref{ana:graph-curvature} gives
contractivity, equivalence of metrics, and the asserted curvature bound
including the Nakano estimate. Proposition~\ref{ana:homogeneous-derivative}
proves the ambient homogeneous derivative estimate, with its indicated
line-bundle target. This proves all the assertions in the same metrics.
\end{proof}

\section{Microscopic division and a positive-twist estimate}
\label{ana:local-and-global}

We use the preceding metrics for two purposes. They normalize local
division uniformly, and their curvature can be dominated by a large
positive line twist.  We first keep the line power and the sampling
radius independent so that the order in which constants are chosen
remains visible.

\subsection{The Dolbeault solution estimate}

\paragraph{Bundle-valued forms and the basic analytic theorem.}
We recall the analytic framework explicitly. On a complex manifold of
dimension $N$, a smooth $(0,j)$-form with values in a holomorphic
bundle $F$ is, in a holomorphic frame, a sum
$\sum_{|J|=j} f_J\,d\bar z_J$, where $J$ runs over increasing
multi-indices and $f_J$ are smooth coefficient columns. We use the
operator
\[
 \bar\partial_F\Big(\sum_Jf_J\,d\bar z_J\Big)
 =\sum_{a,J}\frac{\partial f_J}{\partial\bar z_a}
       d\bar z_a\wedge d\bar z_J.
\]
It is independent of the holomorphic frame and has square zero.
In degree zero its kernel consists of holomorphic sections. A
holomorphic bundle map commutes with this operator, since its
matrix has zero $\bar\partial$ derivative. The Dolbeault cohomology
is the kernel modulo the image at each degree. We use the term
closed for a form in the kernel, not for a topologically closed
subset or a closed operator range.

The metric on $F$, the exterior form metric induced by the fixed
Fubini--Study metric, and the volume $dV$ give the Hilbert spaces
$L^2_{0,j}(\mathbb{P}^N,F)$. The maximal $L^2$ domain of $\bar\partial$ is the
space of forms whose distributional derivative belongs to the next
$L^2$ space. We use its Hilbert adjoint $\bar\partial^*$, which
therefore depends on the metric. The Dolbeault Laplacian in degree
$j$ is $\Box_j=\bar\partial\bar\partial^*+
\bar\partial^*\bar\partial$. On the compact manifold with smooth
positive bundle metrics this is an elliptic nonnegative self-adjoint
operator (with the usual closed quadratic form realization).
Ellipticity here refers to invertibility of the principal symbol at
nonzero cotangent vectors; it gives smooth solutions for smooth
right-hand sides by elliptic regularity.

We use the classical Bochner--Kodaira--Nakano identity on a K\"ahler
manifold. A K\"ahler metric is a Hermitian metric whose associated
real $(1,1)$ form $\omega$ is closed. Let $\Lambda$ be the pointwise
adjoint of exterior multiplication by $\omega$. For a Hermitian
holomorphic bundle $G$ and smooth $G$-valued forms of type $(p,j)$
on a compact K\"ahler manifold without boundary, the identity is
\[
 \|\bar\partial v\|^2+\|\bar\partial^*v\|^2
 =\|D'v\|^2+\|(D')^*v\|^2
     +\int_{\mathbb{P}^N}\langle[i\Theta(G),\Lambda]v,v\rangle\,dV,
\]
where $D'$ is the $(1,0)$ part of its Chern exterior covariant
derivative. The scalar normalization of the volume has no effect.
In top holomorphic degree $p=N$, the curvature term is the sum of
the Nakano forms on the coefficient tuples obtained by inserting
one antiholomorphic index. In particular a scalar line with
positive eigenvalues $\lambda_a$ contributes on a coefficient
with antiholomorphic index set $J$ the sum $\sum_{a\in J}\lambda_a$.
We use this identity and standard compact elliptic theory in the
form explained in the next proof, with the hypotheses checked
there; see \cite{Demailly,DemaillyL2}. No positivity of the sampled
coefficient metric by itself is assumed.

The canonical bundle $K_{\mathbb{P}^N}$ is the top exterior power of the
holomorphic cotangent bundle, locally generated by
$dz_1\wedge\cdots\wedge dz_N$. The metric on the base induces its
metric and that of the dual $K_{\mathbb{P}^N}^{-1}$. We use the canonical
identification $\Lambda^{0,j}\otimes F\simeq
\Lambda^{N,j}\otimes F\otimes K_{\mathbb{P}^N}^{-1}$ to apply the identity in
top holomorphic degree. The curvature of $K_{\mathbb{P}^N}^{-1}$ is fixed and
bounded on compact $\mathbb{P}^N$. We do not identify the $(0,j)$ and $(N,j)$
curvature expressions without this twist.

\begin{proposition}\label{ana:positive-gap}
There are $c_L>0$ and $C_0,C_1\geq0$, independent of the coefficient
stage, such that the following holds.  Let $h\geq0$ be an integer,
$0<\varepsilon\leq\varepsilon_0$, and
\[
 \Delta(h,\varepsilon)
 =c_Lh-C_0-C_1\varepsilon^{-2}>0.
\]
For every distributionally closed $L^2$ form $\zeta$ of type $(0,j)$,
$1\leq j\leq N$, with values in $E_q\otimes L^h$, there is a
unique solution perpendicular to $\ker\bar\partial$ in degree $j-1$,
and it satisfies
\begin{equation}\label{ana:solver}
 \bar\partial u=\zeta,
 \qquad
 \|u\|\leq\Delta(h,\varepsilon)^{-1/2}\|\zeta\|.
\end{equation}
All norms use $H_{q,\varepsilon}$, the original metric on $L^h$,
and the Fubini--Study form metric and volume.  At fixed parameters
the solution operator is linear and sends smooth closed forms to
smooth forms.
\end{proposition}

\begin{proof}
We put $F=E_q\otimes L^h$.  The canonical identification of $F$-valued
$(0,j)$-forms with $(F\otimes K_{\mathbb{P}^N}^{-1})$-valued $(N,j)$-forms is
holomorphic and isometric when the canonical line has the induced
metric.  Locally one tensors a form with a nonvanishing holomorphic
volume form and its dual; their frame changes and metric factors
cancel.  We choose the order of the factors so that this identification
commutes with $\bar\partial$.  Its Hilbert adjoints then correspond
as well.

The Bochner--Kodaira--Nakano identity in top holomorphic degree gives
\[
 \|\bar\partial v\|^2+\|\bar\partial^*v\|^2
 \geq\int_{\mathbb{P}^N}
 \langle[i\Theta(F\otimes K_{\mathbb{P}^N}^{-1}),\Lambda]v,v\rangle\,dV.
\]
We use the sign convention in which the curvature eigenvalues of
the positive line $L$ are positive.  This identity and its $L^2$
consequences are recalled in \cite{Demailly,DemaillyL2}.
Let $c_L$ be a positive lower bound for the eigenvalues of $i\Theta(L)$.
On a top-holomorphic-degree form with antiholomorphic index set $J$,
the line contribution is
$h\sum_{a\in J}\lambda_a\geq hjc_L$.
The fixed curvature of $K_{\mathbb{P}^N}^{-1}$ contributes a bounded operator,
which we allow to cost $C_0$.

For the coefficient curvature, we fix an increasing $(j-1)$-index $S$
and denote by $v_{aS}$ the signed coefficient obtained by inserting
$a$ before $S$, with value zero when an index repeats.  The
top-degree curvature formula is the sum, over $S$, of the Nakano
forms of $(v_{aS})_{a=1}^N$.  By
Proposition~\ref{ana:graph-curvature} its absolute value is bounded
by
\[
 NB\varepsilon^{-2}\sum_{S,a}|v_{aS}|^2
 =NjB\varepsilon^{-2}|v|^2.
\]
Each increasing $j$-index occurs exactly $j$ times in the middle
sum.  We can therefore take $C_1=N^2B$, since $j\leq N$, and obtain
\begin{equation}\label{ana:basic-estimate}
 \|\bar\partial v\|^2+\|\bar\partial^*v\|^2
 \geq\Delta(h,\varepsilon)\|v\|^2.
\end{equation}

For completeness, on this compact manifold we may obtain the
solution directly from the Dolbeault Laplacian.  The smooth positive
metrics at fixed $(h,\varepsilon)$ make that operator elliptic and
self-adjoint.  Estimate \eqref{ana:basic-estimate} removes its kernel
in positive degree and bounds its inverse $G_j$ by $\Delta^{-1}$.
For closed data the Hodge identities give
\[
 u=\bar\partial^*G_j\zeta,
 \qquad \bar\partial u=\zeta.
\]
For smooth data these follow from commutation with the Laplacian;
$\bar\partial G_j\zeta$ is harmonic and hence zero in the next
degree, or is zero by degree at $j=N$.  The same identities extend
to $L^2$ closed data by the spectral calculus, or by smooth heat
regularization, which commutes with $\bar\partial$.  Moreover
\[
 \|u\|^2
 =\langle\zeta,G_j\zeta\rangle
 \leq\Delta^{-1}\|\zeta\|^2.
\]
The image of $\bar\partial^*$ is perpendicular to the closed forms.
Any two solutions differ by a closed form, proving minimality and
uniqueness.  These also imply linearity.  Elliptic regularity gives
smoothness for smooth inputs.
\end{proof}

\subsection{Division on shrinking balls with the original line weight}

We write $B(p,r)$ for a geodesic ball in $\mathbb{P}^N$.  For a section of
$E_q\otimes L^h$ over $U$, let $\|\cdot\|_{q,h;U}$ denote its
$L^2$ norm in $H_{q,\varepsilon}\otimes L^h$.  At $q=-1$ this is
the original scalar line norm.

\begin{proposition}\label{ana:microscopic-division}
For every $s>0$ there are $\Lambda(s)>s$ and
$\varepsilon_0(s)>0$ with the following property.  Fix $K\geq0$.
Suppose $h\geq0$ is an integer,
$0<\varepsilon\leq\varepsilon_0(s)$, and $h\varepsilon^2\leq K$.
Let $b$ be a square-integrable holomorphic section of
$E_{q-1}\otimes L^h$ on $B(p,\Lambda(s)\varepsilon)$.
For $q=0$ assume $b$ is a section of $\cI\otimes L^h$ there;
for $1\leq q\leq N$ assume $d_{q-1}b=0$ there.
There is a holomorphic section $a$ on $B(p,s\varepsilon)$ such that
\begin{equation}\label{ana:micro-bound}
 d_qa=b,
 \qquad
 \|a\|_{q,h;B(p,s\varepsilon)}
 \leq C_s e^{C\Lambda(s)^2K}
 \|b\|_{q-1,h;B(p,\Lambda(s)\varepsilon)}.
\end{equation}
The choice can be made linear in $b$ for each fixed set of parameters.
The constants are uniform in the center, the radius, and the line power.
\end{proposition}

\begin{proof}
If $E_q=0$, the image condition implies that $b=0$,
and we take the lift to be zero.
We may thus assume that $E_q\ne0$.

We first choose the coordinate chart and the radii for local
division. We choose $\ell>0$ so that every geodesic ball
$B(p,2\ell)$ is contained in some smaller chart $V_\alpha$;
such an $\ell$ exists by the Lebesgue number lemma applied
to the finite smaller-chart cover.
The compact buffers in $U_\alpha$ and the finitely many
smooth coordinate maps give, after decreasing $\ell$ if
necessary, a common constant $H\geq1$ such that
for $0<r\leq\ell$ in such a chart,
\[
 B_{\rm coord}(p,r/H)\subset B(p,r)
 \subset B_{\rm coord}(p,Hr),
\]
where $B_{\rm coord}$ denotes a Euclidean coordinate ball.
These inclusions follow from the uniform upper and lower
bounds for the Fubini--Study metric in the buffered
coordinates. We take $\ell$ sufficiently small that the
relevant straight segments and minimizing geodesics remain
in $U_\alpha$.

Proposition~\ref{alg:compact} and
Theorem~\ref{alg:hilbert} supply compact families of
normalized homogeneous matrices with fixed image dimensions
in every degree, after the common target shift is removed.
For the output radius factor $Hs$, let $\lambda_q(s)>Hs$
be the input radius factor supplied by
Theorem~\ref{div:uniform} for the family at stage $q$.
We set
\[
 \Lambda(s)=H\max_{0\leq q\leq N}\lambda_q(s)>s,
\]
assigning any factor larger than $Hs$ at a zero stage.
We decrease $\varepsilon_0(s)$ so that
$\Lambda(s)\varepsilon<\ell$ whenever
$0<\varepsilon<\varepsilon_0(s)$.

For each $p$, we choose a smaller chart $V_\alpha$
containing $B(p,2\ell)$. The preceding inclusions give
\[
 \begin{aligned}
 B(p,s\varepsilon)
 &\subset B_{\rm coord}(p,Hs\varepsilon)\\
 &\subset B_{\rm coord}(p,\lambda_q(s)\varepsilon)\\
 &\subset B(p,\Lambda(s)\varepsilon)
 \subset V_\alpha.
 \end{aligned}
\]
In particular, all the balls on which we perform division
lie in the domain of the uniform metric comparison in
Theorem~\ref{thm:D}.

We next normalize the polynomial matrix in the chosen chart.
We evaluate the local Gram metrics at the center $p$ and
keep them fixed throughout the following division argument.
We put
\[
 R_j=G_{j,\alpha,\varepsilon}(p)^{1/2},
 \qquad R_{-1}=1.
\]
Here we use $p$ also to denote its coordinate representative
in the chosen chart.
In the rescaled coordinate $\zeta=(x-p)/\varepsilon$,
the normalized matrix is
\[
 A_q(\zeta)
 =R_{q-1}d_{q,\alpha}(p+\varepsilon\zeta)R_q^{-1}.
\]
This is the affine representative of a member of the compact
family used above to choose $\lambda_q(s)$.

We now verify the image-sheaf condition on the entire input
ball. At $q=0$, this is the assumed ideal-sheaf membership.
For $1\leq q\leq N$, it follows from exactness of the actual
presentation at $E_{q-1}$, since $b$ is a holomorphic cycle.
Analytification preserves this exactness.
We emphasize that we use membership in the analytic image
sheaf, not merely pointwise membership in the image of an
evaluated matrix.

Before applying Theorem~\ref{div:uniform}, we choose a
holomorphic frame of the common factor $L^h$ in which its
Hermitian weight is uniformly controlled on the input ball.
We write
\[
 |e(x)|_L^2=e^{-\varphi(x)}
\]
in a local holomorphic frame $e$ of $L$, and replace $e$
by the centered holomorphic frame
\[
 e'(x)=
 \exp\!\left(
 \frac12\varphi(p)
 +\sum_a\partial_a\varphi(p)(x_a-p_a)
 \right)e(x).
\]
Here $\partial_a$ denotes differentiation with respect to
the holomorphic coordinate $x_a$.
The metric potential in this new frame is
\[
 \varphi'(x)
 =\varphi(x)-\varphi(p)
 -2\operatorname{Re}
 \sum_a\partial_a\varphi(p)(x_a-p_a).
\]
Since $\varphi$ is real-valued, the expression subtracted
from $\varphi(x)$ is precisely its constant and linear
Taylor polynomial at $p$.
In particular,
\[
 \varphi'(p)=0,\qquad d\varphi'(p)=0.
\]
The real second derivatives of the original potentials
are uniformly bounded in the fixed buffered atlas.
Taylor's formula with integral remainder therefore gives
\[
 |\varphi'(x)|\leq C|x-p|^2
\]
on the relevant coordinate balls, with $C$ independent of
the center and the scale.

We use $(e')^h$ as the frame of the common factor $L^h$,
while keeping the coefficient-bundle frames fixed.
The polynomial matrices $d_{q,\alpha}$ are unchanged by
this operation, because the same additional line frame
occurs in both their sources and their targets.
Since the coordinate input radius is bounded by
$\Lambda(s)\varepsilon$, we obtain
\[
 e^{-C\Lambda(s)^2h\varepsilon^2}
 \leq |(e'(x))^h|_{L^h}^2
 \leq e^{C\Lambda(s)^2h\varepsilon^2}
\]
throughout the input ball.
Thus both the squared line weight and its reciprocal
are uniformly controlled whenever $h\varepsilon^2$
is bounded.

In this centered holomorphic frame of the common factor
$L^h$, we write $b(x)$ for the local coefficient vector
of the datum and define its normalized representative by
\[
 \widetilde b(\zeta)
 =R_{q-1}b(p+\varepsilon\zeta).
\]
The preceding image-sheaf condition and the invertibility
of the constant matrices $R_j$ imply that
$\widetilde b$ belongs to the analytic image sheaf of
$A_q$ throughout $B_{\lambda_q(s)}$.

We apply Theorem~\ref{div:uniform} to obtain a holomorphic
coefficient vector $\widetilde a$ on $B_{Hs}$ satisfying
\[
 A_q\widetilde a=\widetilde b,
 \qquad
 \|\widetilde a\|_{L^2(B_{Hs})}
 \leq C_s
 \|\widetilde b\|_{L^2(B_{\lambda_q(s)})},
\]
where $C_s$ is independent of $p$, $\varepsilon$, and $h$.
We then define the local coefficient vector
\[
 a(x)=R_q^{-1}
 \widetilde a\bigl((x-p)/\varepsilon\bigr)
\]
in the same centered line frame.
By construction,
\[
 \begin{aligned}
 R_{q-1}d_{q,\alpha}(x)a(x)
 &=A_q\bigl((x-p)/\varepsilon\bigr)
   \widetilde a\bigl((x-p)/\varepsilon\bigr)\\
 &=R_{q-1}b(x).
 \end{aligned}
\]
Multiplying by $R_{q-1}^{-1}$ gives $d_qa=b$.
We restrict $a$ to the desired geodesic ball
$B(p,s\varepsilon)$.

It remains to compare the normalized estimate with the
original bundle norms.
We temporarily omit the common line weight and use
coordinate Lebesgue volume.
For any coefficient section $f$ and coordinate radius
$t\varepsilon$, the change of variables
$x=p+\varepsilon\zeta$ gives
\[
 \int_{B_{\rm coord}(p,t\varepsilon)}
 |f(x)|_{G_j(p)}^2\,dx
 =
 \varepsilon^{2N}\int_{B_t}
 |R_jf(p+\varepsilon\zeta)|^2\,d\zeta.
\]
The factor $\varepsilon^{2N}$ is the real Jacobian of
this change of variables in $\mathbb C^N$.
Consequently, the normalized division estimate gives
\[
 \|a\|_{L^2(B_{\rm coord}(p,Hs\varepsilon),G_q(p))}
 \leq C_s
 \|b\|_{L^2(B_{\rm coord}(p,\lambda_q(s)\varepsilon),
                  G_{q-1}(p))}.
\]
In particular, the factor $\varepsilon^N$ occurring in
each unsquared norm cancels.

We now restore the Fubini--Study volume and the moving
coefficient metrics.
The chart volume densities are uniformly bounded above
and below.
Lemma~\ref{ana:relative-gram} compares the frozen and
moving Gram norms on these fixed multiples of
$\varepsilon$, while
Propositions~\ref{ana:fiber-comparison} and
\ref{ana:graph-curvature} compare the moving Gram norms
with the metrics $H_q$ and $H_{q-1}$.
All the resulting comparison constants depend only on
the fixed presentation, the buffered atlas, and the
selected radius factors.

Finally, restoring the common line weight contributes
at most a factor
$e^{C\Lambda(s)^2h\varepsilon^2}$ to the division estimate,
after increasing the fixed constant $C$.
Together with the nested-ball inclusions, we obtain
\[
 \|a\|_{L^2(B(p,s\varepsilon),H_q\otimes L^h)}
 \leq C_s e^{C\Lambda(s)^2h\varepsilon^2}
 \|b\|_{L^2(B(p,\Lambda(s)\varepsilon),
                         H_{q-1}\otimes L^h)}.
\]
This is precisely \eqref{ana:micro-bound}.
Apart from the displayed exponential factor, all constants
are independent of $p$, $\varepsilon$, and $h$.

For each fixed normalized matrix, we use
Theorem~\ref{div:uniform} to choose the minimum-norm
holomorphic primitive on its coordinate output ball.
This choice is linear on admissible data.
Normalization, changes of line frames, and restrictions
all preserve linearity.
No continuity of this choice as $p$ or $\varepsilon$
varies is required.
\end{proof}

\section{Holomorphic overlap data and the global lift}
\label{ana:global-section}

We use local division only on holomorphic sections.  A shrinking
cover converts these holomorphic overlap data into smooth Dolbeault
forms.  The loss from each differentiated cutoff is then recovered
by one use of the positive-twist solution estimate.

\begin{thmE}[Bounded projective lifting]\label{geo:lifting}
Let $\cI\subset\cO_{\mathbb{P}^N}$ be a nonzero proper coherent ideal sheaf on
$\Pj^N$, $N\ge1$.  Let $D=D_0$ be chosen as in Theorem~\ref{thm:B}, with complete
space $V=H^0(\mathbb{P}^N,\cI\otimes L^D)$ and evaluation map
\[
 \Phi:E_0=V\otimes L^{-D}\longrightarrow\cO,
 \qquad\im\Phi=\cI.
\]
Use the metrics of Theorem~\ref{thm:D}, in particular
\eqref{eq:geo-contract}--\eqref{eq:geo-curvature}.
There are fixed constants $R_*>0$, $n_0$, and $C'$ such that for
$n\ge n_0$, $\varepsilon=R_*/\sqrt n$, and $h\in\{n,n+1\}$,
every $F\in H^0(\mathbb{P}^N,\cI\otimes L^h)$ has a global holomorphic lift
$a\in H^0(\mathbb{P}^N,E_0\otimes L^h)$ satisfying
\begin{equation}\label{eq:geo-lift-bound}
 \Phi a=F,\qquad
 \|a\|_{H_{0,\varepsilon}\otimes L^h,L^2}
       \le C'\|F\|_{L^h,L^2}.
\end{equation}
The right-hand norm is the original scalar Fubini--Study norm.
The lift can be chosen linearly in $F$.  All constants and the
presentation are independent of $n,h,F$.  The same coefficient metric
serves the two adjacent powers. The unit ideal sheaf has the same lifting
conclusion with the row $1$ and the constant metric.
\end{thmE}

The derivatives in \eqref{eq:geo-differentiate} are derivatives of the
homogeneous polynomial entries of $\Phi$ in the fixed ambient
coordinates.  Their degree is $D-1$, so they define bundle maps
$E_0\to L^{-1}$, as indicated.  For $D=0$ they vanish.

\paragraph{Historical relation of the descent.}
We use a finite correction through a complex, a version of the
Koszul--Dolbeault procedure in \cite{KelleherTaylor,HormanderGenerators}
and of the projective construction in \cite[Section~6]{AnderssonWulcan}.
In the latter the residue-current identities supply data for successive
Dolbeault equations. Here the cochains are holomorphic overlap lifts,
obtained with uniform estimates from Theorem C before cutoffs are
introduced. We prove the precise norm bounds and closedness at every
stage; no vanishing theorem alone would imply the quantitative conclusion.

\paragraph{The two differentials and the correction.}
We recall the elementary \v{C}ech terminology. For a cover with
ordered indices, a $q$-cochain of sections assigns a section to each
$(q+1)$-fold intersection; it is extended to arbitrary ordered
tuples by alternation. Its coboundary $\delta$ on a tuple of length
$q+2$ is the alternating sum of its restrictions to the faces
obtained by deleting one index. In $\delta^2$ each pair of deleted
indices occurs twice with opposite signs; thus $\delta^2=0$. Our
cochains will be defined on larger buffered balls, not merely the
intersections, so each such restriction must be justified. A
syzygy is a relation among the columns of a presentation. At a
point of our sheaf complex it is a germ in the kernel of the
preceding matrix; exactness says that it is locally the image of
the next matrix. In the construction we first lift the holomorphic
\v{C}ech differences using these syzygy maps and only afterward
insert smooth cutoffs. This order avoids trying to divide
nonholomorphic forms by Theorem C.

For clarity, the first two levels of the construction on a surface
of complex dimension two would be
\[
 \Phi c_{0,i}=F,\quad
 d_1c_{1,ij}=c_{0,j}-c_{0,i},\quad
 d_2c_{2,ijk}=c_{1,jk}-c_{1,ik}+c_{1,ij}.
\]
The second difference is killed by $d_1$, because the three
first differences cancel. After using a partition of unity we
will have smooth forms $\eta_0,\eta_1,\eta_2$ with
$\Phi\eta_0=F$, $\bar\partial\eta_0=d_1\eta_1$ and
$\bar\partial\eta_1=d_2\eta_2$. Since $\eta_2$ is a $(0,2)$-form,
it is closed. We solve $\bar\partial u_2=\eta_2$, put
$\zeta_1=\eta_1-d_2u_2$, and solve $\bar\partial u_1=\zeta_1$.
The second datum is closed since
$\bar\partial\eta_1-d_2\bar\partial u_2=0$. Then
$a=\eta_0-d_1u_1$ has
$\bar\partial a=d_1d_2u_2=0$ and $\Phi a=F$.
The general argument below is this finite descent with estimates
and with domain inclusions at every level. Neither this illustration
nor the general construction assumes injectivity of the final
syzygy map.

\subsection{A cover with buffers for every face}

We choose a maximal $\varepsilon$-separated set of centers
$\{p_i\}$ and fix an ordering.  The balls $B(p_i,\varepsilon)$
cover $\mathbb{P}^N$.  For any fixed $T>0$, there is a bound $M(T)$,
independent of sufficiently small $\varepsilon$, for the number of
centers within $T\varepsilon$ of an arbitrary point.  Indeed the
balls of radius $\varepsilon/2$ about these centers are disjoint
and lie in a ball of radius $(T+1/2)\varepsilon$.
Uniform upper and lower small-ball volume estimates on the compact
manifold give, for example, $M(T)\leq C(1+2T)^{2N}$.

We fix $\kappa>1$.  Smooth bumps equal to one on
$B(p_i,\varepsilon)$ and supported strictly inside
$B(p_i,\kappa\varepsilon)$ can be obtained from a fixed function of
$\dist(p_i,x)^2/\varepsilon^2$ below the injectivity radius.
Their first derivatives are bounded by $C/\varepsilon$.
Writing these bumps as $\beta_i$, the sum $B=\sum_i\beta_i$
is at least one everywhere because the radius-$\varepsilon$ balls
cover $\mathbb{P}^N$. At any point at most $M(\kappa)$ summands or their
first derivatives can contribute, so $|\bar\partial B|\le
C M(\kappa)/\varepsilon$. The functions $\chi_i=\beta_i/B$
are nonnegative and sum to one. Differentiating the quotient gives
$\bar\partial\chi_i=B^{-1}\bar\partial\beta_i-
\beta_iB^{-2}\bar\partial B$, with the asserted bound. Thus we obtain
a partition of unity $\{\chi_i\}$ with
\begin{equation}\label{ana:shrinking-cutoffs}
 \sum_i\chi_i=1,\qquad
 \supp\chi_i\Subset B(p_i,\kappa\varepsilon),\qquad
 |\bar\partial\chi_i|\leq C\varepsilon^{-1}.
\end{equation}
At every point only a uniformly bounded number of terms contributes.
These cutoffs vary with the scale; the atlas partition in
\eqref{ana:patched-metric} was fixed earlier.

For $0\leq q\leq N$, call an increasing tuple
$I=(i_0,\ldots,i_q)$ of degree $q$ \emph{relevant} if the closed
supports of its cutoffs have a common point. The cover may admit
higher-order intersections, but no tuples beyond degree $N$ are
needed or assigned balls or coefficient cochains.  Every nonempty face
of a relevant tuple is relevant.  All its centers lie within
$2\kappa\varepsilon$ of one another.  Anchor $I$ at its smallest
index, and use that same anchor for every permutation of $I$.

We choose fixed output radius factors $s_q^{\mathrm{buf}}$
backward through the presentation. First, we choose
$s_N^{\mathrm{buf}}>2\kappa+2$, and then, for
$q=N,\ldots,1$, choose
\begin{equation}\label{eq:buffer-radii}
 s_{q-1}^{\mathrm{buf}}
 >
 \Lambda(s_q^{\mathrm{buf}})+2\kappa+2.
\end{equation}
For a relevant tuple $I$ of degree $0\le q\le N$,
we define its output and input balls, respectively, by
\[
 V_I=B(p_{i_0},s_q^{\mathrm{buf}}\varepsilon),
 \qquad
 U_I=B(p_{i_0},
       \Lambda(s_q^{\mathrm{buf}})\varepsilon).
\]

\begin{lemma}\label{ana:buffers}
If $J$ is a nonempty proper face of a relevant tuple $I$ of degree
$1\leq q\leq N$, then $U_I\subset V_J$.  A common $\varepsilon_{\rm buf}>0$
makes all these input and output balls admissible in
Proposition~\ref{ana:microscopic-division}, uniformly in the centers.
\end{lemma}

\begin{proof}
We first suppose that $J$ is obtained from $I$ by deleting one index.
The two anchors have distance at most $2\kappa\varepsilon$.
For $x\in U_I$, therefore,
\[
 \operatorname{dist}(x,p_J)
 <
 \bigl(\Lambda(s_q^{\mathrm{buf}})+2\kappa\bigr)
 \varepsilon
 <
 (s_{q-1}^{\mathrm{buf}}-2)\varepsilon.
\]
This includes the face obtained by deleting the smallest index,
whose anchor changes.  For a nonempty proper face of higher codimension, choose a chain of
successive single-index deletions from $I$ to that face. Applying the
codimension-one inclusion at each step and using $V_K\subset U_K$ for
each intermediate tuple $K$ gives the stated inclusion.  There are only finitely many
fixed radii $s_q^{\mathrm{buf}}$ and $\Lambda(s_q^{\mathrm{buf}})$.  Decreasing the admissible
scale places every ball in a suitable buffered chart and makes
all the preceding coordinate, volume, and cutoff estimates hold.
\end{proof}

\subsection{Fixing the scale}

At this point the presentation, metric estimates, local division
radii, and buffers have all been fixed.  In particular the constants
$c_L,C_0,C_1$ of Proposition~\ref{ana:positive-gap} do not depend
on any norm of a division operator.  We choose a fixed $R_*>0$ with
\begin{equation}\label{ana:R-choice}
 C_1/R_*^2<c_L/2.
\end{equation}
We increase an integer $n_0$ so that
\[
 n_0\geq1,\qquad
 R_*/\sqrt{n_0}\leq\varepsilon_{\rm buf},\qquad
 C_0\leq c_Ln_0/4.
\]
We may also require $n_0\geq D_0$.  For $n\geq n_0$ put
\begin{equation}\label{ana:scale}
 \varepsilon=R_*/\sqrt n.
\end{equation}
For both $h=n$ and $h=n+1$,
\begin{equation}\label{ana:simultaneous-scale}
 h\varepsilon^2\leq2R_*^2,
 \qquad
 \Delta(h,\varepsilon)\geq c_Ln/4.
\end{equation}
Local division is consequently uniform at all prescribed radii,
and the minimal global solution satisfies
\begin{equation}\label{ana:scaled-solver}
 \|u\|\leq\frac{2}{\sqrt{c_Ln}}\|\zeta\|
 =b\varepsilon\|\zeta\|,
 \qquad b=\frac{2}{R_*\sqrt{c_L}}.
\end{equation}
We emphasize the order of choices. The uniform regularity degrees
and finite atlas determine the fixed metric and curvature constants.
The compact normalized families determine the local division radii
for each output radius, independently of the line power. We choose
the buffer radii backward from degree $N$ and only then the admissible
scale threshold. We choose $R_*$ using the curvature inequality,
and increase $n_0$ to fit this threshold. Thus no local division
constant can affect the curvature inequality. The constants for
local division may now contain factors such as
$\exp(2C\Lambda(s_q^{\mathrm{buf}})^2R_*^2)$.  These are fixed constants.
They do not enter the earlier curvature constant $C_1$, and the
finite correction below imposes no smallness condition on them.

\subsection{Holomorphic cochains with uniform total norm}

We fix $n\geq n_0$, one of the two powers $h\in\{n,n+1\}$, and
$F\in H^0(\mathbb{P}^N,\cI\otimes L^h)$.  We write $\|F\|$ for its original
scalar $L^2$ norm.  For a cochain of degree $q-1$, we define its
alternating difference on an increasing relevant tuple $I$ of degree
$1\leq q\leq N$ by
\[
 (\delta c)_I
 =\sum_{\ell=0}^q(-1)^\ell
 c_{I\setminus i_\ell}|_{U_I}.
\]
Lemma~\ref{ana:buffers} makes the restrictions meaningful.
These are differences of sections of the same actual bundle;
columns in two chart frames must first be identified as such sections.

\begin{proposition}\label{ana:holomorphic-cochains}
For every increasing relevant tuple $I$ of degree $0\leq q\leq N$, there is a
holomorphic section $c_{q,I}$ of $E_q\otimes L^h$ on $V_I$, linear
in $F$, with
\begin{equation}\label{ana:cochain-identities}
 \Phi c_{0,i}=F|_{V_i},
 \qquad
 d_qc_{q,I}=(\delta c_{q-1})_I|_{V_I}\quad(1\leq q\leq N).
\end{equation}
Moreover
\begin{equation}\label{ana:cochain-mass}
 S_q:=\sum_{I\ \mathrm{increasing,\ relevant}}
       \|c_{q,I}\|_{q,h;V_I}^2
 \leq C_q\|F\|^2,
\end{equation}
with constants independent of $n$, $h$, and $F$.
\end{proposition}

\begin{proof}
At degree zero we divide $F$ on $U_i$ using
Proposition~\ref{ana:microscopic-division}.  At degree one,
$c_{0,j}-c_{0,i}$ is a holomorphic section on $U_{ij}$ whose image
under $\Phi$ is $F-F=0$ on that entire ball.  It belongs to the image sheaf of $d_1$, so by division we obtain $c_{1,ij}$ on $V_{ij}$.

Suppose the construction has reached degree $q-1$, with $2\leq q\leq N$.
All face identities hold on the whole new input $U_I$, by the
buffer inclusions, including those for faces of faces.  Hence
\[
 d_{q-1}(\delta c_{q-1})_I
 =(\delta^2c_{q-2})_I=0\quad\hbox{on }U_I.
\]
Indeed the two orders of deleting positions $r<s$ have signs
$(-1)^{r+s-1}$ and $(-1)^{r+s}$ and restrict the same section to
the same domain. Exactness at $E_{q-1}$ supplies analytic image-sheaf
membership of this holomorphic datum throughout $U_I$. This is within
the exactness range because $q\leq N$.  By local division we obtain $c_{q,I}$ on $V_I$.  We fix all charts and the linear local division
operators before varying $F$; induction proves linear dependence.
For permutations of an increasing tuple we use the corresponding sign
and the original anchored ball, and put the cochain equal to zero
on tuples with repeated indices.

Let $L_q$ be a common local-division norm at output factor $s_q^{\mathrm{buf}}$
for the fixed bound $K=2R_*^2$ on $h\varepsilon^2$.  Every point belongs to at most
$M(\Lambda(s_q^{\mathrm{buf}}))$ initial input balls.  Thus
\[
 S_0\leq L_0^2M(\Lambda(s_q^{\mathrm{buf}}))\|F\|^2.
\]
At degree $1\leq q\leq N$, the squared norm of the alternating sum is
bounded by $q+1$ times the sum of the squared face norms on $U_I$.
Each is at most the full norm on $V_J$.  A fixed nonempty face $J$
has at most $M(2\kappa)$ relevant increasing cofaces: the new
index must have its center within $2\kappa\varepsilon$ of the
anchor of $J$.  Consequently
\[
 S_q\leq L_q^2(q+1)M(2\kappa)S_{q-1}.
\]
Iterating at most $N$ times, we obtain \eqref{ana:cochain-mass}.
The number of balls, which grows as the scale decreases, never
enters this bound.
\end{proof}

\subsection{The cochain-to-form identities}

We set $\alpha_i=\bar\partial\chi_i$ and define, with ordered tuples,
\begin{equation}\label{ana:eta-definition}
 \eta_q=
 \sum_{i_0,\ldots,i_q}
 \chi_{i_0}\alpha_{i_1}\wedge\cdots\wedge\alpha_{i_q}
 c_{q,i_0\ldots i_q},\qquad 0\leq q\leq N.
\end{equation}
For $q=0$ this means $\eta_0=\sum_i\chi_ic_{0,i}$.
An irrelevant tuple contributes zero.  The common cutoff support
is compactly contained in $V_I$, since $s_q^{\mathrm{buf}}>2\kappa+2$.
We can therefore extend every summand by zero to obtain a smooth global form.

\begin{proposition}\label{ana:eta-identities}
These forms have values in $E_q\otimes L^h$ and satisfy
\begin{equation}\label{ana:eta-complex}
 \Phi\eta_0=F,
 \qquad
 d_q\eta_q=\bar\partial\eta_{q-1}\quad(1\leq q\leq N),
 \qquad
 \bar\partial\eta_N=0.
\end{equation}
There are fixed constants $P_q$ such 
There are fixed constants $C_q^{\mathrm{form}}$,
independent of $n$, $h$, and $F$, such that
\begin{equation}\label{ana:eta-estimate}
 \|\eta_q\|_{H_{q,\varepsilon}\otimes L^h}
 \le C_q^{\mathrm{form}}\varepsilon^{-q}\|F\|,
 \qquad 0\le q\le N.
\end{equation}
All forms depend linearly on $F$.
\end{proposition}

\begin{proof}
We obtain the first identity from $\sum_i\chi_i=1$ and
\eqref{ana:cochain-identities}.  For the second, apply $d_q$ to
\eqref{ana:eta-definition} and expand the alternating difference.
The face deleting $i_0$ gives
\[
 \sum_{i_1,\ldots,i_q}
 \alpha_{i_1}\wedge\cdots\wedge\alpha_{i_q}
 c_{q-1,i_1\ldots i_q},
\]
because the sum over $i_0$ of $\chi_{i_0}$ is one.  This is
exactly $\bar\partial\eta_{q-1}$: the local coefficient sections
are holomorphic and $\bar\partial\alpha_i=0$.
For a face deleting $i_\ell$ with $\ell\geq1$, its coefficient is
independent of that index.  Summing the factor $\alpha_{i_\ell}$
in its unchanged wedge position gives zero, since
$\sum_i\alpha_i=0$.  These calculations are pointwise legitimate:
every tuple whose cutoff factors contribute at a point is relevant,
and every coefficient used there is defined on a neighborhood of
that point.  Repeated-index terms obey the alternating convention.
There is no factorial prefactor in \eqref{ana:eta-definition} or
in these identities.  Finally a $(0,N+1)$-form on $\mathbb{P}^N$ vanishes,
which proves the last identity in \eqref{ana:eta-complex} without any
injectivity assumption on $d_N$.

To estimate the norm, we note that at any point at most $M(\kappa)^{q+1}$
ordered tuples can contribute.  The norm of an individual summand
is at most $C_q\varepsilon^{-q}|c_{q,I}|$ on its support, by
\eqref{ana:shrinking-cutoffs} and the exterior-product inequality.
We apply the squared Cauchy--Schwarz inequality to the finite sum and
integrate.  For each increasing relevant tuple there are at most
$(q+1)!$ ordered versions, all using the same anchored section and
ball.  Thus
\[
 \|\eta_q\|^2
 \leq C_q M(\kappa)^{q+1}(q+1)!
        \varepsilon^{-2q} S_q.
\]
Proposition~\ref{ana:holomorphic-cochains} proves
\eqref{ana:eta-estimate}, also at $q=0$, where no differentiated
cutoff appears.  Linearity follows from the construction with
fixed cutoffs.
\end{proof}

\subsection{Finite descending correction}

At a fixed $(n,h)$ all the following spaces carry the same coefficient
metrics and the same line twist.  We begin with $\zeta_N=\eta_N$.
For $q=N,\ldots,1$, we solve in the coefficient bundle $E_q\otimes L^h$
\begin{equation}\label{ana:descent}
 \bar\partial u_q=\zeta_q,
 \qquad
 \zeta_{q-1}=\eta_{q-1}-d_qu_q.
\end{equation}
The form $u_q$ has type $(0,q-1)$, so the two terms defining
$\zeta_{q-1}$ have the same bundle and form degree.

\begin{lemma}\label{ana:descent-closed}
Every input $\zeta_q$ in \eqref{ana:descent} is closed.
The procedure gives a holomorphic section $a=\zeta_0$ with
$\Phi a=F$.
\end{lemma}

\begin{proof}
The top form is closed by \eqref{ana:eta-complex}.  If $\zeta_q$
is closed, Proposition~\ref{ana:positive-gap} and
\eqref{ana:simultaneous-scale} give its smooth minimal solution.
The presentation maps are holomorphic and commute with
$\bar\partial$, so
\[
 \bar\partial\zeta_{q-1}
 =d_q\eta_q-d_q\bar\partial u_q
 =d_q(\eta_q-\zeta_q).
\]
For $q=N$ the last difference is zero.  For $q<N$, the preceding
step says $\eta_q-\zeta_q=d_{q+1}u_{q+1}$, and the expression
vanishes because $d_qd_{q+1}=0$.  This proves closedness at every
step, and hence justifies each subsequent application of the solver.
At degree zero, we conclude from closedness that $a$ is holomorphic.  Finally
\[
 \Phi a=\Phi\eta_0-\Phi d_1u_1=F.
\]
If a coefficient stage is zero, the corresponding data and solutions
are zero; the same induction applies.
\end{proof}

\begin{lemma}\label{ana:descent-norm}
We define constants $B_q$ recursively by
\[
 B_N=C_N^{\mathrm{form}},
 \qquad
 B_{q-1}=C_{q-1}^{\mathrm{form}}+bB_q,
 \qquad q=N,\ldots,1,
\]
where $b$ is the constant in \eqref{ana:scaled-solver}. Then
\[
 \|\zeta_q\|
 \le B_q\varepsilon^{-q}\|F\|,
 \qquad 0\le q\le N,
\]
and
\[
 \|u_q\|
 \le bB_q\varepsilon^{-(q-1)}\|F\|,
 \qquad 1\le q\le N.
\]
In particular, the holomorphic lift $a=\zeta_0$ satisfies
\begin{equation}\label{ana:final-lift-bound}
 \|a\|_{H_{0,\varepsilon}\otimes L^h}
 \le
 \left(\sum_{q=0}^N b^qC_q^{\mathrm{form}}\right)
 \|F\|_{L^h}.
\end{equation}
\end{lemma}

\begin{proof}
Since $\zeta_N=\eta_N$, the initial estimate follows
from \eqref{ana:eta-estimate}.
Suppose that
\[
 \|\zeta_q\|\le B_q\varepsilon^{-q}\|F\|.
\]
By Lemma \ref{ana:descent-closed}, $\zeta_q$ is $\bar\partial$-closed.
Applying the solution estimate \eqref{ana:scaled-solver}, we obtain
\[
 \|u_q\|
 \le b\varepsilon\|\zeta_q\|
 \le bB_q\varepsilon^{-(q-1)}\|F\|.
\]
Since $d_q$ is contractive in the coefficient metrics,
and this contractivity is preserved after tensoring
with $L^h$ and the bundle of differential forms,
we have
\[
 \begin{aligned}
 \|\zeta_{q-1}\|
 &\le \|\eta_{q-1}\|+\|d_qu_q\|\\
 &\le \|\eta_{q-1}\|+\|u_q\|\\
 &\le
 \left(C_{q-1}^{\mathrm{form}}+bB_q\right)
 \varepsilon^{-(q-1)}\|F\|\\
 &=B_{q-1}\varepsilon^{-(q-1)}\|F\|.
 \end{aligned}
\]
This completes the descending induction.

At degree zero, we obtain $\|a\|\le B_0\|F\|$.
Expanding the finite recursion gives
\[
 B_0=\sum_{q=0}^N b^qC_q^{\mathrm{form}}.
\]
All these constants are independent of $n$ and $h$.
In particular, the argument requires no smallness
condition on $b$ or on the local division constants.
The closedness of each correction datum follows from
the complex identities, rather than from estimates
on derivatives of the corrected forms.
\end{proof}

\begin{proposition}\label{ana:global-lift}
For all $n\geq n_0$ and $h\in\{n,n+1\}$ there is a linear map
\[
 T_{n,h}:H^0(\mathbb{P}^N,\cI\otimes L^h)
 \longrightarrow H^0(\mathbb{P}^N,E_0\otimes L^h)
\]
such that $\Phi T_{n,h}=\id$ and
\[
 \|T_{n,h}F\|_{H_{0,R_*/\sqrt n}\otimes L^h}
 \leq C\|F\|_{L^h}.
\]
The datum on the right has its original scalar Fubini--Study norm.
We use the same coefficient metric at the two adjacent line powers.
\end{proposition}

\begin{proof}
For fixed $(n,h)$ we fix the cover, partition, frames, and local
minimum-norm division operators.  Proposition~\ref{ana:holomorphic-cochains}
and formula \eqref{ana:eta-definition} depend linearly on $F$.
Each minimal Dolbeault solution is linear on the closed subspace of
admissible inputs by Proposition~\ref{ana:positive-gap}.
Lemma~\ref{ana:descent-closed} ensures that every corrected input
belongs to that subspace.  The finite recursion consequently gives
a linear $a(F)$.  We set $T_{n,h}F=a(F)$.  Its right-inverse identity
and holomorphy follow from Lemma~\ref{ana:descent-closed}, and its
norm bound is \eqref{ana:final-lift-bound}.  All constants in that
bound were fixed before $n$ varies, and
\eqref{ana:simultaneous-scale} holds at both powers.
We may choose the two linear operators separately; no identification
of their source spaces is asserted.
\end{proof}

\subsection{Completion of the lifting theorem}

\begin{proof}[Proof of Theorem~\ref{geo:lifting}]
For a proper nonzero ideal sheaf in positive projective dimension,
the complete presentation of Proposition~\ref{alg:complete} supplies
the fixed row $\Phi$ and degree $D=D_0$.  The graph metrics
\eqref{ana:graph-definition} are smooth and positive at every
admissible radius, and make $\Phi$ contractive into the original
scalar norm.  Proposition~\ref{ana:homogeneous-derivative} proves
the homogeneous derivative estimate with target $L^{-1}$ in that
same metric.  Proposition~\ref{ana:graph-curvature} proves its
$O(\varepsilon^{-2})$ curvature bound, measured in its own fiber
norm.  Finally Proposition~\ref{ana:global-lift} proves the bounded
holomorphic lifting assertion at $\varepsilon=R_*/\sqrt n$, for
both $h=n$ and $h=n+1$, with the original norm on the datum.
Taking the maximum of the finitely many constants, we obtain the single
constant permitted in the statement.

For the unit ideal we take $D=0$, $V=\C$, $\Phi=1$, and the constant
coefficient metric.  The identity lift has norm one and the row
derivative and curvature vanish.  For completeness, on $\Pj^0$ every nonzero ideal
sheaf is the unit ideal, so its lift is also the identity.
The zero ideal requires no lifting: its discarded coordinate
components vanish.  This proves every asserted clause.
\end{proof}

\begin{remark}\label{ana:constant-order}
We make the choices in a definite order. The complete presentation fixes
the polynomial degrees.  Those degrees and the fixed atlas determine
the Gram estimates and the curvature constants.  Compact normalized
division determines an input radius for each output radius.  The
finite buffer recursion fixes the shrinking-cover geometry.
Only then do we choose $R_*$ to dominate the curvature loss.
The line-weight comparison changes the local-division constants
after this choice, but it does not change the curvature constants.
The finite descent accommodates those constants without any
smallness assumption.  This order accounts for all uniformities
in the lifting theorem.
\end{remark}


\section{Step 3: Return to Drury--Arveson space}
\label{sec:completion}

We have proved Theorem~\ref{geo:lifting}. We now
return to the original scalar norm, obtain the operator estimate, and
transfer it to a vector module before the Hilbert function is used.

\subsection{The scalar estimate in the original Hilbert norm}
\begin{proposition}\label{geo:transfer}
Theorem~\ref{geo:lifting} implies the scalar degree estimate
\eqref{eq:main-degree} for every homogeneous ideal in every finite
number of variables.
\end{proposition}
\begin{proof}
We have already treated the exceptional cases and the finite change under saturation. We work in a degree where
$I_n=I_n^{\mathrm{sat}}$ and $I_{n+1}=I_{n+1}^{\mathrm{sat}}$.
For $F\in I_{n+1}$, we use Theorem~\ref{geo:lifting} with $h=n+1$.
The entries of its global holomorphic lift $a$ are homogeneous
polynomials of degree $n+1-D$.  Differentiating the equality $F=\Phi a$, we obtain
\[
 \partial_{z_j}F=(\partial_{z_j}\Phi)a+\Phi\partial_{z_j}a.
\]
Each entry of the fixed row $\Phi$ belongs to the saturated
ideal in degree $D$. Differentiating any entry of the coefficient
vector produces a polynomial of degree $n-D$, if this is nonnegative.
The product with the corresponding entry of $\Phi$ belongs to
$I_n^{\mathrm{sat}}$ by ideal multiplication, independently of the
size of that derivative. Consequently the sum of these terms
belongs to $I_n^{\mathrm{sat}}=I_n$; if the coefficient degree is
zero it is simply zero.  Using the common orthogonal projection in Lemma~\ref{geo:norm}, we therefore obtain
\[
 P_{Q_n}\partial_{z_j}F
   =P_{Q_n}\bigl((\partial_{z_j}\Phi)a\bigr).
\]
We apply \eqref{eq:geo-differentiate} pointwise and
\eqref{eq:geo-lift-bound} after integration.  Since the projection is
contractive in the scalar norm,
\begin{equation}\label{eq:geo-projected-derivative}
 \|P_{Q_n}\partial_{z_j}F\|_{\FS,n}
       \le C''\varepsilon^{-1}\|F\|_{\FS,n+1}.
\end{equation}
On degree $n+1$, the Drury--Arveson adjoint is
$Z_j^*=(n+1)^{-1}\partial_{z_j}$, and
$B_{j,n}^*F=P_{Q_n}Z_j^*F$.  Consequently
\begin{align}
 \|B_{j,n}^*F\|_{\DA,n}
 &\le \frac{C''}{\varepsilon(n+1)}
       \sqrt{\frac{c_n}{c_{n+1}}}\|F\|_{\DA,n+1}\notag\\
 &=\frac{C''}{\varepsilon\sqrt{(n+1)(n+m)}}
       \|F\|_{\DA,n+1}.\label{eq:geo-exact-transfer}
\end{align}
At $\varepsilon=R_*/\sqrt n$, summing the squared operator bounds gives
\[
 A_n\le
 \frac{m(C'')^2n}{R_*^2(n+1)(n+m)}I_{Q_n}
 \le\frac{C_I}{n+1}I_{Q_n}.
\]
We absorb the finitely many omitted degrees into $C_I$.
\end{proof}

We use no derivative of the metric or of the chosen lift in the last
calculation.  The derivative of the lift has been removed by ideal
membership.  Curvature and derivatives of the metric are needed earlier,
to construct a holomorphic lift with a compatible norm bound.  This
distinction separates the analytic construction from its operator
consequence.

In particular, Proposition~\ref{geo:transfer}, with the finite-degree
adjustment of Lemma~\ref{geo:saturation}, proves for every homogeneous
scalar ideal, in any finite number of variables, the estimate
\begin{equation}\label{eq:complete-scalar}
 A_n^I\le\frac{C_I}{n+1}I_{Q_n}.
\end{equation}

\subsection{An auxiliary-variable realization}
Let $M\subset R\otimes\C^r$, where $R=\C[z_1,\ldots,z_m]$, and
let $T=(R\otimes\C^r)/M$.  Adjoin $r$ variables and define
\[
 R'=R[y_1,\ldots,y_r],\qquad
 \lambda(f)=\sum_{b=1}^r f_b(z)y_b\quad(f\in R\otimes\C^r).
\]
We choose homogeneous generators $f^{(1)},\ldots,f^{(s)}$ of $M$ and set
\begin{equation}\label{eq:complete-encoding}
 J_M=(\lambda(f^{(1)}),\ldots,\lambda(f^{(s)}))+(y_1,\ldots,y_r)^2.
\end{equation}
The auxiliary-variable realization and its degreewise normalization
occur in Shalit's stable-division argument
\cite[Section~5, equation~(5.1)]{ShalitDivision}.  The square-zero
relations in \eqref{eq:complete-encoding} let us record the complete
positive defect exactly.

For an $R$-module $T$, the square-zero extension $R\ltimes T(-1)$
is the algebra on the vector space $R\oplus T(-1)$ with product
$(a,t)(a',t')=(aa',at'+a't)$. Its second summand is an ideal whose
product with itself is zero, but multiplication by $R$ retains the
original module action. This construction is the algebraic reason
for imposing every product $y_by_c$ as a relation. The degree
shifts matter: the class of $f\in T_n$ becomes degree $n+1$ when
represented by $\lambda(f)$. We will work in that specific
auxiliary-degree sector and not estimate the dimension of the
original module by the dimension of this larger algebra.

\begin{lemma}\label{complete:ring}
The ideal $J_M$ is homogeneous separately in the $z$ and $y$ variables.
Its component of $y$-degree one is exactly $\lambda(M)$.  There is an
isomorphism of graded algebras
\[
 R'/J_M\simeq R\ltimes T(-1),
 \qquad (a,t)(a',t')=(aa',at'+a't).
\]
In particular $\sqrt{J_M}=(y_1,\ldots,y_r)$ and
$\dim R'/J_M=m$.
\end{lemma}
\begin{proof}
Each generator of \eqref{eq:complete-encoding} has definite degree in
both groups of variables.  In $y$-degree zero no relation occurs; in
$y$-degree one only coefficients in $R$ can multiply the first set of
generators.  Their image is exactly $\lambda(M)$ by $R$-linearity.
All higher $y$-degrees vanish.  Thus a class has the displayed pair
description, and multiplication has the claimed form.  The shift is
$T(-1)_\ell=T_{\ell-1}$, because each auxiliary variable adds one to
total degree.  Finally
$(y)^2\subset J_M\subset(y)$ and $(y)$ is prime, with quotient $R$.
\end{proof}

Let $Q^{J_M}\subset H_{m+r}^2$ be the scalar quotient
associated with the auxiliary ideal $J_M$.
Since $J_M$ is homogeneous separately in the $z$ and
$y$ variables, $Q^{J_M}$ inherits an orthogonal
decomposition by the two degrees.
We denote by $Q_n^{\mathrm{aux}}$ its sector of
$z$-degree $n$ and $y$-degree one.
In particular, $Q_n^{\mathrm{aux}}$ has total degree
$n+1$.
We write $Q_n^M$ for the degree-$n$ component of the
original vector-valued quotient.
\begin{lemma}\label{complete:unitary}
The map
\begin{equation}\label{eq:complete-unitary}
 U_n:Q_n^M\longrightarrow Q_n^{\mathrm{aux}},
 \qquad
 U_nq=\sqrt{n+1}\sum_{b=1}^r q_b(z)y_b
\end{equation}
is unitary. In total degree $\ell\ge1$,
\[
 Q_\ell^{J_M}
 =
 H_{m,\ell}^2\oplus Q_{\ell-1}^{\mathrm{aux}}
\]
orthogonally; in degree zero, $Q_0^{J_M}=\mathbb C$.
\end{lemma}
\begin{proof}
For $|\alpha|=n$,
\[
 \|z^\alpha e_b\|_{H_m^2\otimes\C^r}^2=\frac{\alpha!}{n!},
 \qquad
 \|z^\alpha y_b\|_{H_{m+r}^2}^2=\frac{\alpha!}{(n+1)!}.
\]
Distinct displayed monomials are orthogonal.  Therefore the same
formula defines a unitary $\widetilde U_n$ between the full vector
degree space and the full scalar sector linear in $y$.  It maps $M_n$
onto the corresponding sector of $J_M$ by Lemma~\ref{complete:ring},
and hence maps their orthogonal complements onto one another.
The other surviving sector is the pure $z$ sector, whose monomial norms
are unchanged.  All higher auxiliary degrees lie in the ideal.
\end{proof}

Let $A^{J_M}_\ell$ be the scalar defect square, summed over all
$m+r$ coordinates.  The following equality supplies the quantitative
transfer.

\begin{proposition}\label{complete:defect}
The operator $A^{J_M}_\ell$ preserves the two surviving auxiliary-degree
sectors.  On the linear sector,
\begin{equation}\label{eq:complete-defect}
 U_n^*A^{J_M}_{n+1}U_n
 =\frac{n+1}{n+2}A_n^M+\frac{r+1}{n+2}I_{Q_n^M}.
\end{equation}
On the pure $z$ sector in total degree $\ell$,
\[
 0\le A^{J_M}_\ell|_{H_{m,\ell}^2}
       \le\frac r{\ell+1}I.
\]
\end{proposition}
\begin{proof}
Let $q\in Q_n^M$.  For each original coordinate,
\[
 z_jU_nq
 =\sqrt{\frac{n+1}{n+2}}\,\widetilde U_{n+1}(z_jq).
\]
The scalar projection to the ideal preserves auxiliary degree, and
$\widetilde U_{n+1}$ intertwines its restriction in auxiliary degree
one with projection onto $M_{n+1}$.  The sum of the discarded squared
norms over the $z$ coordinates is therefore
\[
 \frac{n+1}{n+2}\langle A_n^Mq,q\rangle.
\]
Every product $y_bU_nq$ has auxiliary degree two and belongs to $J_M$,
so it is discarded in full.  On a normalized monomial of total degree
$\ell$ and auxiliary multi-degree $\beta$, the squared forward weight
in $y_b$ is $(\beta_b+1)/(\ell+1)$.  Summing over $b$ gives
$(|\beta|+r)/(\ell+1)$.  For each fixed coordinate, distinct monomials
have orthogonal images.  Hence the same identity holds on any vector
in a sector of fixed auxiliary degree. On $Q_n^{\mathrm{aux}}$, this sum is
$(r+1)/(n+2)$ times the squared norm of the input vector. Adding and polarizing proves
\eqref{eq:complete-defect}.

For a pure $z$ input, multiplication by an original coordinate remains
in the pure sector and has zero defect.  The sum of squared unprojected
auxiliary products is $r/(\ell+1)$ times the input norm; projection can
only decrease it.  For each coordinate, the discarded images of the
two input sectors have distinct auxiliary degrees, so their cross
inner products vanish.  This proves the reducing statement as well.
\end{proof}

\begin{proof}[Proof of Theorem~\ref{thm:main}]
We apply \eqref{eq:complete-scalar} to $J_M$ in $m+r$ variables
and restrict the resulting estimate to
$Q_n^{\mathrm{aux}}$ in total degree $n+1$. Dropping the positive
scalar summand in \eqref{eq:complete-defect} gives
\[
 \frac{n+1}{n+2}A_n^M\le\frac{C_{J_M}}{n+2}I,
 \qquad
 A_n^M\le\frac{C_{J_M}}{n+1}I.
\]
We have thus proved \eqref{eq:main-degree} in the original vector norm.
We deduce \eqref{eq:main-commutator} and the restriction block bound from Lemma~\ref{op:commutators}. Using Lemma~\ref{op:hilbert} and the singular-value calculation in Lemma~\ref{op:counting}, we obtain all the claimed weak and strong Schatten classes, including the defect column and the projection commutators.
We also accounted for the finite-dimensional quotient case there.
Lemma~\ref{op:graded} identifies the closed graded formulation with
the polynomial formulation.
\end{proof}

The order of these steps is essential for the dimension-sensitive
conclusion.  The scalar algebra in Lemma~\ref{complete:ring} has
dimension $m$ even when $T$ has smaller support.  Counting in that
scalar algebra would discard the sharper exponent.  Equation
\eqref{eq:complete-defect} first returns the degree estimate to
$Q_n^M$, whose dimension is the Hilbert function of $T$.

The square-zero part also explains the relevance of nonreduced ideals
in the scalar theorem.  Replacing $J_M$ by its radical removes the
entire sector that encodes the vector quotient.  Algebraic support
alone does not retain the information needed for this transfer.


\section{Other ball norms, analytic ideals, and the endpoint}
\label{sec:extensions}

\subsection{The standard kernel family}
For $\lambda>0$, let $H^{(\lambda)}$ be the Hilbert function space on
$\mathbb B_m$ with reproducing kernel
\[
 K_\lambda(z,w)=(1-\langle z,w\rangle)^{-\lambda}.
\]
Its orthogonal monomial norms are
\begin{equation}\label{eq:ext-norm}
 \|z^\alpha\|_\lambda^2
 =\frac{\alpha!\,\Gamma(\lambda)}{\Gamma(|\alpha|+\lambda)}.
\end{equation}
These formulas follow by the binomial expansion of the kernel.  The
cases $\lambda=1$, $m$, and $m+1$ are respectively the
Drury--Arveson, sphere Hardy, and normalized ball Bergman spaces.

\begin{corollary}\label{ext:kernel-family}
The quotient and restriction conclusions of Theorem~\ref{thm:main}
hold for the closure of the same homogeneous submodule in
$H^{(\lambda)}\otimes\C^r$, for every $\lambda>0$, with the same
dimension exponents.  The constants may depend on $\lambda$.
\end{corollary}
\begin{proof}
The inner products on each homogeneous degree differ by a scalar,
so the degreewise unitary between normalized monomial bases maps the
ideal or submodule spaces, and their orthogonal complements, onto
one another.  Under this identification the quotient coordinate in
degree $n\to n+1$ is
\[
 S_{j,n}^{(\lambda)}=r_nS_{j,n},\qquad
 r_n=\sqrt{\frac{n+1}{n+\lambda}}.
\]
For $n\ge1$ this gives
\begin{align*}
 [(S_i^{(\lambda)})^*,S_j^{(\lambda)}]|_{Q_n}
 &=r_n^2[S_i^*,S_j]|_{Q_n}
   +(r_n^2-r_{n-1}^2)S_jS_i^*|_{Q_n},\\
 r_n^2-r_{n-1}^2
 &=\frac{\lambda-1}{(n+\lambda)(n+\lambda-1)}=O(n^{-2}).
\end{align*}
The sequence $r_n$ is bounded and the quotient coordinates in the
Drury--Arveson norm are contractions.  Hence the new commutator blocks
have norm $O((n+1)^{-1})$.  The same calculation applies to the
restriction coordinates.  The defect blocks are multiplied by $r_n$,
so their norms remain $O((n+1)^{-1/2})$.  The finitely many initial
terms are harmless, and the original Hilbert function supplies the
same singular-value counts.  Finite quotients remain finite.
\end{proof}

In this transfer we use actual scalar weights between adjacent degrees. We do not appeal to invariance of essential normality under an
arbitrary equivalent renorming.

\subsection{Ideals with holomorphic coefficients}
Ordinary homogeneity is particularly restrictive for a holomorphic
function at the origin.

\begin{proposition}\label{ext:holomorphic}
If a nonzero holomorphic function near $0\in\C^m$ satisfies
$f(tz)=t^df(z)$ for positive $t$ near $1$, then $d$ is a nonnegative
integer and $f$ is a homogeneous polynomial of degree $d$.
If $p_1,\ldots,p_s$ are homogeneous polynomials,
$J=(p_1,\ldots,p_s)\subset R$, and
$\mathfrak a=\sum_jp_j\cO(\mathbb B_m)$, then
\[
 \overline{\mathfrak a\cap H^{(\lambda)}}^{\,H^{(\lambda)}}
   =\overline J^{\,H^{(\lambda)}}.
\]
\end{proposition}
\begin{proof}
Differentiating in $t$ at $1$ gives
$\sum_jz_j\partial_jf=df$.  For the convergent Taylor expansion
$f=\sum_\alpha c_\alpha z^\alpha$, uniqueness of coefficients gives
$(|\alpha|-d)c_\alpha=0$.  Some coefficient is nonzero, so
$d\in\N$, and only the finitely many monomials of that degree occur.
Equality extends throughout a connected domain of holomorphy by
the identity theorem.

For the second assertion, we write $F=\sum_jp_jg_j\in H^{(\lambda)}$
with the $g_j$ holomorphic.  If $\deg p_j=d_j$, the homogeneous
component is $F_n=\sum_jp_j(g_j)_{n-d_j}\in J_n$.  Its finite
homogeneous sums belong to $J$ and converge in the Hilbert norm,
because the homogeneous spaces are orthogonal.  Thus
$\mathfrak a\cap H^{(\lambda)}\subset\overline J$.
The reverse inclusion after closure follows from
$J\subset\mathfrak a\cap H^{(\lambda)}$.
\end{proof}

The coefficients $g_j$ in this proposition need not be bounded
multipliers.  If one starts with infinitely many homogeneous
holomorphic generators, the polynomial ring is Noetherian, so a
finite subset generates the same polynomial ideal and its analytic
extension.  The same component argument applies to homogeneous
vector-valued generators and to global sections of the generated
analytic ideal sheaf, since a representation near the origin already
determines every Taylor component.

\subsection{Why the endpoint is stated weakly}
The weak exponent cannot in general be replaced by its strong
counterpart, even for a simple homogeneous ideal.

\begin{example}\label{ext:sharp}
Let $I=(x^2)\subset\C[x,y]$.  Its quotient has support dimension one.
In degree $n\ge1$ an orthonormal basis is
\[
 e_n=y^n,\qquad f_n=\sqrt n\,xy^{n-1},
\]
and degree zero is spanned by $e_0=1$.  The compressed $x$ coordinate
satisfies
\[
 S_xe_n=\frac1{\sqrt{n+1}}f_{n+1},\qquad S_xf_n=0,
\]
and therefore
\[
 [S_x^*,S_x]e_n=\frac1{n+1}e_n,
 \qquad [S_x^*,S_x]f_n=-\frac1n f_n\quad(n\ge1).
\]
The eigenvalue on $e_0$ is $1$.  Its singular values thus belong to
$\ell^{1,\infty}$ but not to $\ell^1$.  In particular the quotient
commutator is in $\cS_{1,\infty}$ and in every $\cS_p$, $p>1$,
but is not trace class.  This example retains a transverse first
order nilpotent direction; replacing the ideal by its radical would
remove that direction.
\end{example}

\subsection{Exact defect blocks for homogeneous ideals}

We test the degree estimate and its normalization by the following calculations
without using the projective lifting theorem. They also exhibit the
necessity of the order $(n+1)^{-1}$ in the general defect-square bound.

\begin{example}[An ideal with an embedded component]\label{ext:embedded-defect}
Let $I=(x^2,xy)\subset\C[x,y,z]$. Its radical is $(x)$, while the class
of $x$ in $R/I$ is nonzero and has annihilator $(x,y)$. In particular
$(x,y)$ is an embedded associated prime. For $n\geq1$, an orthogonal
monomial basis of $Q_n$ consists of
\[
 z^n,\quad y^a z^{n-a}\ (1\leq a\leq n),\quad xz^{n-1}.
\]
We normalize these monomials in the Drury--Arveson norm. In any coordinate
$z_j$, the squared forward weight of a normalized degree-$n$ monomial
$z^\alpha$ is $(\alpha_j+1)/(n+1)$. Distinct monomials have orthogonal
images under each fixed coordinate, and projection onto the ideal
keeps exactly the monomials in the ideal. Thus each $B_{j,n}^*B_{j,n}$
is diagonal in this basis. No product of $z^n$ by a coordinate lies in
$I$. For $y^a z^{n-a}$, $a\geq1$, only the product by $x$ lies in $I$;
its discarded squared weight is $1/(n+1)$. For $xz^{n-1}$, the products
by $x$ and $y$ lie in $I$, with discarded squared weights $2/(n+1)$
and $1/(n+1)$, respectively. Consequently the eigenvalues of $A_n$,
with multiplicities indicated, are
\begin{equation}\label{eq:embedded-spectrum}
 \begin{array}{c|ccc}
 \text{eigenvalue}&0&1/(n+1)&3/(n+1)\\
 \text{multiplicity}&1&n&1.
 \end{array}
\end{equation}
In particular $\|A_n\|=3/(n+1)$ and $\Tr A_n=(n+3)/(n+1)$.
The trace agrees with \eqref{eq:op-trace}, since $h_n=n+2$ for
$n\geq1$. At degree zero $A_0=0$, as all coordinates survive.
This example has $\delta=2$ and shows that one cannot replace
$O((n+1)^{-1})$ by $o((n+1)^{-1})$ for every homogeneous ideal.
\end{example}

\begin{example}[A nonmonomial quadric]\label{ext:quadric-defect}
Let $m\geq3$, $q(z)=\sum_{j=1}^m z_j^2$, and $I=(q)$. On the
polynomial ring equip monomials with the Fischer inner product
$\langle z^\alpha,z^\beta\rangle_{\rm F}=\alpha!\,\delta_{\alpha\beta}$.
Within degree $n$ it is $n!$ times the Drury--Arveson inner product;
therefore the orthogonal projections coincide. The Fischer adjoint of multiplication by $z_j$ is
$\partial_{z_j}$: on monomials with $\beta=\alpha+e_j$, both pairings
are $\beta!$, and otherwise both are zero. Taking the adjoint of
each squared multiplication and summing, we see that multiplication
by $q$ has adjoint $\Delta=\sum_j\partial_{z_j}^2$. Here harmonic
means annihilated by this complex polynomial differential operator,
not by the real Laplace operator on $\C^m$. Since the orthogonal
complement of $qR_{n-2}$ is the kernel of its adjoint, we obtain
$Q_n=\ker(\Delta|_{R_n})$. We take $p\in Q_n$, $n\geq1$. The projection
of $z_jp$ onto $I_{n+1}=qR_{n-1}$ is
\begin{equation}\label{eq:quadric-projection}
 B_{j,n}p=\frac{q\,\partial_{z_j}p}{2n+m-2}.
\end{equation}
Indeed, $\partial_{z_j}p$ is harmonic and homogeneous of degree $n-1$.
The product rule and Euler identity give
\[
 \Delta(z_jp)=2\partial_{z_j}p,\qquad
 \Delta(q\,\partial_{z_j}p)=2(2n+m-2)\partial_{z_j}p.
\]
Hence the right side of \eqref{eq:quadric-projection} belongs to
$I_{n+1}$ and its difference from $z_jp$ belongs to $Q_{n+1}$.
Fischer adjointness gives
\[
 \|q\,\partial_{z_j}p\|_{\rm F}^2
   =2(2n+m-2)\|\partial_{z_j}p\|_{\rm F}^2,
 \qquad
 \sum_j\|\partial_{z_j}p\|_{\rm F}^2=n\|p\|_{\rm F}^2,
\]
where the second equality is the Euler identity paired with $p$.
Converting the output squared norm by $(n+1)!$ and the input squared
norm by $n!$, we obtain
\[
 \sum_j\|B_{j,n}p\|_{\DA}^2
   =\frac{2n}{(n+1)(2n+m-2)}\|p\|_{\DA}^2.
\]
By polarization we obtain the exact positive block
\begin{equation}\label{eq:quadric-defect}
 A_n=\frac{2n}{(n+1)(2n+m-2)}I_{Q_n}\qquad(n\geq1).
\end{equation}
At $n=0$ all degree-one polynomials survive and $A_0=0$. In particular,
$n\|A_n\|\to1$. This gives a nonmonomial test of the adjacent-degree
conversion and again shows sharpness of the degree-decay order.

We can also check the coefficient independently by the trace identity. Since
multiplication by the nonzero quadratic polynomial is injective,
\[
 h_n=\binom{n+m-1}{m-1}-\binom{n+m-3}{m-1}
 =\frac{(2n+m-2)(n+m-3)!}{n!(m-2)!},
\]
where the second binomial is zero if $n<2$. Thus for $n\geq0$,
\[
 \frac{\Tr A_n}{h_n}
 =\frac{n+m}{n+1}-\frac{h_{n+1}}{h_n}
 =\frac{2n}{(n+1)(2n+m-2)}.
\]
The Fischer calculation established scalarity, so the trace gives
exactly the coefficient in \eqref{eq:quadric-defect}.
\end{example}

\begin{example}[The $2\times3$ determinantal ideal]\label{ext:segre-defect}
Let $R=\C[z_{ij}:1\leq i\leq2,\ 1\leq j\leq3]$ and let $I$ be the
ideal generated by all $2\times2$ minors of the matrix $(z_{ij})$.
The substitution $z_{ij}\mapsto u_i v_j$ induces in each degree $n$
an isomorphism
\[
 (R/I)_n\simeq \C[u_1,u_2]_n\otimes\C[v_1,v_2,v_3]_n.
\]
Indeed, each monomial of bidegree $(n,n)$ in the two sets of variables
is the image of a product of $n$ pairs $z_{ij}$. Two such products have
the same image precisely when the multisets of row and column indices
agree. After ordering the row indices, their column assignments are
related by transpositions. A transposition of column indices between
two different rows and two different columns is precisely the relation
$z_{ij}z_{k\ell}=z_{i\ell}z_{kj}$ modulo a minor; transpositions
with equal row or column indices leave the product unchanged. Hence
all products with the same image agree modulo $I$. Distinct image
monomials are linearly independent, proving the isomorphism and
\[
 h_n=(n+1)\binom{n+2}{2}=\frac{(n+1)^2(n+2)}2.
\]
The tensor product action of $U(2)\times U(3)$ on the six coordinate
variables is unitary. Its induced action on the Drury--Arveson space
is unitary, since unitary changes of variables preserve the kernel
$(1-\langle z,w\rangle)^{-1}$. The substitution is equivariant
for the natural actions on the two factors, so $I$ is invariant;
thus the quotient and its orthogonal realization $Q_n$ are equivalent
representations. This bidegree space is irreducible: an invariant
subspace is a sum of weight spaces for the diagonal tori, which are
one-dimensional monomial spaces. From any one monomial, the operators
$u_i\partial_{u_k}$ and $v_j\partial_{v_\ell}$ in the complexified
Lie algebra generate every monomial of the same bidegree by successive
transfers of exponents. Consequently any nonzero invariant subspace
is the whole bidegree space.

The operator $A_n$ is invariant under this action. Indeed the
orthogonal projections onto $I$ and $Q$ commute with it, and a unitary
linear change of the six coordinates rotates the column of discarded
coordinate maps by a unitary matrix on the column index. The sum of
their adjoint squares is unchanged. Schur's lemma now shows that
$A_n$ is scalar. Since $m=6$, the trace identity
\eqref{eq:op-trace} gives for every $n\geq0$
\begin{equation}\label{eq:segre-defect}
 A_n=\left(\frac{n+6}{n+1}
       -\frac{(n+2)(n+3)}{(n+1)^2}\right)I_{Q_n}
     =\frac{2n}{(n+1)^2}I_{Q_n}.
\end{equation}
Here $A_0=0$ and $n\|A_n\|\to2$. This calculation tests the rate for
a nonprincipal, nonmonomial ideal with nontrivial syzygies; it uses
no projective lifting estimate.
\end{example}

We impose no smoothness or transversality assumption on
the projective support or its link in the sphere. We deduce the analytic estimates from the homogeneous presentation rather than impose them as additional hypotheses on that support.  The use of the standard
grading is nevertheless essential to the argument.  A
nonhomogeneous relation such as $z_2-z_1^2$ need not contain its
homogeneous components in the ideal.  Such an ideal does not give
the projective ideal sheaf with the exact degree decomposition used
in \eqref{eq:geo-sections}; we make no conclusion about that general situation.

\hspace{0.9em}

\noindent{\bf Acknowledgments.} We thank Ronald G. Douglas, Kunyu Guo, Kai Wang, Yi Wang, and Guoliang Yu for inspiring discussions. Ding's research was supported in part by the National Natural Science Foundation of China (No.~12471127) and the Natural Science Foundation of Henan Province (No.~252300421208); Fang's research was supported in part by the National Science and Technology Council (No.~114-2115-M-A49-003-MY3); Tang's research was supported in part by DMS-2350181 and Simons Foundation Collaboration Grant MPS-TSM-00007714.\\  

\noindent\textbf{AI Declaration.} We used AI to assist with the discovery and development of the main proofs, the exploration and checking of arguments, and the preparation and revision of the exposition. We reviewed and rewrote the resulting material and take full responsibility for the mathematical claims, references, and final manuscript. Fang and Tang used ChatGPT-5.6, 6.0-6.2 (OpenAI), and Ding used DeepSeek V4, V4.1.

\raggedbottom

\end{document}